\documentclass[11pt]{article}
\usepackage[pagebackref,colorlinks=true,urlcolor=blue,linkcolor=red,citecolor=red]{hyperref}

\usepackage{epsfig, graphics, graphicx}
\usepackage{subcaption, comment, bm}
\usepackage[percent]{overpic}
\usepackage{float}
\usepackage{amssymb,bm}
\usepackage{amsthm}
\usepackage{color}
\usepackage{amsmath}
\usepackage[]{amsfonts}
\usepackage{fancyhdr}
\usepackage[]{graphicx,wrapfig}
\usepackage{enumitem}
\usepackage{array}
\usepackage{graphicx}
\usepackage[dvipsnames]{xcolor}
\newcommand\norm[1]{\left\lVert#1\right\rVert}

\newcommand{\D}{\mathrm{d}}

\newcommand{\Lc}{\mathcal{L}}
\newcommand{\Mc}{\mathcal{M}}

\newcommand{\Qc}{\mathcal{Q}}
\newcommand{\Rc}{\mathcal{R}}
\newcommand{\Sc}{\mathcal{S}}
\newcommand{\Tc}{\mathcal{T}}

\newcommand{\Vc}{\mathcal{V}}

\newcommand{\Xc}{\mathcal{X}}

\newcommand{\Zc}{\mathcal{Z}}

\newcommand{\Rb}{\mathbb{R}}
\newcommand{\Sb}{\mathbb{S}}

\newcommand{\vev}{\textbf{\textit{e}}}
\newcommand{\vf}{\textbf{\textit{f}}}
\newcommand{\vg}{\textbf{\textit{g}}}
\newcommand{\vx}{\textbf{\textit{x}}}

\newcommand{\vu}{\textbf{\textit{u}}}
\newcommand{\vv}{\textbf{\textit{v}}}

\newcommand{\vgamma}{\pmb{\gamma}}

\newcommand{\vxi}{\pmb{\xi}}

\usepackage{amsthm}
\usepackage{amsfonts}
\newcommand{\Beq}{\begin{equation}}
\newcommand{\Eeq}{\end{equation}}
\newcommand{\beq}{\begin{equation*}}
\newcommand{\eeq}{\end{equation*}}
\newcommand{\bal}{\begin{align}}
\newcommand{\eal}{\end{align}}

\graphicspath{{./Figs/}}

\newtheorem{thr}{Theorem}
\newtheorem{defn}{Definition}
\newtheorem{rem}{Remark}

\theoremstyle{definition}
\newtheorem{definition}{Definition}[section]

\newtheorem{remark}[definition]{Remark}

\title{\vspace{-1cm} V-line tensor tomography: numerical results}
\author{Gaik Ambartsoumian\thanks{Department of Mathematics, University of Texas at Arlington, Arlington, TX, United States of America.  \url{gambarts@uta.edu}}\and Rohit Kumar Mishra\thanks{Department of Mathematics, Indian Institute of Technology, Gandhinagar, Gujarat, India. \url{rohit.m@iitgn.ac.in}, \url{rohittifr2011@gmail.com}}
\and Indrani Zamindar\thanks{Department of Mathematics, Indian Institute of Technology, Gandhinagar, Gujarat, India. \url{indranizamindar@iitgn.ac.in}}}

\begin{document}
\date{}
\maketitle
\begin{abstract}
This article presents the numerical verification and validation of several inversion algorithms for V-line transforms (VLTs) acting on symmetric 2-tensor fields in the plane. The analysis of these transforms and the theoretical foundation of their inversion methods were studied in a recent work \cite{Gaik_Indrani_Rohit_24}. We demonstrate the efficient recovery of an unknown symmetric 2-tensor field from various combinations of the longitudinal, transverse, and mixed VLTs, their corresponding first moments, and the star VLT. The paper examines the performance of the proposed algorithms in different settings and illustrates the results with numerical simulations on smooth and non-smooth phantoms.
\end{abstract}
\vspace{2mm}

\vspace{-2mm}
\section{Introduction}\label{Introduction}
 Various generalizations of the Radon transform, 
 involving integrals over piecewise linear trajectories, have been considered in recent years by different research groups. Some typical examples include operators mapping a scalar function to its integrals along broken rays/lines (also called ``V-lines'') or over stars (finite unions of rays emanating from a common vertex). 
 Such integral transforms appear in single scattering optical tomography \cite{FMS-PhysRev-10, FMS-09}, single scattering X-ray tomography \cite{krylov2015inversion, walker2021iterative}, fluorescence imaging \cite{florescu2018}, Compton camera imaging \cite{Basko_et_al-97, Basko_V-proj}, Compton scattering emission tomography with collimated receivers \cite{morvidone2010, Truong_2011-Vline}, and gamma ray transmission/emission imaging  \cite{Rigaud2013}. 
 A detailed discussion of the history and the state of the art in studies of these generalized Radon transforms, as well as their applications in various imaging modalities can be found in the recent book \cite{amb_2023_book}. 

From the mathematical point of view, the broken ray/V-line transforms and their generalizations 
can be split into two distinct groups: those with the vertices of integration submanifolds inside the support of the integrand, and those where the vertices are outside (or on the boundary) of the support. 
This paper deals with transforms of the first type, and we briefly review below the relevant results on the subject. Several VLTs defined through a rotation invariant family of V-lines were studied in 
\cite{amb_2012, Gaik_Moon, Amb_Roy, Sherson}, where the authors came up with various inversion formulas of the transforms, as well as implemented and analyzed the performance of the resulting numerical algorithms. VLTs using translation invariant families of V-lines were studied in \cite{amb-lat_2019, Florescu-Markel-Schotland_2011, gouia2014exact, Kats_Krylov-13, krylov2015inversion, Sherson}. The obtained results included different inversion formulas, their numerical verification and validation, range description of the transforms, and support theorems.
VLTs arising in imaging models using curvilinear detectors (corresponding to V-lines with focused rays) have been analyzed in \cite{Kats_Krylov-13, krylov2015inversion, Sherson}. 
Microlocal analysis of VLTs has been discussed in \cite{amb-chapter, Sherson}. The study of various properties and the inversion of the star transform has been conducted in \cite{Amb_Lat_star, ZSM-star-14}. Several extensions of the VLT to  three and higher dimensions (the conical Radon transforms) have been studied in \cite{amb-lat_2019, gouia2014analytical, gouia2014exact, Palamodov2017}. The V-line transform has also been generalized to more general objects, such as vector and tensor fields in $\Rb^2$. In \cite{Gaik_Mohammad_Rohit_2020}, the authors introduced a new set of generalized V-line transforms (longitudinal, transverse, and their first integral moments), studied their properties, and derived various inversion algorithms to recover a vector field in $\Rb^2$ from different combinations of those transforms. The numerical verification and validation of the inversion methods proposed in \cite{Gaik_Mohammad_Rohit_2020} were demonstrated in the follow-up article \cite{Gaik_Mohammad_Rohit_2024_numerics}. Motivated by \cite{amb_2012, Gaik_Mohammad_Rohit_2020, Gaik_Moon}, the authors of \cite{Rohit_Rahul_Manmohan} studied rotationally invariant V-line transforms of vector fields and came up with two approaches for recovering the unknown vector field supported in a disk.

Another class of integral transforms considered in the literature is concerned with broken rays reflecting from the boundary of an obstacle (sometimes known as a reflector). In \cite{Ilmavirta_Salo_2016}, the authors studied this problem for scalar fields on a Riemannian surface in the presence of a strictly convex obstacle. Later, that work was extended to symmetric $m$-tensor fields in a similar setting in \cite{ilmavirta2022broken}. The broken ray transform of two tensors arises from the linearization of the length function of broken rays \cite{Ilmavirta_Salo_2016} and can be used in the field of seismic imaging. Notably, in this setting the authors showed that the kernel of the longitudinal transform consists of all potential tensor fields, which is aligned with the theory of straight-line transforms. In a recent work \cite{Shubham_Manas}, the authors considered a similar setup and proved uniqueness results for the broken ray transform acting on a combination of functions and vector fields on smooth surfaces. 

As a natural extension of works \cite{Gaik_Mohammad_Rohit_2020, Gaik_Mohammad_Rohit_2024_numerics}, we studied the V-line transforms of symmetric 2-tensor fields in $\Rb^2$ and discovered multiple interesting results about those operators  
\cite{Gaik_Indrani_Rohit_24}. 

In particular, it was shown that
the kernel descriptions of the longitudinal and transverse V-line transforms are very different from their counterparts appearing in the theory of straight-line transforms. We obtained exact inversion formulas to reconstruct a symmetric 2-tensor field from various combinations of longitudinal, transverse, and mixed V-line transforms, and their first integral moments. We also derived an inversion formula for the star transform of symmetric 2-tensor fields. In this paper, our aim is to verify and validate all inversion algorithms arising from the theoretical discoveries of \cite{Gaik_Indrani_Rohit_24}, and analyze their performance on various phantoms. All our simulations are done in MATLAB. To show the robustness of our numerical algorithms, we have done the reconstructions in the presence of various levels of noise.

The rest of the article is organized as follows. In Section \ref{sec: definitions and notations}, we introduce the integral transforms of interest and the required notations that we use throughout this article. At the end of that section we present two tables (Tables \ref{Table1}, \ref{Table2}), providing a summary of the algorithms studied in the paper and shortcuts to the appropriate sections discussing the reconstructions from particular combinations of the transforms. Section \ref{sec: data formation and numerical methods} starts with a description of the phantoms used in numerical simulations and explains the details of generating the forward data. 
A discussion about numerical solutions of certain partial differential equations (PDEs) required in some of our algorithms is presented at the end of that section. All numerical simulations and image reconstructions are demonstrated in Section \ref{sec: Numerical implementation}. This section is divided into subsections, each of which focuses on a specific combination of the transforms used for the reconstruction of the tensor fields (for details, see Tables \ref{Table1}, \ref{Table2}). We conclude the article with acknowledgments in Section \ref{sec:acknowledge}.

\section{Definitions and notations}\label{sec: definitions and notations}
This section is devoted to introducing the notations and definitions used throughout the article. The bold font letters are used to denote vectors and 2-tensors in $\mathbb{R}^2$ (e.g., $\vx$, $\textbf{\textit{u}}$, $\textbf{\textit{v}}$, $\vf$, etc.), and the regular font letters are used to denote scalars (e.g. $x_i$, $h$, $f_{ij}$, etc). The space of symmetric $2$-tensor fields defined in some disc $D\subset \Rb^2$ is denoted by $S^2(D)$, and the space of twice differentiable, compactly supported tensor fields is denoted by $C_c^2\left(S^2;D \right)$. The inner product in $S^2(D)$ is given by
\begin{align}\label{eq:inner product in Sm}
\langle \vf, \textbf{\textit{g}}\rangle =  \sum_{i, j = 1}^2 f_{ij}g_{ij} = f_{11}g_{11} + 2 f_{12}g_{12} + f_{22}g_{22}.
\end{align}

\noindent Next, we recall various well-known differential operators on scalar functions, vector fields, and tensor fields, which are needed in the upcoming discussions. For a scalar function $\varphi(x_1, x_2)$ and a vector field $\vf =(f_1,f_2)$, we use the notations 
\begin{align}\label{eq: definition of div and curl}
\D \varphi = \left(\frac{\partial \varphi}{\partial x_1}, \frac{\partial \varphi}{\partial x_2}\right), \ \  \D^\perp\varphi = \left(-\frac{\partial \varphi}{\partial x_2}, \frac{\partial \varphi}{\partial x_1}\right), \ \    \delta \vf =  \frac{\partial f_1}{\partial x_1}+ \frac{\partial f_2}{\partial x_2},\ \  
\delta^\perp \vf =  \frac{\partial f_2}{\partial x_1}- \frac{\partial f_1}{\partial x_2}.
\end{align}
These operators are naturally generalized to higher-order tensor fields in the following way: 
\begin{align}
(\D \vf)_{ij} &= \frac{1}{2}\left(\frac{\partial f_i}{\partial x_j} +\frac{\partial f_j}{\partial x_i}\right), \quad \vf \mbox{ is a vector field,}\label{eq:symm-der} \\
(\D^\perp \vf)_{ij} &= \frac{1}{2}\left((-1)^j\frac{\partial f_i}{\partial x_{3-j}} + (-1)^i \frac{\partial f_j}{\partial x_{3-i}}\right), \quad \vf \mbox{ is a vector field,}\\
(\delta \vf)_i &= \frac{\partial f_{i1}}{\partial x_1} + \frac{\partial f_{i2}}{\partial x_2} =  \frac{\partial f_{ij}}{\partial x_j}, \quad \vf \mbox{ is a symmetric 2-tensor field,}\\
(\delta^\perp \vf)_i &= - \frac{\partial f_{i1}}{\partial x_2} + \frac{\partial f_{i2}}{\partial x_1} =  (-1)^j\frac{\partial f_{ij}}{\partial x_{3-j}}, \quad \vf \mbox{ is a symmetric 2-tensor field.} \label{eq:def of deltaprep}
\end{align}
\noindent The directional derivative of a function $h$ in the direction $\vu\in \Sb^1$ is denoted by $D_{\vu}$, i.e. 
\begin{equation}
D_{\vu} h = \vu \cdot \D h.
\end{equation}
\begin{rem}
For a detailed discussion of the operators $\D$, $\D^\perp$, $\delta$, and $\delta^\perp$, we refer to \cite{derevtsov3,Sharafutdinov_Book}. 
\end{rem}
Consider a pair of fixed linearly independent unit vectors $\vu$ and $\vv$. We denote the rays emanating from $\vx \in \Rb^2$ in the directions $\vu$ and $\vv$ by
$$ L_{\vu}(\textbf{\textit{x}}) = \left\{\textbf{\textit{x}} +t \textbf{\textit{u}}: 0 \leq t < \infty\right\} \quad \mbox{ and } \quad  L_{\vv}(\textbf{\textit{x}}) = \left\{\textbf{\textit{x}} +t \textbf{\textit{v}}: 0 \leq t < \infty\right\}.$$
A V-line with the vertex $\textbf{\textit{x}}$ is the union of rays $L_{\vu}(\textbf{\textit{x}})$ and $L_{\vv}(\textbf{\textit{x}})$. Since the ray directions are fixed, each V-line can be uniquely identified by the coordinates of its vertex $\vx$.

Now, we are ready to recall the definitions of the transforms of our interest. These operators were introduced in \cite{Gaik_Indrani_Rohit_24}, where the authors presented various inversion algorithms to recover a symmetric $2$-tensor field from different combinations of the considered  transforms. 
\begin{defn}\label{eq: def of the divergent beam and its first moments}
\begin{enumerate}
    \item[(a)]  The \textbf{divergent beam transform} $\mathcal{X}_{\vu}$ of a function $h$ at $\vx\in \mathbb{R}^2$ in the direction $\vu\in\Sb^1$ is defined as:
 \begin{equation}\label{def:DivBeam}
   \mathcal{X}_{\vu}h(\vx) =  \int_{0}^{\infty} h(\vx+t \vu)\,dt.
 \end{equation}
\item[(b)] The \textit{\textbf{$1^{st}$ moment divergent beam transform}} of a function $h$ in the direction $\vu\in\Sb^1$ is defined as follows 
 \begin{equation}\label{def:moment DivBeam}
  \Xc^1_{\vu} h =  \int_0^\infty h(\vx + t \vu)\, t\,  dt.
 \end{equation}
 \item[(c)] The \textbf{\textit{V-line transform}} of a function $h$ with branches in the directions $\vu, \ \vv \in\Sb^1$ is defined as follows 
 \begin{equation}\label{def:V-line transform}
  \mathcal{V}h(\vx) = \Xc_{\vu} h(\vx) + \Xc_{\vv} h(\vx).
 \end{equation}
 \end{enumerate}  
\end{defn}

\begin{defn}
For  two vectors $\vu=(u_1, u_2)$ and $\vv = (v_1, v_2)$, the tensor product $\vu \otimes \vv$ is a rank-2 tensor with its $ij$-th component defined by
\begin{align}\label{eq:tensor product}
(\vu \otimes \vv)_{ij}   =  u_iv_j.
\end{align}
The symmetrized tensor product $\vu \odot \vv$
is then defined as
\begin{align}\label{eq:symmetrized tensor product}
\vu \odot \vv  =  \frac{1}{2}\left(\vu \otimes \vv + \vv \otimes \vu \right).
\end{align}
\end{defn}
\noindent We use the notation $\vu^2$ for the (symmetrized) tensor product of a vector $\vu$ with itself, that is, 
$$ \vu^2 = \vu \odot \vu = \vu \otimes \vu.$$
\begin{defn} \label{def:generalied V line transforms} Let $\vf\in C_c^2\left(S^2;\mathbb{R}^2\right)$.
\begin{enumerate}
    \item The \textbf{longitudinal V-line transform} of $\textbf{\textit{f}}$ is defined as
\begin{align}\label{eq:def longitudinal V-line transform}
\mathcal{L}_{\vu, \vv}\, \vf\  = \mathcal{X}_{\vu}  \left( \left\langle \vf , \vu^2 \right\rangle \right)  + \mathcal{X}_{\vv}  \left( \left\langle \vf , \vv^2 \right\rangle \right).
\end{align}
\item The \textbf{transverse V-line transform} of $\textbf{\textit{f}}$ is defined as
\begin{align}\label{eq:def transverse V-line transform}
\mathcal{T}_{\vu, \vv}\, \vf\  = \mathcal{X}_{\vu}  \left( \left\langle \vf , (\vu^\perp)^2 \right\rangle \right)  + \mathcal{X}_{\vv}  \left( \left\langle \vf , (\vv^\perp)^2 \right\rangle \right).
\end{align}
Here $\vu^\perp = (u_1, u_2)^\perp = (-u_2, u_1)$ is the the normal vector to $\vu$.
\item The \textbf{mixed V-line transform} of $\textbf{\textit{f}}$ is defined as
\begin{align}\label{eq:def mixed V-line transform}
\mathcal{M}_{\vu, \vv}\, \vf\  = \mathcal{X}_{\vu}  \left( \left\langle \vf , \vu \odot \vu^\perp \right\rangle \right)  + \mathcal{X}_{\vv}  \left( \left\langle \vf , \vv \odot \vv^\perp \right\rangle \right).
\end{align}
\item The \textbf{$1^{st}$ moment longitudinal V-line transform} of $\textbf{\textit{f}}$ is defined as
\begin{align}\label{eq:1-th moment longitudinal V-line transform}
\mathcal{L}^1_{\vu, \vv}\, \vf\  = \mathcal{X}^1_{\vu}  \left( \left\langle \vf , \vu^2 \right\rangle \right)  + \mathcal{X}^1_{\vv}  \left( \left\langle \vf , \vv^2 \right\rangle \right).
\end{align}
\item The \textbf{$1^{st}$ moment transverse V-line transform} of $\textbf{\textit{f}}$ is defined as
\begin{align}\label{eq:1-th moment transverse V line transform}
\mathcal{T}^1_{\vu, \vv}\, \vf\  = \mathcal{X}^1_{\vu}  \left( \left\langle \vf , (\vu^\perp)^2 \right\rangle \right)  + \mathcal{X}^1_{\vv}  \left( \left\langle \vf , (\vv^\perp)^2 \right\rangle \right).
\end{align}
\item The \textbf{$1^{st}$ moment mixed V-line transform} of $\textbf{\textit{f}}$ is defined as
\begin{align}\label{eq:1-th moment mixed V line transform}
\mathcal{M}^1_{\vu, \vv}\, \vf\  = \mathcal{X}^1_{\vu}  \left( \left\langle \vf , \vu \odot \vu^\perp \right\rangle \right)  + \mathcal{X}^1_{\vv}  \left( \left\langle \vf ,\vv \odot \vv^\perp \right\rangle \right).
\end{align}
\end{enumerate}
\end{defn}
\begin{rem}\label{rem:choice of u and v}
To simplify various calculations, we assume $\vu = (u_1, u_2)$ and $\vv =  (-u_1, u_2)$, that is, the V-lines are symmetric with respect to the $y$-axis. This choice does not change the analysis of the general case, since the data obtained in one setup of $\vu$ and $\vv$ can be transformed into the data obtained for the other setup, and vice versa.
\end{rem}

Next, we introduce the star transform of symmetric 2-tensor fields.  
In the upcoming definition, we identify a symmetric 2-tensor $\vf = \left(f_{ij}\right)$ with the vector $\vf = \left(f_{11}, f_{12},  f_{22}\right)\in\Rb^3$. Similarly, the tensor products $\vu^2$ and $\vu \odot \vv$ are identified with the vectors $\left(u_1^2, u_1 u_2, u_2^2\right)$ and $\left(u_1v_1, \frac{1}{2}(u_1 v_2 + u_2 v_1), u_2v_2\right)$.  

\vspace{1mm}
\begin{defn}\label{def: star transform}
Let $\vf\in C_c^2\left(S^2;\Rb^2\right)$,
and let $\vgamma_{1},\vgamma_{2},\dots ,\vgamma_{m}$ be distinct unit vectors in $\mathbb{R}^{2}$. The \textbf{star transform} of $\vf$ is defined as
\begin{align}\label{eq: star transform}
\Sc\vf &=  \sum_{i=1}^{m} c_{i} {\Xc}_{{\vgamma}_{i}} 
\begin{bmatrix}
 \vf\cdot {{\vgamma}_{i}}^{2}\\
 \vf\cdot {{\vgamma}_{i} \odot {\vgamma}_{i}^{\perp} }\\
 \vf\cdot ({\vgamma}_{i}^{\perp})^{2}
\end{bmatrix} ,
\end{align} 
where $c_{1}, c_{2},\dots , c_{m}$ are non-zero constants in $\mathbb{R}$.
\end{defn}

In our theoretical work \cite{Gaik_Indrani_Rohit_24}, we used various combinations of these transforms and derived multiple inversion formulas to recover a symmetric 2-tensor field in $\mathbb{R}^2$. All inversion methods and their numerical implementations are discussed in great detail in the upcoming sections, illustrating the effectiveness of the proposed algorithms. The tables below present a summary of these methods, identifying the transform combinations used in each specific method and the section discussing the corresponding reconstruction results. The first table is concerned with the recovery of special kinds of symmetric 2-tensor fields, while the second table deals with the recovery of general symmetric 2-tensor fields.\\
\vspace{1mm}
\begin{table}[hbt]
\centering
 \begin{tabular}{||c | l | c ||}   
 \hline
 Form of $\vf$ & Combinations of transforms used to recover $\vf$ & Sections and Figures\\ [0.5ex] 
 \hline\hline
$\D^2 \varphi$ & $\Lc \vf$ or $\Mc \vf$ (formula \eqref{eq:f=d^2varphi}) & \ref{subsec: d2phi}, Figure \ref{fig:f=d^2varphi}  \\ 
 \hline
$\D\D^\perp \varphi$ & $\Lc \vf$ or $\Tc \vf$ or $\Mc \vf$ (formulas \eqref{eq:f=dd^perp_varphi_LT} and \eqref{eq:f=dd^perp_varphi_M}) & \ref{subsec: d2phi}, Figures \ref{fig:f=dd^perp_varphi elliptic with noise}, \ref{fig:f=dd^perp_varphi parabolic with noise}, \ref{fig:f=dd^perp_varphi hyperbolic with noise}\\ 
 \hline

 $\D \vg$ PH1 & $\Lc \vf$ and $\Mc \vf$ (formulas \eqref{eq:g2}, \eqref{eq:g1_ellip}, and \eqref{eq:g1_para}) & \ref{subsec:dg and dperpg}, Figures \ref{fig:ph1 dg (elliptic)}, \ref{fig:ph1 dg (parabolic)}, \ref{fig:ph1 dg (hyperbolic)},  \ref{fig:ph1 dg with noise} \\ 
 \hline
$\D \vg$ PH2 &  $\Lc \vf$ and $\Mc \vf$ (formulas \eqref{eq:g2}, \eqref{eq:g1_ellip}, and \eqref{eq:g1_para}) & \ref{subsec:dg and dperpg},  Figures \ref{fig:ph2 dg (elliptic)}, \ref{fig:ph2 dg (parabolic)}, \ref{fig:ph2 dg (hyperbolic)},  \ref{fig:ph2 dg with noise} \\
   \hline
\end{tabular}
 \caption{Reconstructions of special kinds of symmetric 2-tensor fields.}
\label{Table1}
\end{table}
\begin{table}[hbt]
\centering
 \begin{tabular}{||l | c ||} 
 \hline
Combinations of transforms used to recover $\vf$ & Sections and Figures   \\ [0.5ex] 
 \hline\hline
$\mathcal{L}\vf$, $\mathcal{T}\vf$, and $\Mc \vf$  ($u_1 = u_2$) (formulas \eqref{eq:f11(u1=u2)}, \eqref{eq:f12(u1=u2)}, and \eqref{eq:f22(u1=u2)}) & \ref{subsec: full recovery Lf, Mf, and Tf}, Figures \ref{fig:ph2 (u1 = u2) Lf, Tf, and Mf}, \ref{fig:ph2 (u1 = u2) Lf, Tf, and Mf with noise}\\ 
 \hline
$\mathcal{L}\vf$, $\mathcal{T}\vf$, and $\Mc \vf$  ($u_1 \neq u_2$)  (formulas \eqref{eq:elliptic equation for f12}, \eqref{eq:f_11}, and \eqref{eq:f_22}) & \ref{subsec: full recovery Lf, Mf, and Tf}, Figures \ref{fig:ph1 (u1 neq u2) Lf, Tf, and Mf}, \ref{fig:ph1 (u1 neq u2) Lf, Tf, and Mf with noise}, \ref{fig:ph2 (u1 neq u2) Lf, Tf, and Mf}, \ref{fig:ph2 (u1 neq u2) Lf, Tf, and Mf with noise}\\ 
 \hline
$\mathcal{L}\vf$, $\Lc^1 \vf$, and $\mathcal{T}\vf$ ($u_1 \neq u_2$)  
(formulas  \eqref{eq: f11 from moments}, \eqref{eq: f12 from moments}, and \eqref{eq: f22 from moments})  & \ref{subsec: full recovery with moments}, Figures \ref{fig:ph1 from L, L1, T}, \ref{fig:ph1 from L, L1, T with noise}, \ref{fig:ph2 from L, L1, T}, \ref{fig:ph2 from L, L1, T with noise}\\ 
 \hline
$\mathcal{L}\vf$, $\Lc^1 \vf$, and $\Mc \vf$, (formulas \eqref{eq: f12 from moments L, L^1}, \eqref{eq: f11 from moments L, L^1}, and \eqref{eq: f22 from moments L, L^1}) & \ref{subsec: full recovery with moments}, Figures \ref{fig:ph1 from L, L1, M}, \ref{fig:ph1 from L, L1, M with noise}, \ref{fig:ph2 from L, L1, M}, \ref{fig:ph2 from L, L1, M with noise}\\
\hline
$\mathcal{S}\vf$ (formula \eqref{eq: star inversion})
& \ref{subsec: star transform}, Figures \ref{fig:ph2 from star transform}, \ref{fig:ph2 from star transform with noise} \\ 
 \hline
\end{tabular}
 \caption{Full recovery of a general symmetric 2-tensor field.}
\label{Table2}
\end{table}
\vspace{-5mm}

\section{Phantoms, data formation, and numerical solution of PDEs}\label{sec: data formation and numerical methods}
In this section, we introduce the phantoms used in our numerical experiments and delineate the generation of the forward data.  
We also discuss the methods used for numerical solution of PDEs appearing in some of the inversion algorithms developed in \cite{Gaik_Indrani_Rohit_24}.

\subsection{Description of phantoms}
The performance of all algorithms will be tested using two phantoms defined on $[-1,1] \times [-1,1]$. Since $\vf$ is a symmetric $2$-tensor field in $\mathbb{R}^2$, each phantom has 3 components $f_{11}$, $f_{12}$, and $f_{22}$. In the case of special tensor fields, we use one or two of these functions as scalar or vector potentials of the field.\\

\begin{figure}[H]
     \centering
   \includegraphics[width= 0.8\textwidth]{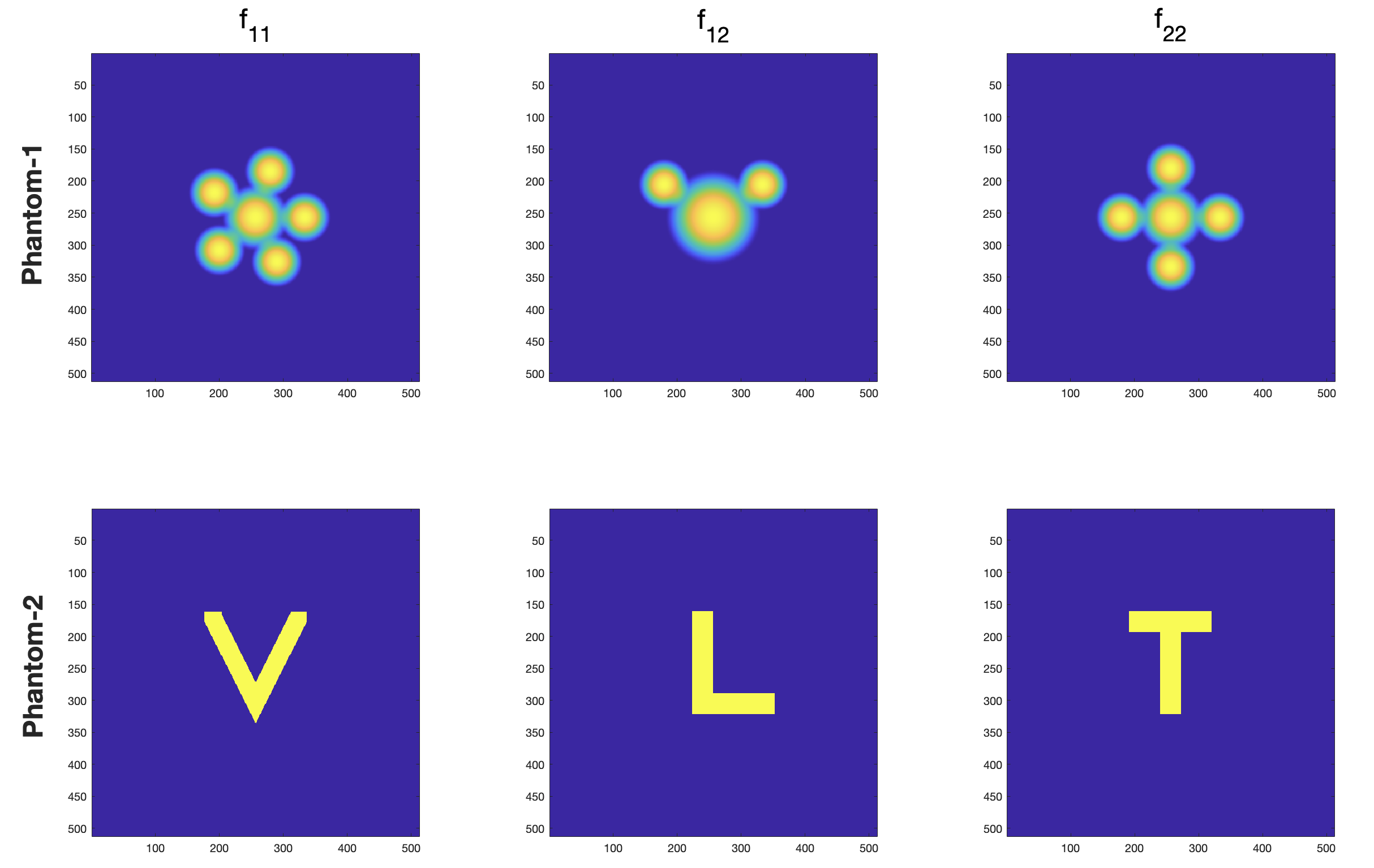}
    \caption{The components of the phantoms considered in numerical reconstructions.}\label{fig: phantoms}
 \end{figure}

\noindent \textbf{Phantom 1:} $\vf$ is a smooth field, and its components $f_{11}$, $f_{12}$, $f_{22}$ are different combinations of cutoff functions defined below (see Figure \ref{fig: phantoms}). Let $\displaystyle C^{(a,b)}_r $ be the smooth cut-off function, supported in a disk of radius $r$ and center at $(a, b)$:
\begin{align*}
C^{(a,b)}_r(x,y) = \begin{cases}
        e^{-r^2/\left(r^2-[(x-a)^2+(y-b)^2]\right)}, &(x-a)^2+(y-b)^2 < r^2,\\
     0,  & (x-a)^2+(y-b)^2 \ge r^2.
\end{cases}
\end{align*}
\begin{align*}
f_{11}(x,y)&= \left[C^{(0,0)}_{\sqrt{0.05}} 
  + C^{(0.09,0.28)}_{\sqrt{0.03}}  + C^{(-0.25,0.15)}_{\sqrt{0.03}}  +C^{(-0.22,-0.2)}_{\sqrt{0.03}}+C^{(0.13,-0.27)}_{\sqrt{0.03}}+C^{(0.3,0)}_{\sqrt{0.03}}\right] (x,y),\\
f_{12}(x,y) &= \left[C^{( 0,0)}_{\sqrt{0.1}}+ C^{(0.3, 0.2)}_{\sqrt{0.03}} + C^{(-0.3,0.2)}_{\sqrt{ 0.03}}\right](x,y), \\
    f_{22}(x,y) &= \left[C^{(0,0)}_{\sqrt{0.05}}(x,y)+ C^{(0,0.3)}_{\sqrt{0.03}}+C^{(0,-0.3)}_{\sqrt{ 0.03}} +C^{(-0.3,0)}_{\sqrt{0.03}}+C^{( 0.3,0)}_{\sqrt{0.03}}\right](x,y).
\end{align*}
            
\noindent \textbf{Phantom 2:} $\vf$ is non-smooth, with its components $f_{11}$, $f_{12}$, $f_{22}$ in the shape of the letters \textbf{``V'',} \textbf{``L''}, and \textbf{``T'', }respectively (see Figure \ref{fig: phantoms}). 
\begin{remark}
    It is well known that the numerical inversions of generalized Radon transforms work much better (produce less artifacts) on smooth phantoms than non-smooth ones. As it will be shown below, some of our reconstruction algorithms work extremely well even for the non-smooth phantoms. Therefore, in those cases, we choose to show the reconstructions only for Phantom 2 to save space and avoid repetition.  
\end{remark}

\subsection{Data formation}\label{data formation}
In the numerical simulations of V-line transforms, we use the unit vectors $\vu=(\cos\phi, \sin\phi)$ and $\vv = (-\cos \phi, \sin \phi)$,  with various choices of the angle $\phi = \pi/3$, $\pi/4$, or $\pi/6$. In the case of the star transform, we use the setup comprising three branches with polar angles $\phi_1 = 0$, $\phi_2=2\pi/3$, $\phi_3=4\pi/3$, and weights $c_i=1$ for $i=1,2,3.$

Recall, that all transforms introduced in Definitions \ref{def:generalied V line transforms} and \ref{def: star transform} are linear combinations of the divergent beam transform \eqref{def:DivBeam} and its first moment \eqref{def:moment DivBeam} of various projections of the unknown tensor field $\vf.$  Therefore, to generate the generalized V-line transforms and the star transform of $\vf$, we need a numerical method to compute the divergent beam transform and its first moment of a scalar field. The required algorithm has already been discussed in a recent work \cite[Section 3.2]{Gaik_Mohammad_Rohit_2024_numerics}, and we briefly present it here for the sake of completeness.

\vspace{2mm}
 
\noindent \textbf{Numerical evaluation of the divergent beam transform and its first moment:}
Let $F$ be an $n\times n$ pixelized image defined on $[-1,1] \times [-1,1]$. Then, the divergent beam transform of $F$ is also of the size $n\times n$, as the discretized set of rays are parametrized by the coordinates of their vertices and are assumed to emanate from the centers of the pixels. Let $\vx=(x,y)$ be the vertex, and $\vu=(\cos\phi, \sin\phi)$ be the direction of the ray. The first step of the computation is to find the intersections of the ray emanating from $\vx$ in the direction $\vu$ and the boundaries of square pixels appearing on the path of this ray. Next, we find the product of $F(i,j)$ and the length of the ray inside the pixel $(i,j)$. By summing up this product over all pixels, one generates the divergent beam transform of $F$, i.e. $\Xc_{\vu} F(\vx).$  

\vspace{2mm}
To evaluate the first moment of the divergent beam transform of $F$, we need to consider the product of three quantities $ F(i,j)$, the distance from the center of the pixel $(i,j)$ to the vertex $\vx$, and the length of the line segment of the ray inside the pixel $(i,j)$. Then, by adding this product over all pixels, we get the required  $\Xc^1_{\vu} F(\vx).$

\vspace{2mm}

As discussed above, to generate $\Lc\vf,\Tc\vf,$ and $\Mc\vf$ we evaluate the divergent beam transform  $\mathcal{X}_{\vu}$ of the projections  
$\left\langle \vf,\vu^2 \right\rangle,$ $\left\langle \vf,(\vu^\perp)^2 \right\rangle$, $\left\langle \vf,\vu \odot \vu^\perp \right\rangle$ and  $\mathcal{X}_{\vv}$ of $ \left\langle\vf,\vv^2 \right\rangle$, $\left\langle \vf,(\vv^\perp)^2 \right\rangle$, $\left\langle \vf,\vv \odot \vv^\perp \right\rangle$, and combine them according to equations \eqref{eq:def longitudinal V-line transform} \eqref{eq:def transverse V-line transform}, and \eqref{eq:def mixed V-line transform}. Similarly, $\Lc^1\vf,\Tc^1\vf,$ and $\Mc^1\vf,$ are generated by considering $\mathcal{X}^1_{\vu}$, $\mathcal{X}^1_{\vv}$ of the appropriate projections and combining them using formulas \eqref{eq:1-th moment longitudinal V-line transform}, \eqref{eq:1-th moment transverse V line transform}, and \eqref{eq:1-th moment mixed V line transform}. Finally, to generate $\Sc\vf$, we consider the divergent beam transform $\mathcal{X}_{\vgamma_{i}}$ of $\vf\cdot {{\vgamma}_{i}}^{2}$, 
 $\vf\cdot {{\vgamma}_{i} \odot {\vgamma}_{i}^{\perp} }$, 
 $\vf\cdot ({\vgamma}_{i}^{\perp})^{2}$ for $i= 1,2 \dots m$ and use formula \eqref{eq: star transform}  to combine them. 
 
 The graphical representations of $\Lc\vf,\Tc\vf,$ and $\Mc\vf$ corresponding to Phantoms 1 and Phantom 2 can be found in Figures \ref{fig:ph1 (u1 neq u2) Lf, Tf, and Mf} and \ref{fig:ph2 (u1 = u2) Lf, Tf, and Mf}, respectively. $\Sc\vf$ of Phantom 2 is depicted in Figure \ref{fig:ph2 from star transform}.

Another operator that is used repeatedly throughout the article is the directional derivative of a scalar function $h$ along $\vu$ or $\vv$. Numerically $D_{\vu}h$ is computed in the following two steps ($D_{\vv}h$ is computed in exactly the same way): 
\begin{itemize}
    \item Calculate the gradient $\left(\partial_{x} h,\partial_{y}h \right)$ by Matlab function \textbf{gradient}.
    \item Compute  $\displaystyle D_{\vu}h= u_1 \partial_{x} h+ u_2 \partial_{y}h.$
\end{itemize}

\begin{rem}\label{rem: different pixel size for experiments}
In some inversion formulas of \cite{Gaik_Indrani_Rohit_24}, one needs to solve initial/boundary value problems for PDEs, which requires an inversion of an $n^2\times n^2$ matrix. In all such experiments, we use images with a resolution of $160 \times 160$ pixels to reduce the computational time. For all other experiments, we use images with a resolution of $512 \times 512$ pixels.
\end{rem}

\subsection{Solving PDEs numerically}
As mentioned above in Remark \ref{rem: different pixel size for experiments}, in some cases \cite[Theorems 4, 5, and 6]{Gaik_Indrani_Rohit_24} the recovery of the unknown symmetric $2$-tensor field is achieved by solving initial/boundary value problems for PDEs. The equations appearing in \cite[Theorems 4, 5, and 6]{Gaik_Indrani_Rohit_24} are of the following form:
\begin{align}\label{eq:partial differential equation for u}
a u_{xx}+ b u_{yy} = -f, \quad \mbox{ in } \Omega=[-1,1]\times[-1,1].
\end{align}
This PDE can be elliptic, parabolic, or hyperbolic in nature, depending on the coefficients $a$ and $b$. We solve it numerically, with the appropriate initial/boundary conditions discussed below. 
\subsubsection*{Elliptic PDE ($a>0$ and $b>0$)}
\noindent

In this case we need to solve numerically the following Dirichlet boundary value problem.
\begin{align}\label{eq:elliptic equation for u}
    \left\{\begin{array}{rll}
   a u_{xx}+ b u_{yy}&= -f      & \mbox{ in } \Omega=[-1,1]\times[-1,1],  \\
     u &= g    & \mbox{ on } \partial \Omega.
    \end{array}\right.
\end{align}
Dividing $\Omega$ into $N \times N$ uniform pixels with the pixel size $h \times h,$ we use the central difference approximation to write the second-order derivatives at an interior grid point $(x_i, y_j )$ as:
\begin{align} \label{FDM for d^u/dx^2 and d^u/dy^2}
    \frac{\partial^2u}{\partial x^2}(x_i, y_j)= \frac{u_{i+1,j}-2u_{i,j}+u_{i-1,j}}{h^2}, \qquad  \frac{\partial^2u}{\partial y^2}(x_i, y_j)= \frac{u_{i,j+1}-2u_{i,j}+u_{i,j-1}}{h^2},
\end{align}
where $u_{i,j} = u(x_i, y_j ).$ Then the differential operator in equation \eqref{eq:elliptic equation for u} can be approximated at $(x_i, y_j)$ as follows:
$$\left[a u_{xx}+ b u_{yy}\right](x_i, y_j )= \frac{au_{i-1,j}-2(a+b)u_{i,j}+au_{i+1,j}+bu_{i,j-1}+bu_{i,j+1}}{h^2}.$$
Thus, the finite difference version of the partial differential equation at the interior points is given by
\begin{align}\label{FDM for the elliptic pde}
 -\frac{au_{i-1,j}-2(a+b)u_{i,j}+au_{i+1,j}+bu_{i,j-1}+bu_{i,j+1}}{h^2}=f_{i,j}.
\end{align}
By defining the index map $(i, j) \rightarrow k = (N-2)(i-2)+(j-1)$ for $2 \le i,j \le N-1$, we arrange  the interior $(N-2)\times (N-2)$ grid points in one row using a single index $k = 1$ to $(N-2)^2$ and use notations $u_k = u_{i(k),j(k)}$ and $f_k=f_{i(k),j(k)}$. With this choice of indexing, equation \eqref{FDM for the elliptic pde} takes the following form:
\begin{align}\label{FD matrix for elliptic pde}
  AU = F,
  \end{align}
where $U=(u_{k})_{1\le k\le (N-2)^2}$ and $A$ is an $(N-2)^2 \times (N-2)^2$ matrix that has block tridiagonal structure given by 
\begin{align*}
A= \begin{pmatrix}
    B & C & 0 & \dots & 0 & 0 & 0\\
    C & B & C & \dots & 0 & 0 & 0\\
   \vdots & \vdots & \vdots &\ddots & \vdots & \vdots & \vdots\\
   0 &0 &0 & \dots \quad   &  C & B & C\\
   0& 0 & 0 & \dots    & 0 & C & B
\end{pmatrix} \quad \mbox{ with } \quad C= -bI_{(N-2)\times (N-2)},
\end{align*}
\begin{align*}
B= \begin{pmatrix}
    2(a+b) & -a & 0 & \dots & 0 & 0 & 0\\
    -a & 2(a+b) & -a & \dots & 0 & 0 & 0\\
   \vdots & \vdots & \vdots &\ddots & \vdots & \vdots & \vdots\\
   0 &0 &0 & \dots \quad   &  -a & 2(a+b) & -a\\
   0& 0 & 0 & \dots    & 0 & -a & 2(a+b)
\end{pmatrix}_{(N-2)\times (N-2)}, 
\end{align*}
and $F = h^2(\Tilde{f_k})_{1\le k \le (N-2)^2}$. Here $\Tilde{f}$ represents the modified source term, which satisfies $ \Tilde{f}_{i,j} = f_{i,j} $ for $ 3 \le i, j \le N-2,$ and involves the boundary terms (i.e. the given data $g(i, j)$) for other indices. More specifically, we have
\begin{align*}
    \begin{array}{lrll}
 \Tilde{f}_{2,2} = f_{2,2} + \frac{1}{h^2}\left(ag_{1,2} + bg_{2,1}\right),  & &\Tilde{f}_{N-1,2} = f_{N-1,2} + \frac{1}{h^2}\left(ag_{N,2} + bg_{N-1,1}\right), \\
  \Tilde{f}_{2,j} = f_{2,j} + \frac{a}{h^2}g_{1,j}, ~~ \text{for} ~3 \le j \le N-2,&       &  \Tilde{f}_{N-1,j} = f_{N-1,j} + \frac{a}{h^2}g_{N,j}, ~\text{for} ~ 3 \le j \le N-2,\\
 \Tilde{f}_{i,2} = f_{i,2} + \frac{b}{h^2}g_{i,1},~\text{for}~ 3 \le j \le N-2, & & \Tilde{f}_{i,N-1} = f_{i,N-1} + \frac{b}{h^2}g_{i,N}, ~~\text{for} ~ 3 \le j \le N-2,\\
 \Tilde{f}_{2,N-1} = f_{2,N-1} + \frac{1}{h^2}\left(a g_{1,N-1} + b g_{2,N}\right), & & \Tilde{f}_{N-1,N-1} = f_{N-1,N-1} + 
\frac{1}{h^2}\left(a g_{N,N-1} + b g_{N-1,N}\right).
    \end{array}
\end{align*}
Finally, we solve the system of linear equations $AU = F$ to get $U$ as a numerical solution $u$ of the required boundary value problem \eqref{eq:elliptic equation for u} in the case when $a>0$ and $b>0$. 

\subsubsection*{Parabolic PDE ($a = 0$ or $b = 0$)}
For $b=0,$ the PDE \eqref{eq:partial differential equation for u} reduces to 
\begin{align} \label{eq:hyperbolic pde}
    u_{xx}=-f/a.
\end{align}
This equation is solved by repeated integration along $\vev_{1} = (1, 0)$. 

More specifically, we apply the divergent beam transform $\Xc_{\vev_{1}}$ twice to equation \eqref{eq:hyperbolic pde} to obtain $u$.

The case $a =0$ is treated analogously, by integrating twice in the direction  ${\vev}_{2}=(0,1)$  to get $u$.
 
\subsubsection*{Hyperbolic PDE ($a$ and $b$ are of opposite signs) }
Let us consider the case when $a>0$ and $b<0$ (the other case is analogous). 
For $b=-\Tilde{b},$ where $\Tilde{b}>0$, we get the following initial value problem from \eqref{eq:partial differential equation for u} (see the discussion on \cite[Page 10]{Gaik_Indrani_Rohit_24}): 
\begin{align}\label{eq:hyperbolic equation for u}
    \left\{\begin{array}{rll}
   a u_{xx}- \Tilde{b} u_{yy}& = -f      & \mbox{ in } \Omega=[-1,1]\times[-1,1],  \\
     u(-1,y) &= g,\\
     u_x(-1,y) &=\Tilde{g} .
    \end{array}\right.
\end{align}
Dividing $\Omega$ into $N \times M$ uniform pixels with the pixel size $h \times k$, we use the forward difference approximation for the first-order derivative at the left boundary and the central difference approximation for the second-order derivatives at an interior grid point $(x_i, y_j )$. 
\begin{align*}
    \frac{\partial u}{\partial x}(x_i, y_j)&= \frac{u_{i+1,j}- u_{i,j}}{h}, \qquad \qquad \qquad
    \frac{\partial^2u}{\partial x^2}(x_i, y_j)= \frac{u_{i+1,j}-2u_{i,j}+u_{i-1,j}}{h^2}, \\
    \frac{\partial^2u}{\partial y^2}(x_i, y_j)& = \frac{u_{i,j+1}-2u_{i,j}+u_{i,j-1}}{k^2}, \qquad  \mbox{ where } u_{i,j} = u(x_i, y_j) .
\end{align*}
The finite difference version of the hyperbolic PDE is given by
\begin{align}\label{FDM for the hyperbolic pde}
\left\{\begin{array}{lll}
  u_{i+1,j} =2u_{i,j}-u_{i-1,j}+\lambda\left( u_{i,j+1}-2u_{i,j}+u_{i,j-1}\right) -\frac{h^2}{a}f_{i,j},    &  2\le i, j\le (N-1)\\
  u_{1,j}=g_{1,j}   & 1\le j\le N  \\
  u_{2,j}= u_{1,j} + h\Tilde{g}_{1,j}, & 1\le j\le N
\end{array}\right.
\end{align}
where $\displaystyle \lambda= \frac{\Tilde{b}}{a}\frac{h^2}{k^2}.$\\
\vspace{2mm}

\noindent Relation \eqref{FDM for the hyperbolic pde} provides a method for iterative solution of the initial value problem \eqref{eq:hyperbolic equation for u}.

\begin{rem}
    The stability of the finite difference method discussed above  
    (for the hyperbolic case) depends on the proper balance between the step sizes $h$ and $k$, given by the Courant-Friedrichs-Lewy (CFL) condition. The CFL condition ensures that the numerical domain of dependence contains the analytical domain of dependence for a given initial condition. In our setting, this condition can be written as:
    $$C\doteq\sqrt{\lambda}
    \le 1.$$
    The constant $C$ is often called the \textit{Courant number}. For a more detailed discussion on the CFL condition, please refer to \cite[Section 2.2.3]{durran2013numerical}
\end{rem}
\begin{rem}
Note that all integral transforms (the forward data) are numerically computed on grids with $h = k$. However, for a stable reconstruction using hyperbolic PDEs, in some cases we need to take $h \neq k$. To address this issue we use interpolation to generate the data on a finer grid of size $4 N\times N $, that is $k = 4h$, and use this data to solve the given hyperbolic equation in a stable way with the Courant number $C=1/ 4$.  
\end{rem}

\section{Numerical Implementation}\label{sec: Numerical implementation}
In this section, we present the results of the numerical implementation of the inversion algorithms proposed in \cite{Gaik_Indrani_Rohit_24} for the transforms introduced in Section \ref{sec: definitions and notations}. We divide our presentation into subsections, each of which deals with either a particular type of tensor field or a particular combination of transforms used for reconstruction. Each subsection starts with a brief review of the theoretical result implemented there. Below is a summary of the results discussed in this section.
\begin{enumerate}
\item \textbf{Special tensor fields (here $\varphi$ is a function, and $\vg$ is a vector field)}
\begin{enumerate}[label=(\roman*)]
\item The tensor field $\vf=\D^2\varphi$ is recovered from 
either $\Lc\vf$ or $\Mc\vf$.
\item The tensor field $\vf=(\D^\perp)^2\varphi$ is recovered from 
either $\Tc\vf$ or $\Mc\vf$.
\item The tensor field $\vf=\D\D^\perp\varphi$ is recovered from 
either $\Lc\vf$, or $\Tc\vf$, or $\Mc\vf$.
\item The tensor field $\vf=\D\vg$ is recovered from 
the combination of $\Lc\vf$ and $\Mc\vf$.
\item The tensor field $\vf=\D^{\perp}\vg$ is recovered from 
the combination of $\Tc\vf$ and $\Mc\vf$.
\end{enumerate}
\item $\vf$ is recovered from the combination of $\Lc\vf,\Tc\vf$ and $\Mc\vf.$
\item  When $u_1\ne u_2$, $\vf$ is recovered from any of the  following combinations 
\begin{enumerate}
\item $\Lc\vf,\Lc^1\vf$ and $\Tc\vf,$
\item $\Tc\vf,\Tc^1\vf$ and $\Lc\vf,$
\item $\Mc\vf,\Mc^1\vf$ and $\Lc\vf,$
\item $\Mc\vf,\Mc^1\vf$ and $\Tc\vf.$
\end{enumerate}
\item $\vf$ is recovered either from the  combination of $\Lc\vf,\Lc^1\vf,\Mc\vf$  or $\Tc\vf,\Tc^1\vf,\Mc\vf$.
\item $\vf$ is recovered from $\Sc\vf$.
\end{enumerate}

Our choice of the special tensor fields listed above is not arbitrary. Their forms are motivated by the decomposition $\vf=(\D^\perp)^2 \varphi + (\D\D^\perp)\chi + \D^2\psi$ of a symmetric 2-tensor field $\vf$ discussed in \cite[Theorem 4.2]{derevtsov3}. In the straight-line setup considered in \cite{derevtsov3},  the scalar functions $\varphi,\chi$, and $\psi$ are recovered explicitly from the knowledge of the longitudinal, mixed, and transverse ray transforms of $\vf$, respectively. The situation is different in the V-line setting, as discussed in Theorem \ref{th: special tensor fields 1} below. Namely, the function $\varphi$ (when $\vf=\D^2\varphi$) can be recovered explicitly from the longitudinal or mixed V-line transform of $\vf$, and similar results hold for the other types of special tensor fields.  

On the other hand, the special tensor fields $\vf = \D \vg$ and $\vf = \D^\perp \vg$ correspond to the decomposition $\vf =  \D^\perp \tilde{\vg} + \D \vg$, which can be obtained by adjusting the terms of the decomposition from \cite[Theorem 4.2]{derevtsov3} mentioned above. It is known that $\D\vg$ is in the kernel of the longitudinal (straight line) ray transform. However, $\D\vg$ is NOT in the kernel of the longitudinal V-line transform. Moreover, as we have shown in \cite{Gaik_Indrani_Rohit_24}, any potential tensor fields $\vf=\D\vg$ can be recovered from the combination of the longitudinal and mixed VLTs. Similarly, $\D^\perp \vg$ is in the kernel of the transverse (straight line) ray transform, but not in the kernel of the transverse V-line transform. Instead, one can recover $\vf=\D^\perp \vg$ from the combination of the transverse and mixed VLTs of $\vf$.

\subsection{Numerical implementation for special kinds of tensor fields}
\subsubsection{Tensor fields of the form \texorpdfstring{$\vf$}{f} = \texorpdfstring{$\D^2\varphi$}{d2 varphi}, \texorpdfstring{$\vf$}{f} = \texorpdfstring{$\D \D^\perp \varphi$}{ddperp varphi}, or \texorpdfstring{$\vf$}{f} = \texorpdfstring{$(\D^\perp)^{2} \varphi$}{(dperp2 varphi} }\label{subsec: d2phi}
 
\begin{figure}[H]
     \centering
   \includegraphics[width= 1\textwidth]{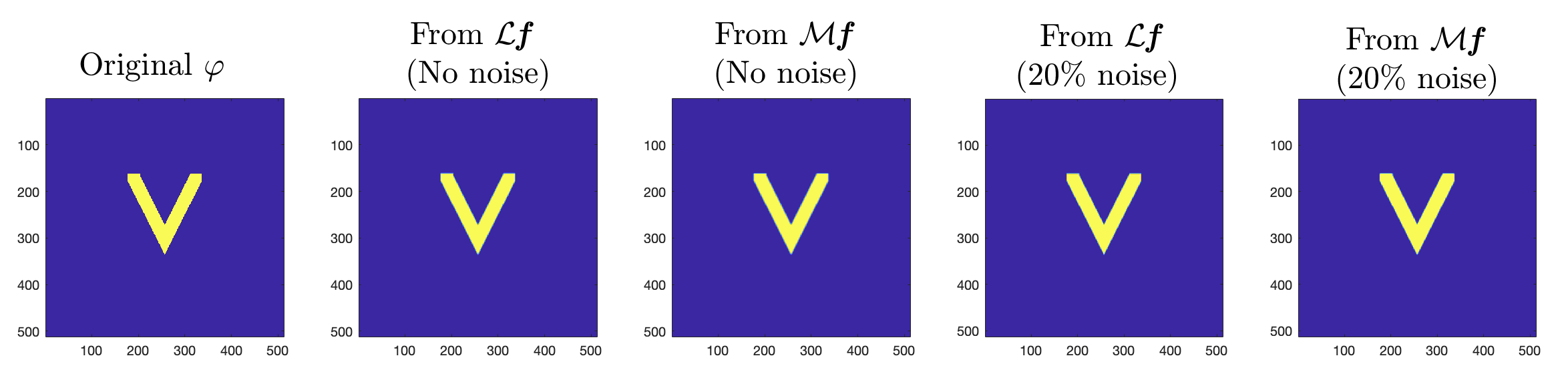}
    \caption{Recovery of $\varphi$ when $\vf$=$\D^2\varphi$, using formula \eqref{eq:f=d^2varphi}.}\label{fig:f=d^2varphi}
 \end{figure}

In this subsection, we consider symmetric 2-tensor fields of the form $\vf=\D^2\varphi, (\D^\perp)^2\varphi$, or $\D\D^\perp\varphi,$ where $\varphi$ is a scalar function. We reconstruct $\varphi$ from specific VLTs, depending on the form of $\vf$. 
Our reconstructions of $\varphi$ when $\vf=\D^2\varphi$ and when $\vf=(\D^\perp)^2\varphi$ are almost identical. Therefore, to save space and avoid repetition, we present below only the numerical results corresponding to the first case. We also demonstrate the numerical reconstruction of $\varphi$ when $\vf=\D\D^\perp\varphi$. To avert redundancy, we have not included reconstructions without noise in this case.   
\begin{thr}\label{th: special tensor fields 1} Let $\varphi$ be a twice differentiable, compactly supported function, that is, $\varphi \in C^2_c(D_1)$. 
\begin{itemize}
\item[(a)] If  $\vf$ is a symmetric 2-tensor field of the form $\vf = \D^2 \varphi$, then it can be reconstructed explicitly in terms of $\Lc \vf$ or $\Mc \vf$,
using the following formulas:  
\begin{align}\label{eq:f=d^2varphi}
 \varphi(\vx) = \frac{1}{2u_2}  \int_0^\infty \Lc \vf(\vx + s\vev_2)ds  = - \frac{1}{2u_2}  \int_0^\infty \Mc \vf(\vx + s\vev_1)ds .
\end{align} 
\item[(b)]  If  $\vf$ is a symmetric 2-tensor field of the form $\vf = \left(\D^\perp\right)^2 \varphi$, then it can be reconstructed explicitly in terms of $\Tc \vf$ or $\Mc \vf$,
using the following formulas: 
\begin{align}\label{eq:f=(d^perp)^2varphi}
 \varphi(\vx)= \frac{1}{2u_2}  \int_0^\infty \Tc \vf(\vx + s\vev_2)ds =  \frac{1}{2u_2}  \int_0^\infty \Mc \vf(\vx + s\vev_1)ds.  
\end{align}

\item[(c)] If  $\vf$ is a symmetric 2-tensor field of the form $\vf = \D\D^\perp \varphi$, then it can be reconstructed explicitly in terms of $\Lc \vf$ or $\Tc \vf$, using the following formulas: 
 \begin{align}\label{eq:f=dd^perp_varphi_LT}
  \varphi(\vx)= \frac{1}{2u_2}  \int_0^\infty \Lc \vf(\vx + s\vev_1)ds = -\frac{1}{2u_2}  \int_0^\infty \Tc \vf(\vx + s\vev_1)ds.
 \end{align}
In this case, $\vf$ can also be reconstructed from $ \Mc\vf$ by solving the following second order partial differential equation for $\varphi$ (with appropriate initial/boundary conditions):
 \begin{align}\label{eq:f=dd^perp_varphi_M}
        \left(1+2u_1^2\right)\frac{\partial^2 \varphi}{\partial x_1^2} + \left(u_1^2 - u_2^2\right) \frac{\partial^2 \varphi}{\partial x_2^2} &= -\frac{1}{u_2}\Xc_{\vev_2}\left(D_\vu D_\vv \Mc \vf \right)(\vx).
 \end{align}
\end{itemize}
\end{thr} 
\noindent \underline{Let $\vf=\D^2\varphi.$} Figure \ref{fig:f=d^2varphi} shows the results of the numerical implementation of formula \eqref{eq:f=d^2varphi} for a non-smooth phantom, with V-lines corresponding to $\vu=(\cos{\pi/4},   \sin{\pi/4})$.

Here and in all other numerical implementations presented in the paper, we test the effectiveness of the inversion algorithm by adding to the data different amounts of Gaussian noise. Figure \ref{fig:f=d^2varphi} includes reconstructions with $20\%$ of additive Gaussian noise. Also, here and in all upcoming sections, we include a table describing the relative error of the reconstructions,  calculated as follows:
\begin{align*}
    \text{Relative error of }\varphi= \frac{\norm{\varphi_{\text{original}} -\varphi_{\text{rec}}}_2}{\norm{\varphi_{\text{original}}}_2}\times 100\%,
\end{align*}
where $\| \varphi \|_2$ denotes the spectral norm, i.e. the maximum singular value of the matrix $\varphi$.

\begin{table}[h!]
 \begin{center}
\begin{tabular}{ |c|c|c|c| } 
 \hline
   From $\Lc\vf$(No noise) & From $\Mc\vf$(No noise)& From $\Lc\vf$($20\%$ noise)& From $\Mc\vf$($20\%$ noise)\\ 
 \hline
5.60\% &4.99\% & 10.91\% & 36.61\%\\ 
 \hline
\end{tabular}
\caption{Relative errors of reconstruction of $\varphi$ from $\Lc\vf$ or $\Mc\vf$ when $\vf= \D^2\varphi$.}\label{tab:V-d2phi}
\end{center}
\end{table}

Since visually all reconstructions in Figure \ref{fig:f=d^2varphi} look the same, we include below the images of the differences of the original image and the reconstructions, explaining the relative errors in Table \ref{tab:V-d2phi}. We skip this part for the other reconstructions in this paper.

\begin{figure}[H]
     \centering
   \includegraphics[width= 0.9\textwidth]{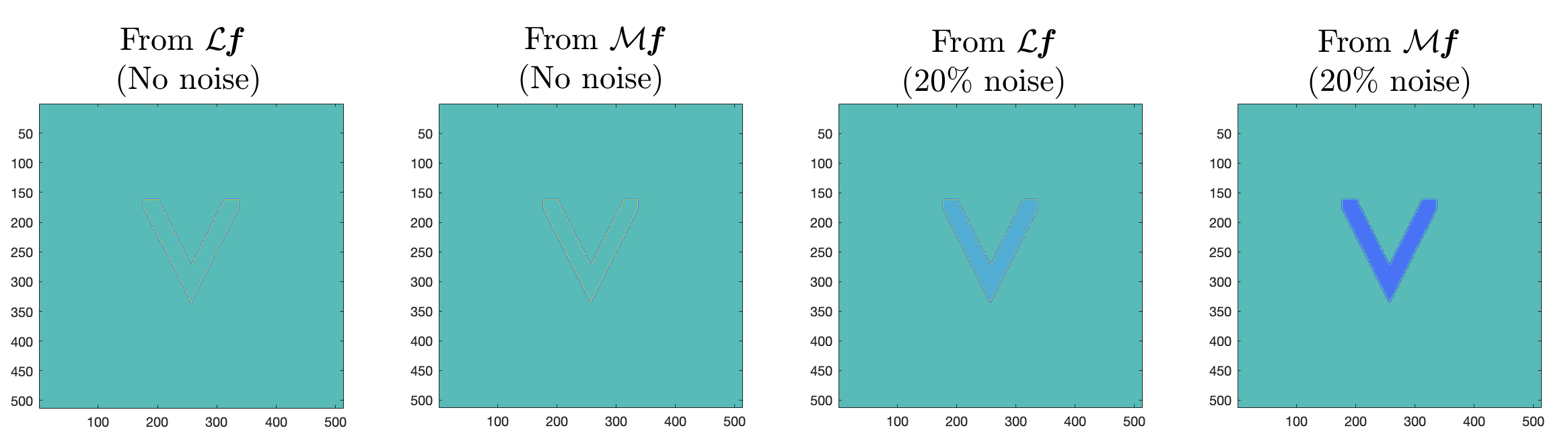}
    \caption{Graphical representation of $\varphi_{\text{original}} -\varphi_{\text{rec}}$. The colormap is scaled according to the minimum and maximum values of all four pictures.}\label{fig:err} 
 \end{figure}

\begin{figure}[H]
     \centering
   \includegraphics[width= 0.84\textwidth]{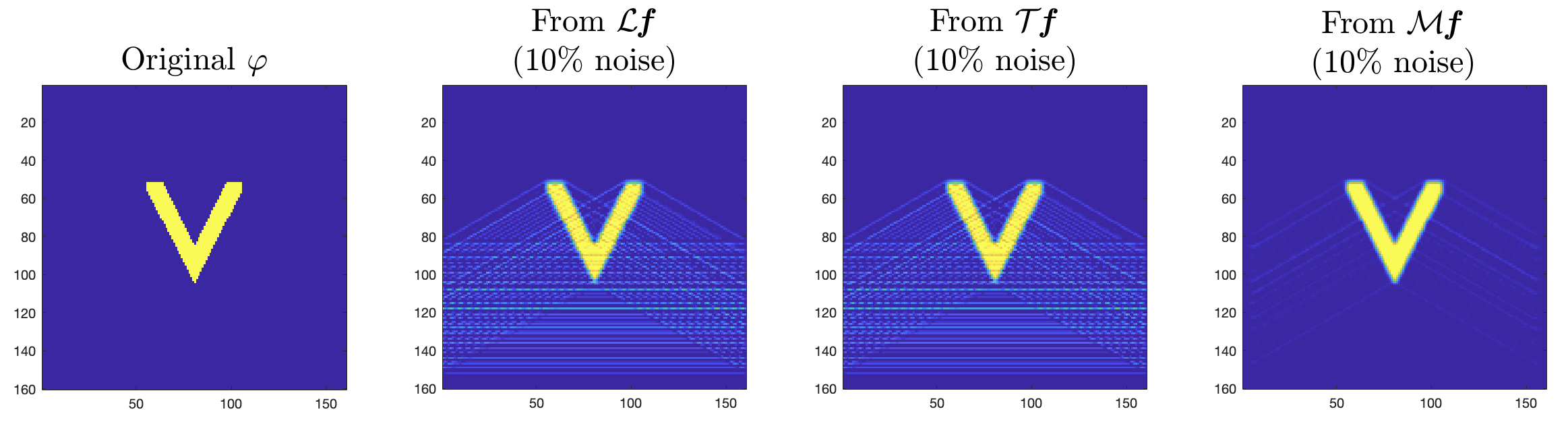}
    \caption{Recovery of $\varphi$ when $\vf$=$\D\D^\perp\varphi$, $u_1 >u_2$ (elliptic), with $10\%$ noise, using \eqref{eq:f=dd^perp_varphi_LT}, \eqref{eq:f=dd^perp_varphi_M}.}\label{fig:f=dd^perp_varphi elliptic with noise}
 \end{figure}

\begin{figure}[H]
     \centering
   \includegraphics[width= 0.9\textwidth]{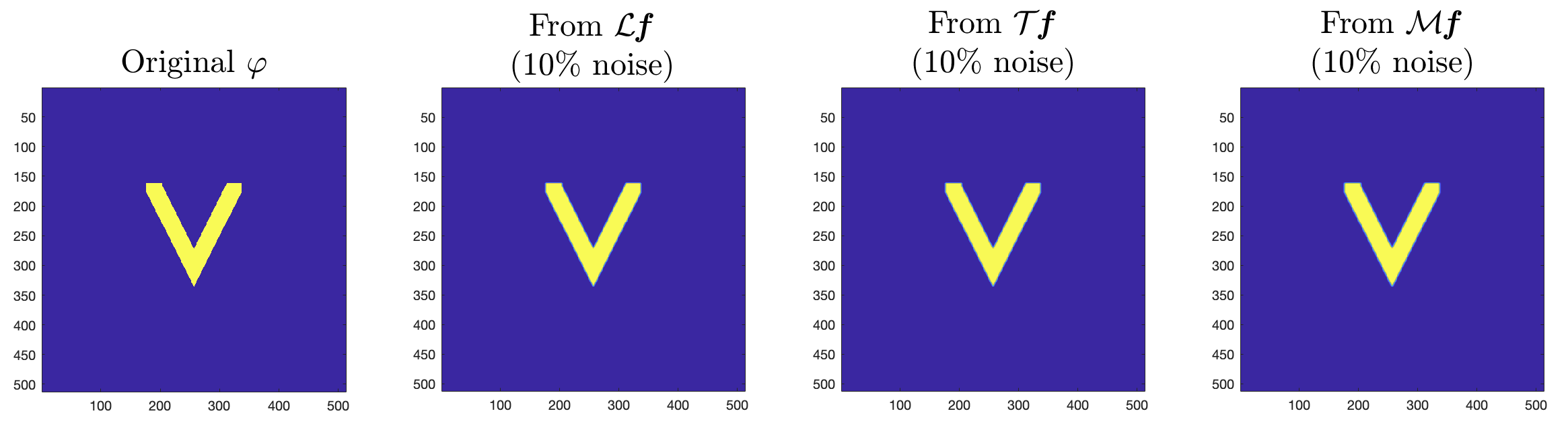}
    \caption{Recovery of $\varphi$ when $\vf$=$\D\D^\perp\varphi$, $u_1=u_2$ (parabolic), with $10\%$ noise,  using  \eqref{eq:f=dd^perp_varphi_LT}, \eqref{eq:f=dd^perp_varphi_M}.} \label{fig:f=dd^perp_varphi parabolic with noise}
 \end{figure}

\begin{figure}[H]
     \centering
   \includegraphics[width= 0.9\textwidth]{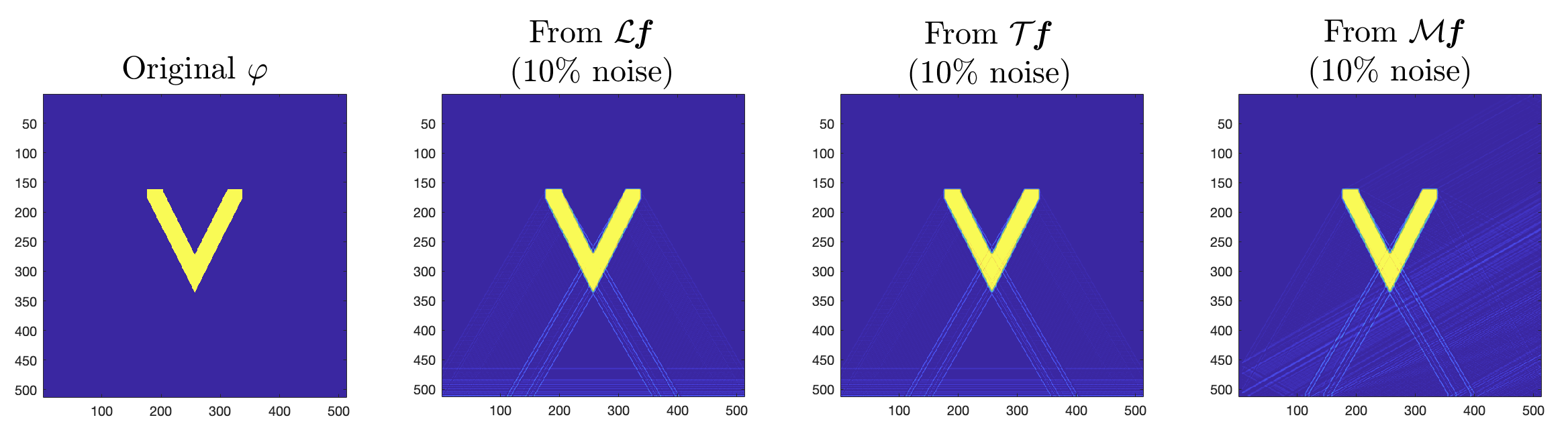}
    \caption{Recovery of $\varphi$ when $\vf$=$\D\D^\perp\varphi$, $u_1 < u_2$ (hyperbolic), with $10\%$ noise, using \eqref{eq:f=dd^perp_varphi_LT}, \eqref{eq:f=dd^perp_varphi_M}.}\label{fig:f=dd^perp_varphi hyperbolic with noise}
 \end{figure}

\noindent \underline{Let $\vf=\D\D^\perp \varphi.$} The reconstructions using formulas \eqref{eq:f=dd^perp_varphi_LT} and \eqref{eq:f=dd^perp_varphi_M} are shown in Figures \ref{fig:f=dd^perp_varphi elliptic with noise}, \ref{fig:f=dd^perp_varphi parabolic with noise}, \ref{fig:f=dd^perp_varphi hyperbolic with noise}. Depending on the choice of $\vu$, the PDE appearing in \eqref{eq:f=dd^perp_varphi_M} is elliptic $\displaystyle\left(\vu=(\cos{\pi/6},\sin{\pi/6})\right)$, parabolic $\displaystyle\left(\vu=(\cos{\pi/4},\sin{\pi/4})\right)$, or hyperbolic $\displaystyle\left(\vu=(\cos{\pi/3},\sin{\pi/3})\right)$. The numerical solutions of such PDEs are discussed in Section \ref{sec: data formation and numerical methods}. Since the reconstructions are quite good even with $10\%$ of added noise, we have not included the images of reconstructions using data without noise. The relative errors of reconstructions without noise are included in the appropriate tables.

\begin{rem}\label{rem:5-pixels}
    To recover $\varphi$ from $\Mc\vf$ when $\vf= \D \D^\perp \varphi,$ we need to solve the PDE \eqref{eq:f=dd^perp_varphi_M}, where the source term $D_\vu D_\vv \Mc\vf$ is known to us. In particular, it is known (see \cite{Gaik_Indrani_Rohit_24}) that the source term is supported inside a disc contained in the square domain of the data. We use this fact to assign the value of 0 to all entries of $D_\vu D_\vv \Mc\vf$ within 5 pixels from the boundary of that domain. This helps to substantially reduce the artifacts in the reconstruction due to the errors of numerical differentiation at the boundary of the domain.
\end{rem}

 The horizontal artifacts in the reconstructions in Figure \ref{fig:f=dd^perp_varphi elliptic with noise} are due to the cut of the transform data at the boundary (and its absence outside of) the square $[-1,1]\times[-1,1]$. In fact, the essential part here is the cut of the data along the right side of the square (see formula \ref{eq:f=dd^perp_varphi_LT}). Such artifacts can be easily eliminated by adjusting the size of the data domain according to the support of the image field, so that the part of the boundary, where the data is cut, is located below the support of the image field. 

\begin{table}[h!]
\begin{center}
\begin{tabular}{ |c|c|c|c| } 
 \hline
   & From $\Lc\vf$ & From $\Tc\vf$ & From $\Mc\vf$\\ 
 \hline
  No Noise & 79.09\%& 79.09\%& 12.37\%\\ 
 \hline
 $10\%$ Noise & 81.40\%& 80.66\%& 17.45\%\\ 
 \hline
\end{tabular}
\caption{Relative errors of reconstruction of $\varphi$ from $\Lc\vf$, $\Tc\vf$ or $\Mc\vf$ when $\vf= \D\D^{\perp}\varphi$ and $u_1 >u_2$.}
\end{center}
\end{table}

\begin{table}[h!]
\begin{center}
\begin{tabular}{ |c|c|c|c| } 
 \hline
   & From $\Lc\vf$ & From $\Tc\vf$ & From $\Mc\vf$\\ 
 \hline
  No Noise & 4.99\%& 4.99\%& 6.03\%\\ 
 \hline
 $10\%$ Noise & 7.48\%& 5.13\%& 7.65\%\\ 
 \hline
\end{tabular}
\caption{Relative errors of reconstruction of $\varphi$ from $\Lc\vf$, $\Tc\vf$ or $\Mc\vf$ when $\vf= \D\D^{\perp}\varphi$ and $u_1 =u_2$.}
\end{center}
\end{table}

\begin{table}[h!]
\begin{center}
\begin{tabular}{ |c|c|c|c| } 
 \hline
   & From $\Lc\vf$ & From $\Tc\vf$ & From $\Mc\vf$\\ 
 \hline
  No Noise & 12.09\%& 12.09\%& 7.89\%\\ 
 \hline
 $10\%$ Noise & 30.39\%& 14.02\%& 14.66\%\\ 
 \hline
\end{tabular}
\caption{Relative errors of reconstruction of $\varphi$ from $\Lc\vf$, $\Tc\vf$ or $\Mc\vf$ when $\vf= \D\D^{\perp}\varphi$ and $u_1 < u_2$.}
\end{center}
\end{table}

\subsubsection{Tensor fields of the form \texorpdfstring{$\vf$}{f} = \texorpdfstring{$\D \vg$}{dg} or \texorpdfstring{$\vf$}{f} = \texorpdfstring{$\D^\perp \vg$}{dperp g}}\label{subsec:dg and dperpg}
In this subsection, we consider 2-tensor fields of the form $\vf= \D\vg ~\text{or}~ \vf=\D^\perp \vg$, where $\vg$ is a vector field. We validate the methods obtained in \cite{Gaik_Indrani_Rohit_24} for efficient reconstructions of $\vg$ using certain combinations of the VLTs of $\vf$ with and without noise. Since the reconstructions in the case $ \vf = \D^\perp \vg$ are almost identical to those in the case $ \vf = \D\vg$, we have not included the reconstruction images for $\D^\perp \vg.$

\begin{thr}\label{th:special tensor recovery 1} Let $\vg = (g_1, g_2)$ be a vector field with components $g_i(\vx) \in C_c^2(D_1)$, for $i =  1, 2$.
\begin{itemize}
\item[(a)] If $\vf$ is a symmetric 2-tensor field of the form $\vf = \D \vg$, then it can be recovered explicitly in terms of $\Lc \vf$ and $\Mc \vf$ as follows:
\begin{align}\label{eq:g2}
   g_2(\vx) =  - \frac{1}{2u_2}\Lc \vf(\vx),  
\end{align}
and $g_1$ can be recovered by solving a second-order partial differential equation
 \begin{align}\label{eq:g1_ellip} 
 2 u_1^2 \frac{\partial^2 g_1}{\partial x_1^2}(\vx)  + (u_1^2 - u_2^2)  \frac{\partial^2 g_1}{\partial x_2^2}(\vx)  =   - \frac{1}{2u_2}D_\vu D_\vv h (\vx),
\end{align}
where $\displaystyle h (\vx)=2\Mc \vf (\vx) + \Vc \left(\frac{\partial g_2}{\partial x_1}\right) (\vx)$
with additional homogeneous boundary (or initial) conditions. More  explicitly, we have the following three cases:
\begin{enumerate}[label=(\roman*)]
\item When $u_1^2 >  u_2^2$, we have an elliptic  PDE \eqref{eq:g1_ellip}
with homogeneous boundary conditions.
    \item When $u_1^2 =  u_2^2$, the partial differential equation \eqref{eq:g1_ellip} becomes:
    \begin{align}\label{eq:g1_para}
         2 u_1^2 \frac{\partial^2 g_1}{\partial x_1^2}(\vx)   =   - \frac{1}{2u_2}D_\vu D_\vv h (\vx),
    \end{align}
    which can be solved for $g_1$ by integrating twice along $\vev_1$-direction. 

 \item When $u_1^2 <  u_2^2$, we have a hyperbolic PDE \eqref{eq:g1_ellip}, which can be solved by choosing appropriate initial conditions on $g_1$ (for instance, we can take $g_1(a, y) = 0$ and $\displaystyle \frac{\partial g_1}{\partial x_1} (a, y) = 0$, for $y\in \mathbb{R}$ and any fixed $a \in \mathbb{R}\setminus [-1, 1]$).   
\end{enumerate}

\item[(b)] If $\vf$ is a symmetric 2-tensor field of the form $\vf = \D^\perp \vg$, then it can be recovered explicitly in terms of $\Tc \vf$ and $\Mc \vf$ as follows:
\begin{align}\label{eq:g1(dperpg)}
  g_1(\vx) = \frac{1}{2u_2}\Tc \vf(\vx) 
  \end{align}
 and $g_2$ can be recovered by solving the following second-order partial differential equation (equipped with appropriate boundary or initial conditions):
\begin{align}\label{eq:g2(dperpg)}
    2 u_1^2 \frac{\partial^2 g_2}{\partial x_1^2}(\vx)  + (u_1^2 - u_2^2)  \frac{\partial^2 g_2}{\partial x_2^2}(\vx)  =   - \frac{1}{2u_2}D_\vu D_\vv \Tilde{h} (\vx),
\end{align}
where $\displaystyle \Tilde{h}(\vx) = -2\Mc \vf(\vx) - \Vc\left(\frac{\partial g_1}{\partial x_1}\right)(\vx)$.
\end{itemize}
\end{thr}

\begin{figure}[H]
    \centering
     \begin{subfigure}[b]{0.89\textwidth}
         \centering
         \includegraphics[width=\textwidth]{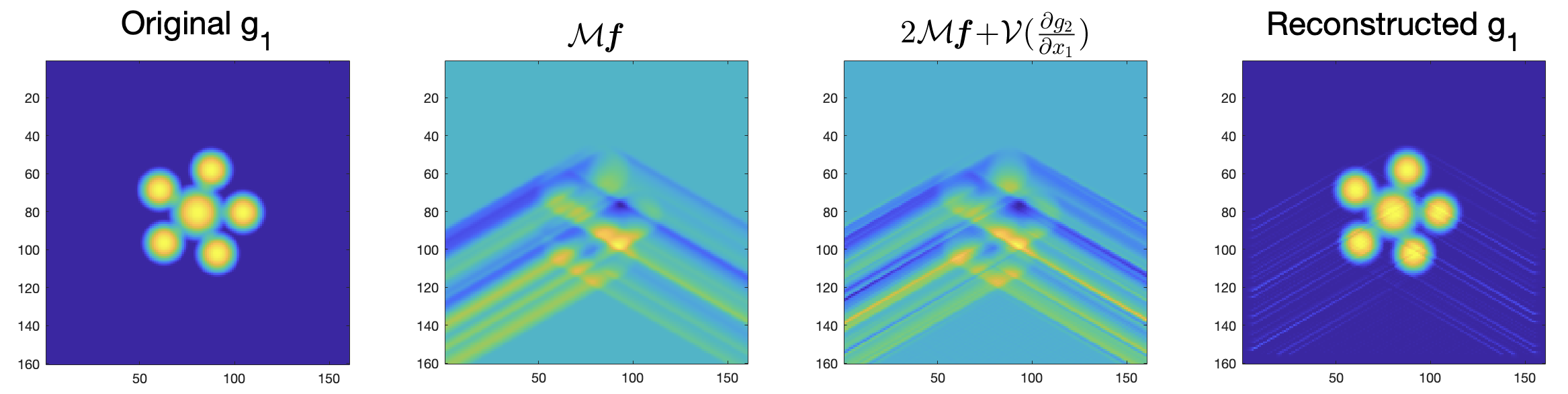}
     \end{subfigure}
     \vfill
     \begin{subfigure}[b]{0.89\textwidth}
         \centering
         \includegraphics[width=\textwidth]{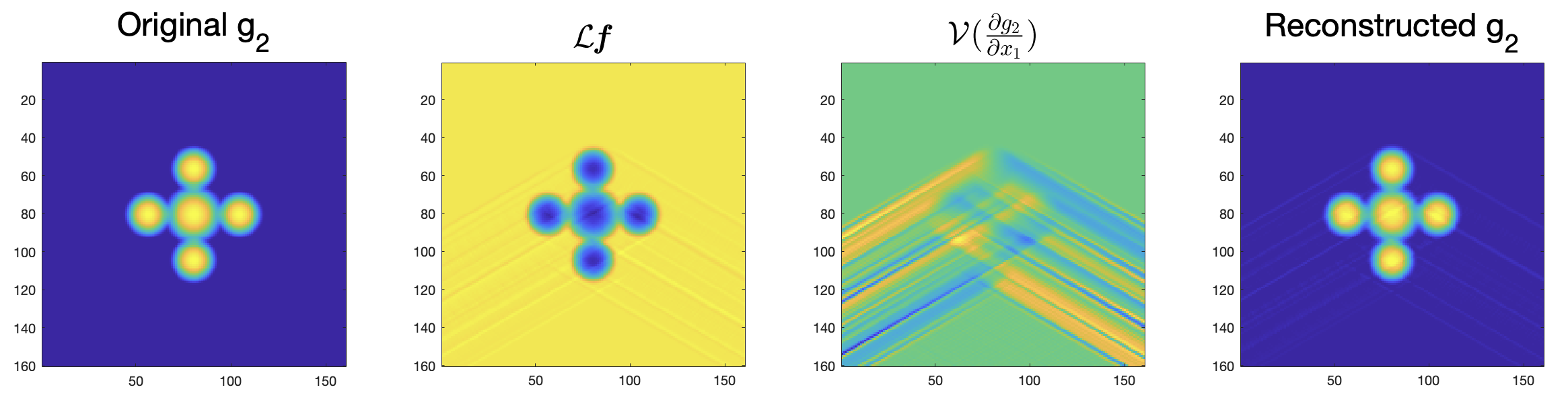}
     \end{subfigure}
        \caption{Recovery of $\vg$ from $\Lc\vf$ and $\Mc\vf$ when $\vf=\D\vg$, $u_1 >u_2$ (elliptic), using \eqref{eq:g2}, \eqref{eq:g1_ellip}.}\label{fig:ph1 dg (elliptic)}
\end{figure}

\begin{figure}[H]
     \centering
     \begin{subfigure}[b]{0.89\textwidth}
         \centering
         \includegraphics[width=\textwidth]{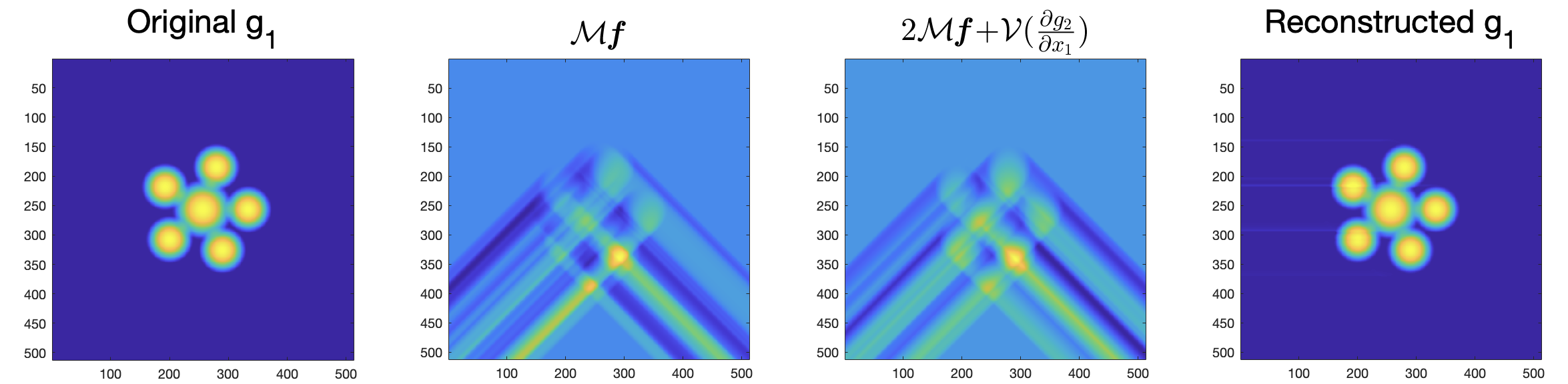}
     \end{subfigure}
     \vfill
     \begin{subfigure}[b]{0.89\textwidth}
         \centering
         \includegraphics[width=\textwidth]{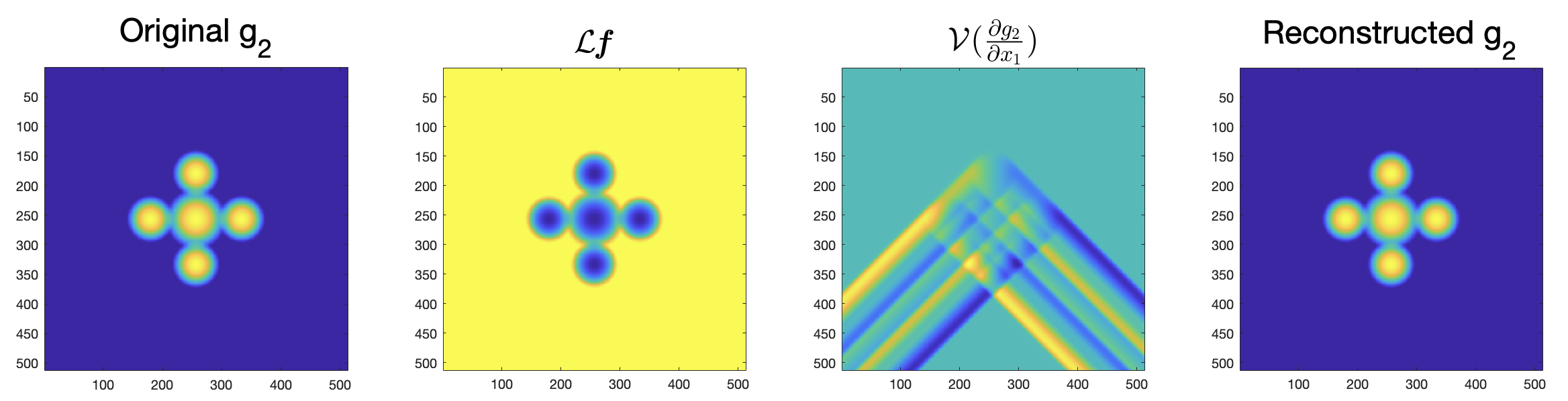}
     \end{subfigure}
        \caption{Recovery of $\vg$ from $\Lc\vf$ and $\Mc\vf$ when $\vf=\D\vg$, $u_1 = u_2$ (parabolic), using  \eqref{eq:g2}, \eqref{eq:g1_para}.}\label{fig:ph1 dg (parabolic)}
\end{figure}

\begin{figure}[H]
     \centering
     \begin{subfigure}[b]{0.89\textwidth}
         \centering
         \includegraphics[width=\textwidth]{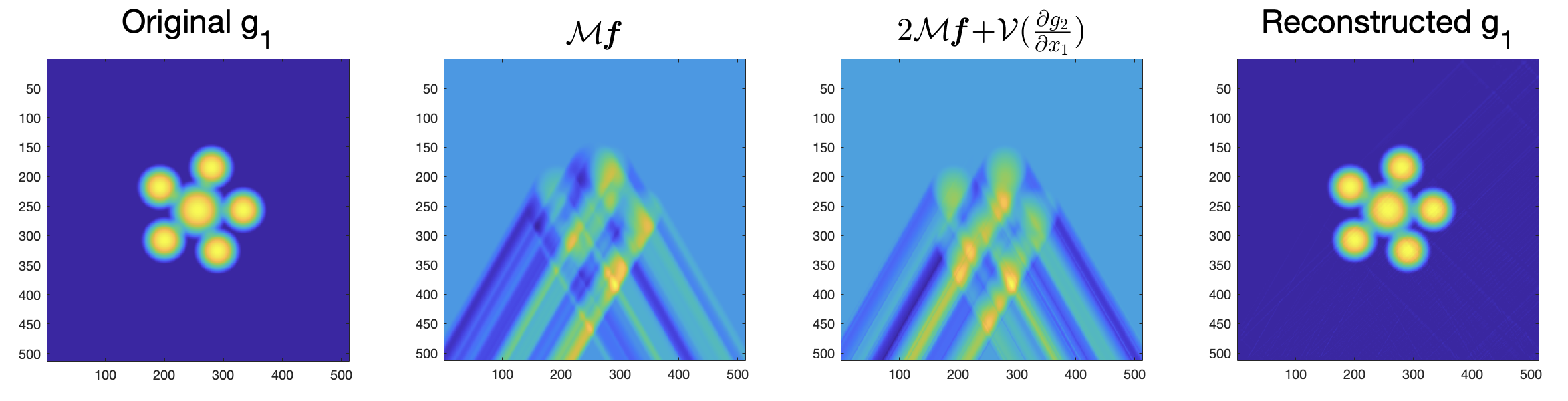}
     \end{subfigure}
     \vfill
     \begin{subfigure}[b]{0.89\textwidth}
         \centering
         \includegraphics[width=\textwidth]{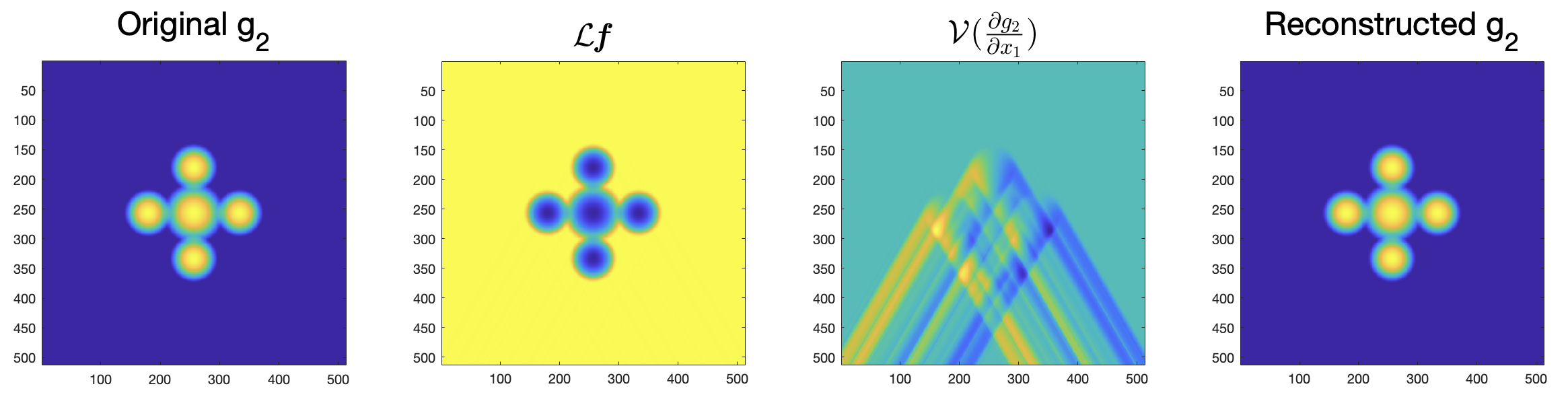}
     \end{subfigure}
        \caption{Recovery of $\vg$ from $\Lc\vf$ and $\Mc\vf$ when $\vf=\D\vg$, $u_1 < u_2$ (hyperbolic), using \eqref{eq:g2}, \eqref{eq:g1_ellip}.}\label{fig:ph1 dg (hyperbolic)}
\end{figure}

\begin{figure}[H]
     \centering
     \begin{subfigure}[b]{0.19\textwidth}
         \centering
         \includegraphics[width=\textwidth]{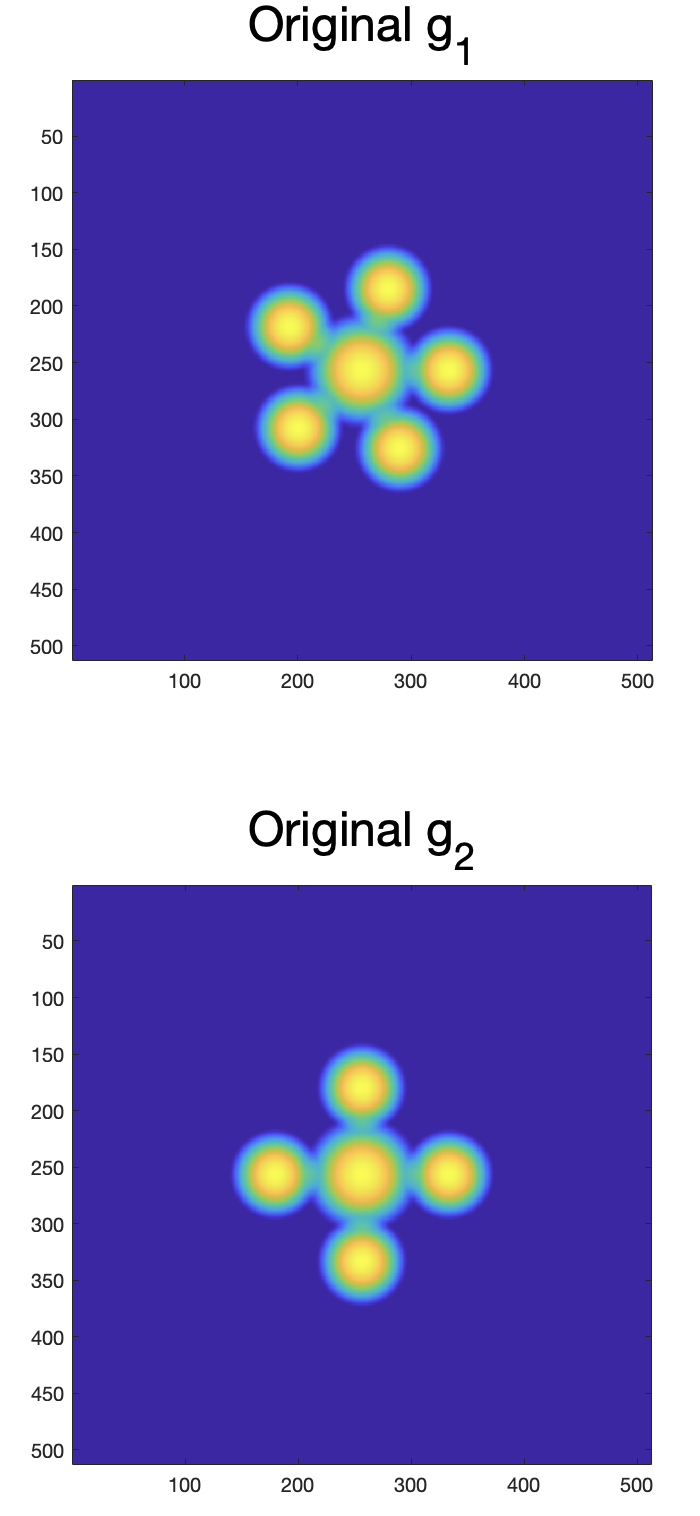}
     \end{subfigure}
     \hfill
     \begin{subfigure}[b]{0.19\textwidth}
         \centering
         \includegraphics[width=1.02\textwidth]{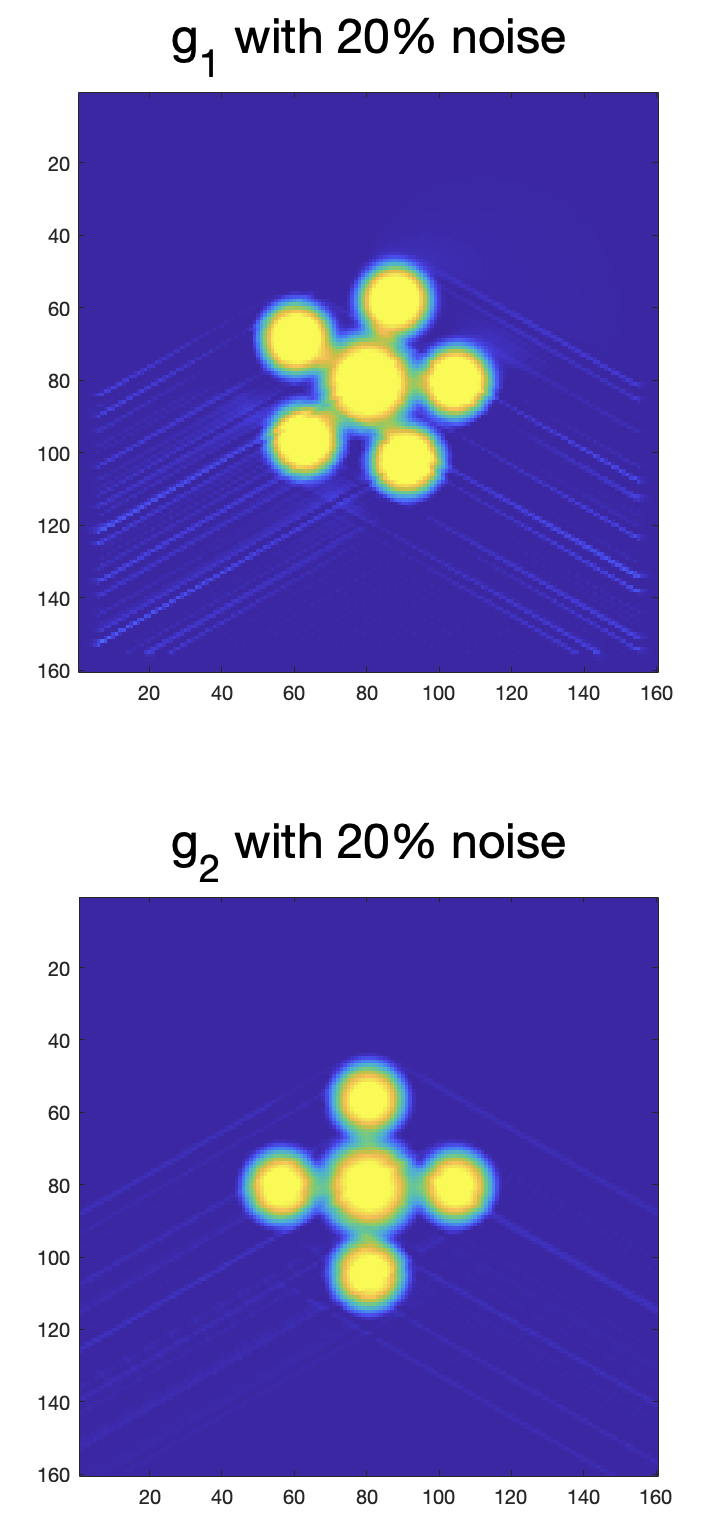}
     \end{subfigure}
     \hfill
     \begin{subfigure}[b]{0.19\textwidth}
         \centering
         \includegraphics[width=0.99\textwidth]{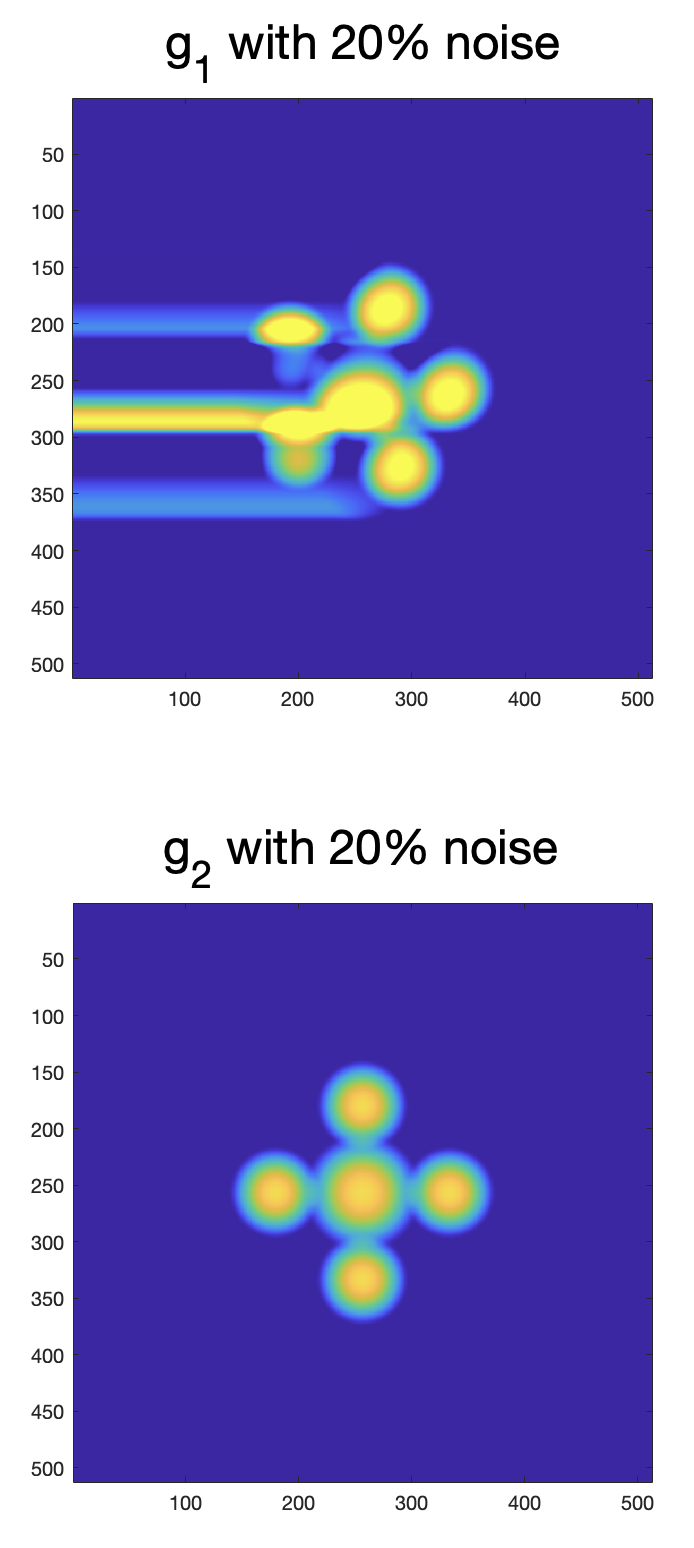}
     \end{subfigure}
     \hfill
     \begin{subfigure}[b]{0.19\textwidth}
         \centering
         \includegraphics[width=1.01\textwidth]{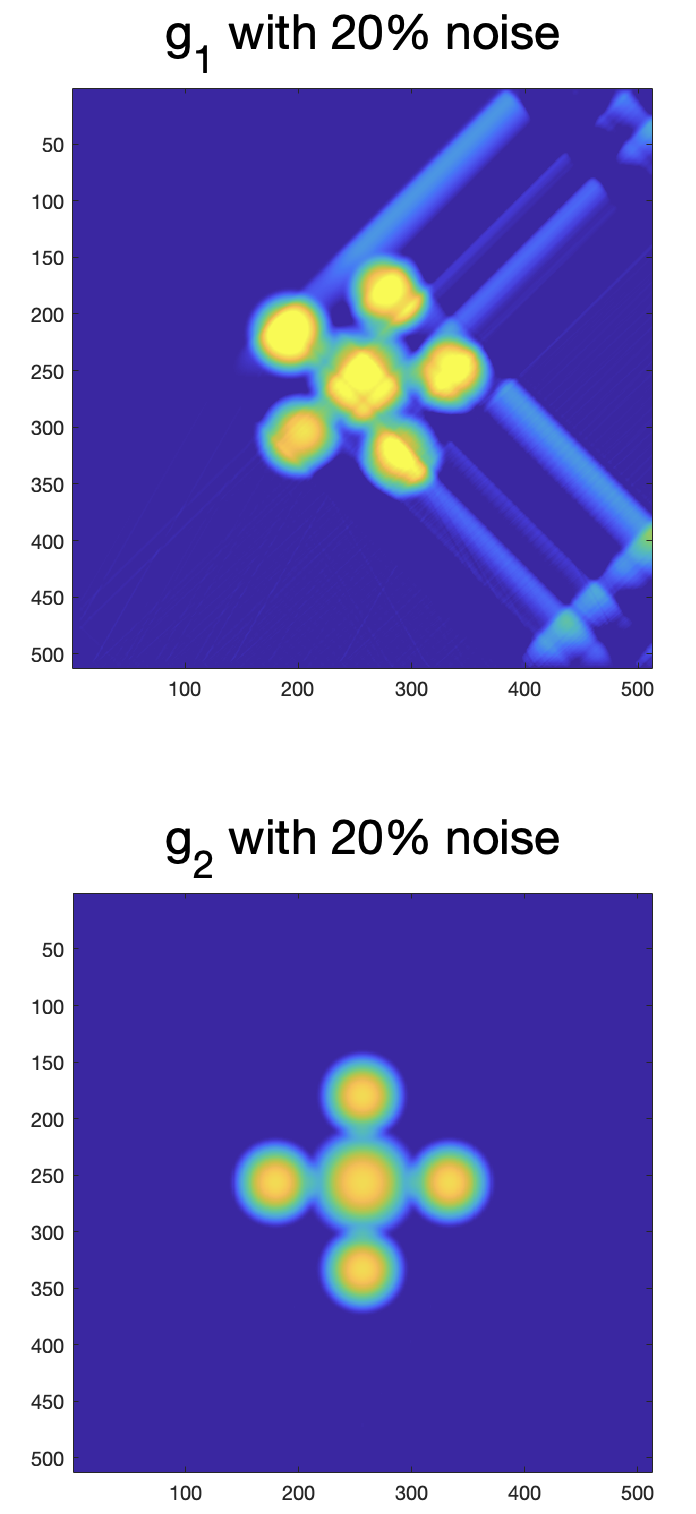}
     \end{subfigure}
     \caption{Column 1: original $g_1, g_2$; Column 2: reconstr. $g_1, g_2$, $u_1>u_2$ (elliptic), $20\%$ noise; Column 3: reconstr. $g_1, g_2$, $u_1 = u_2$ (parabolic), $20\%$ noise; Column 4: reconstr. $g_1, g_2$, $u_1 < u_2$ (hyperbolic),  $20\%$ noise.}\label{fig:ph1 dg with noise}
\end{figure}

\begin{figure}[H]
     \centering
     \begin{subfigure}[b]{0.9\textwidth}
         \centering
         \includegraphics[width=\textwidth]{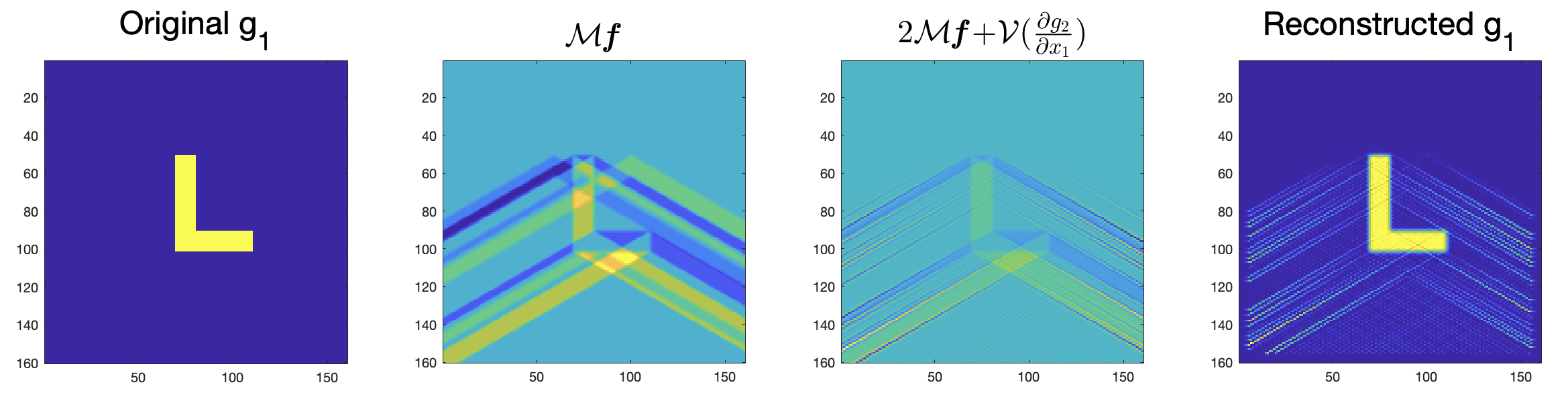}
     \end{subfigure}
     \vfill
     \begin{subfigure}[b]{0.9\textwidth}
         \centering
         \includegraphics[width=\textwidth]{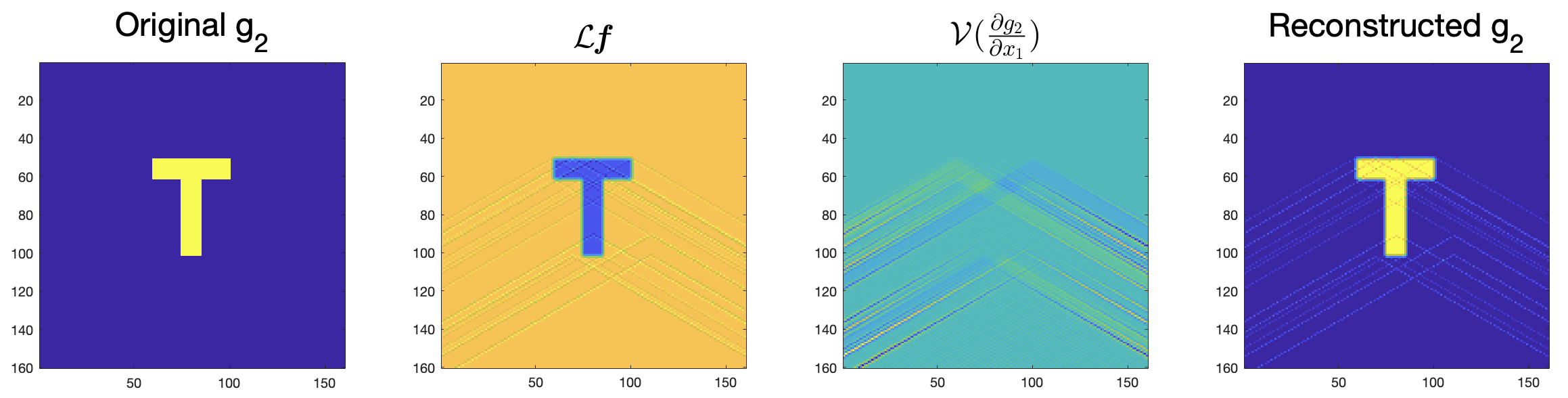}
     \end{subfigure}
        \caption{Recovery of $\vg$ from $\Lc\vf$ and $\Mc\vf$  when $\vf=\D\vg$, $u_1 >u_2$ (elliptic), using  \eqref{eq:g2}, \eqref{eq:g1_ellip}.}\label{fig:ph2 dg (elliptic)}
\end{figure}

\begin{figure}[H]
     \centering
     \begin{subfigure}[b]{0.9\textwidth}
         \centering
         \includegraphics[width=\textwidth]{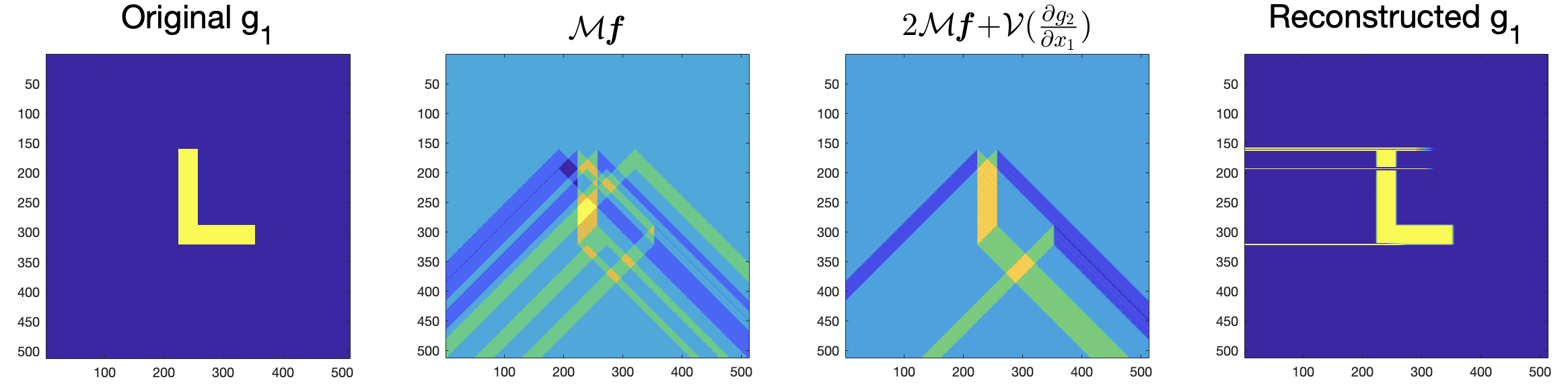}
     \end{subfigure}
     \vfill
     \begin{subfigure}[b]{0.9\textwidth}
         \centering
         \includegraphics[width=\textwidth]{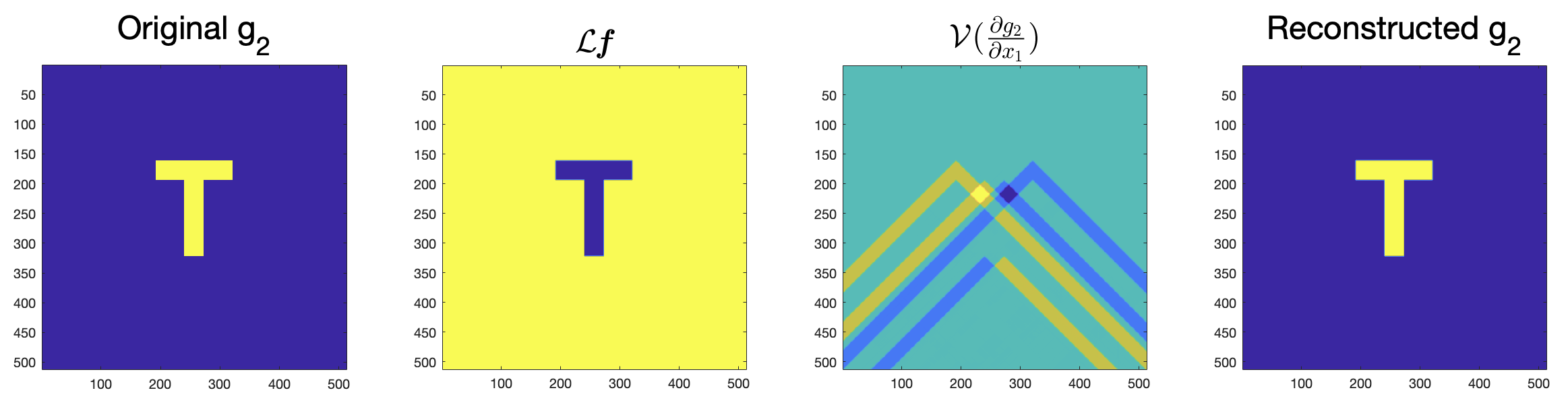}
     \end{subfigure}
        \caption{Recovery of $\vg$ from $\Lc\vf$ and $\Mc\vf$ when $\vf=\D\vg$, $u_1 = u_2$ (parabolic), using \eqref{eq:g2}, \eqref{eq:g1_para}.}\label{fig:ph2 dg (parabolic)}
\end{figure}

\begin{figure}[H]
     \centering
     \begin{subfigure}[b]{0.87\textwidth}
         \centering
         \includegraphics[width=\textwidth]{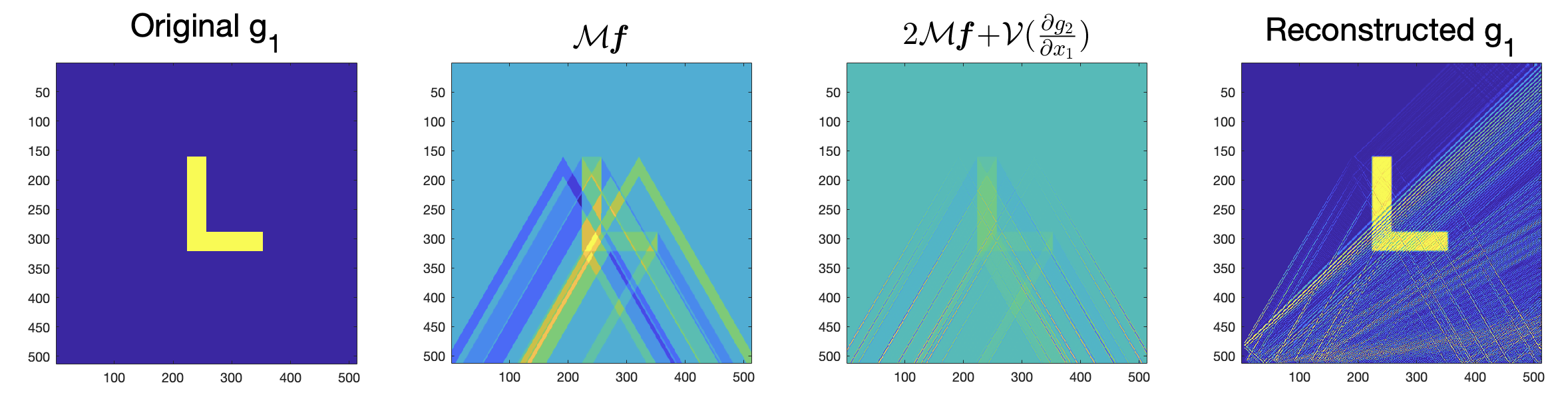}
     \end{subfigure}
     \vfill
     \begin{subfigure}[b]{0.87\textwidth}
         \centering
         \includegraphics[width=\textwidth]{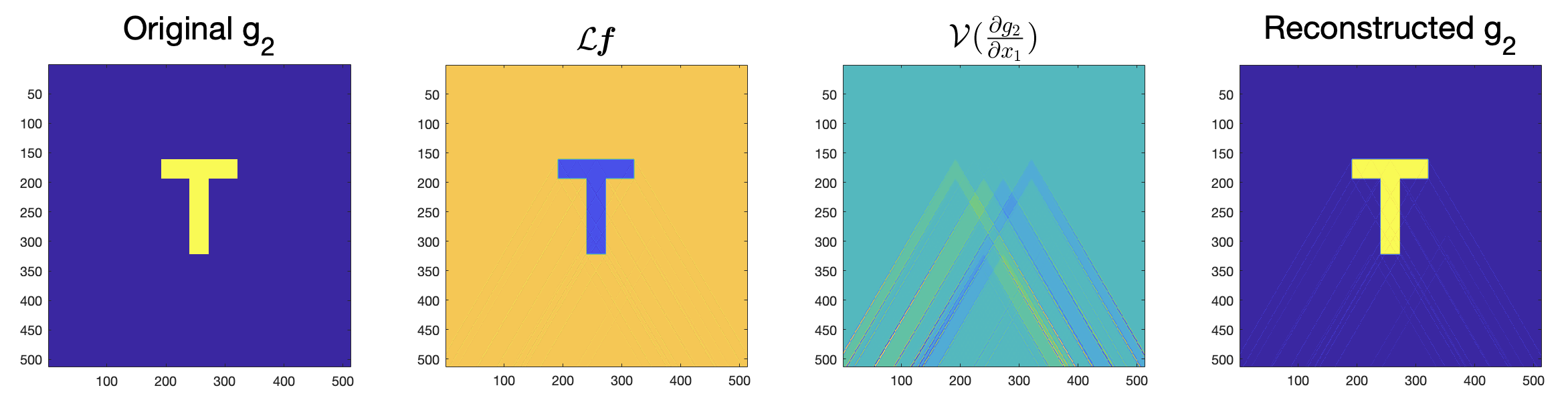}
     \end{subfigure}
        \caption{Recovery of $\vg$ from $\Lc\vf$ and $\Mc\vf$ when $\vf=\D\vg$, $u_1 < u_2$ (hyperbolic), using \eqref{eq:g2},\eqref{eq:g1_ellip}.}\label{fig:ph2 dg (hyperbolic)}
\end{figure}

\begin{figure}[H]
     \centering
     \begin{subfigure}[b]{0.21\textwidth}
         \centering
         \includegraphics[width=0.91\textwidth]{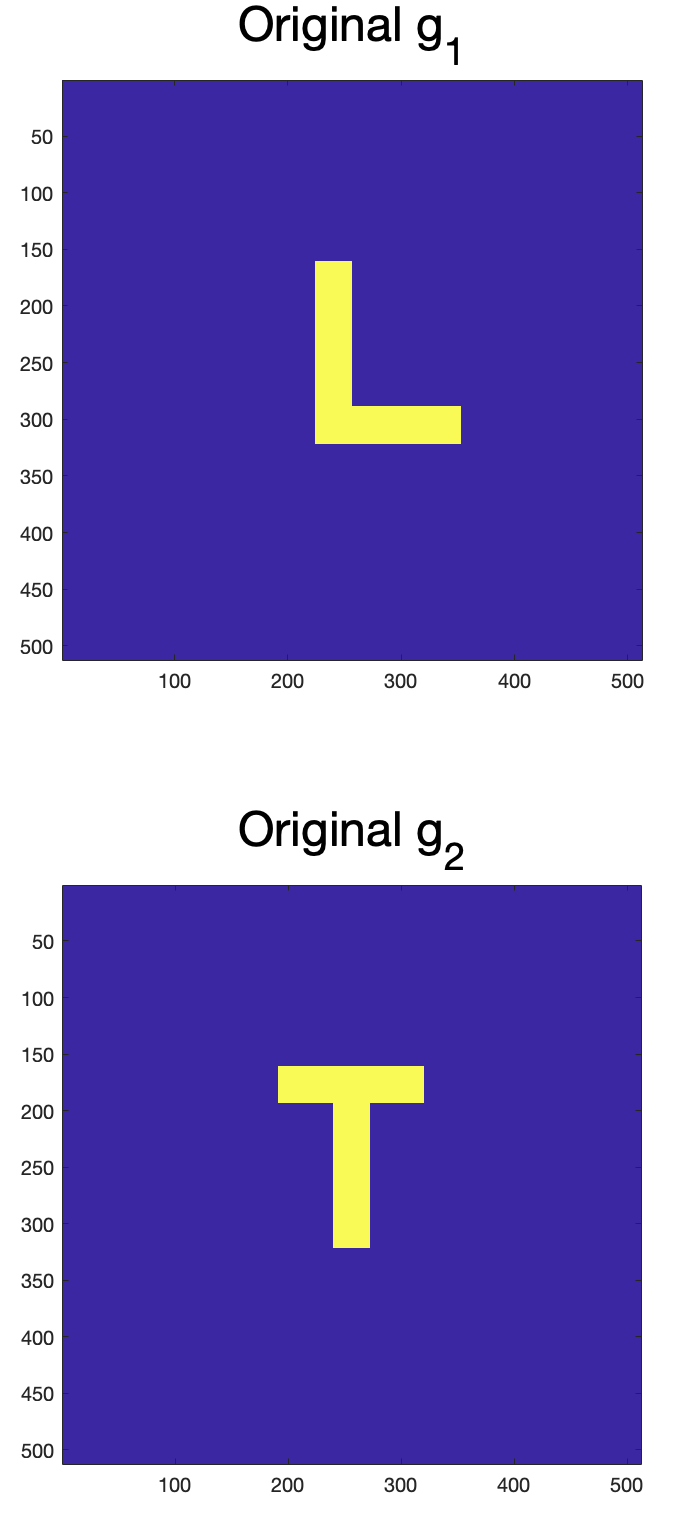}
     \end{subfigure}
     \hfill
     \begin{subfigure}[b]{0.21\textwidth}
         \centering
         \includegraphics[width=0.93\textwidth]{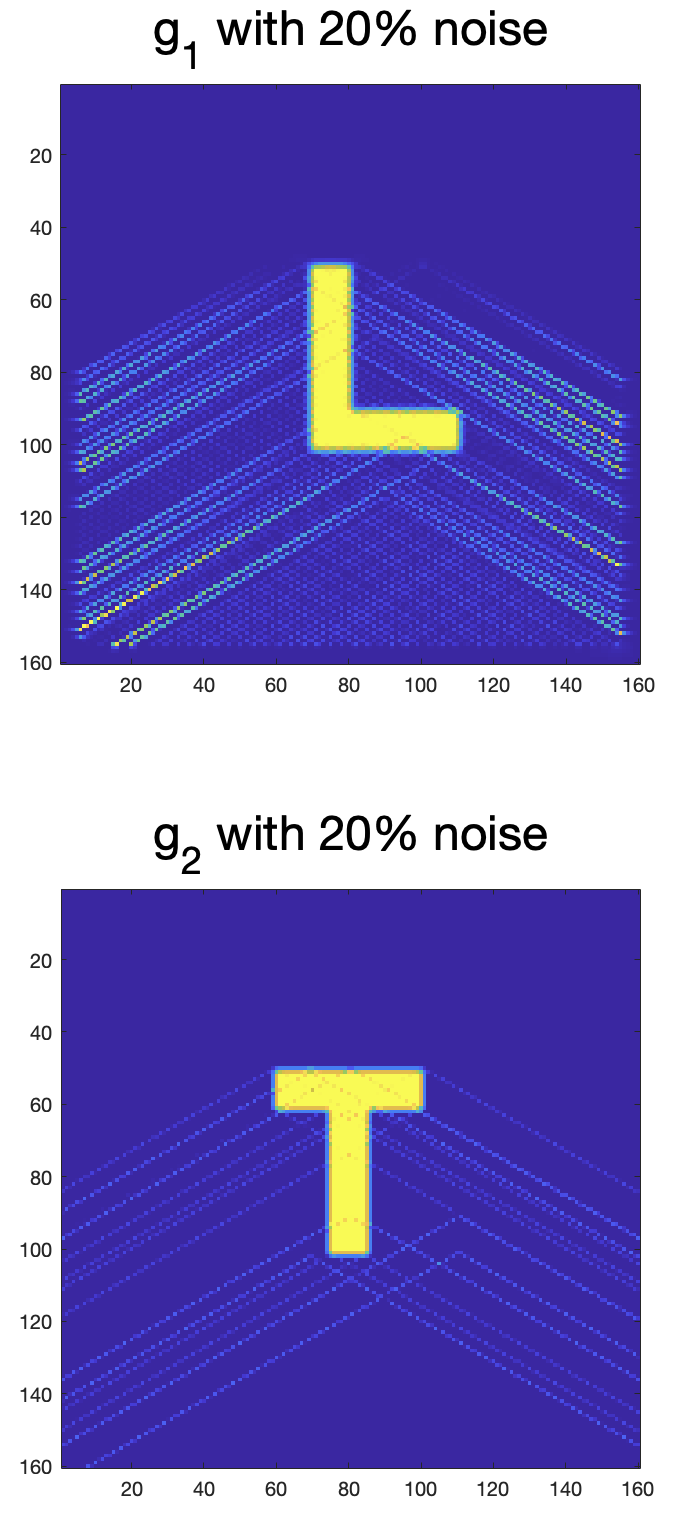}
     \end{subfigure}
     \hfill
     \begin{subfigure}[b]{0.21\textwidth}
         \centering
         \includegraphics[width=0.9\textwidth]{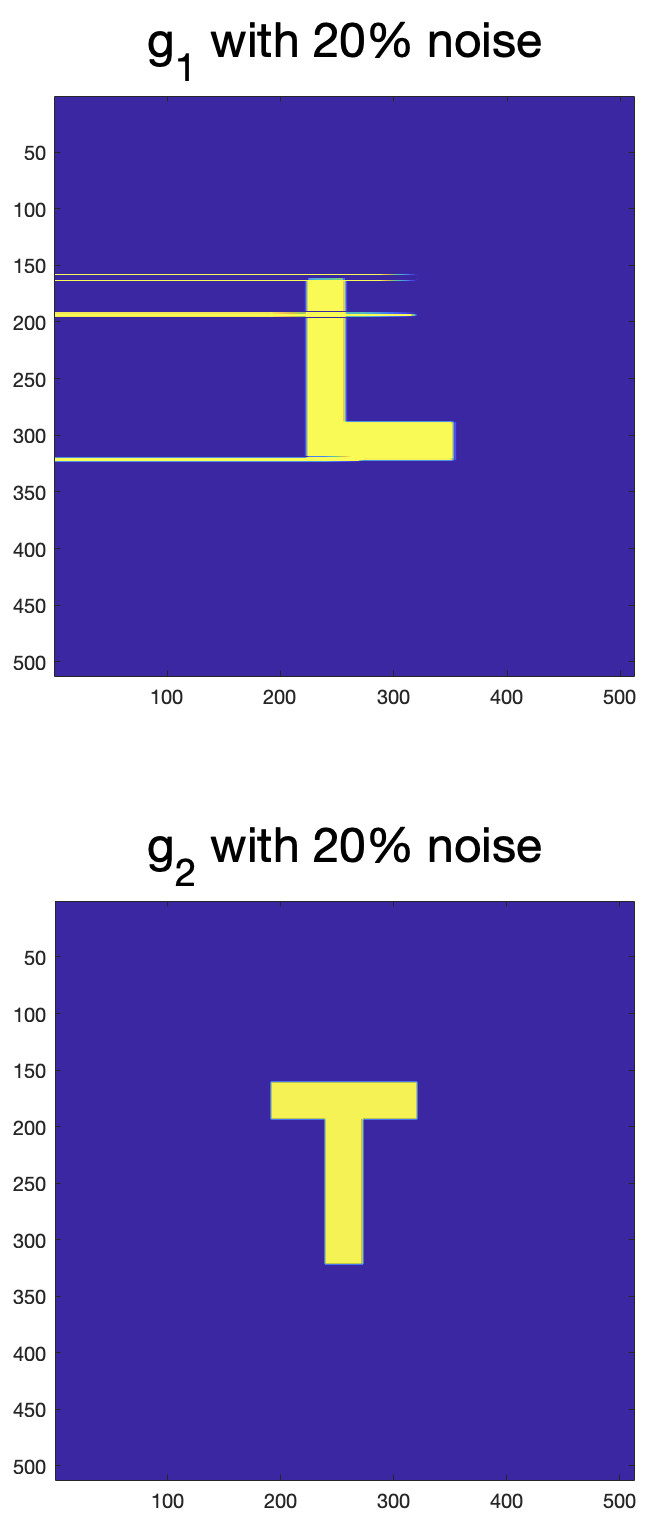}
     \end{subfigure}
     \hfill
     \begin{subfigure}[b]{0.21\textwidth}
         \centering
         \includegraphics[width=0.91\textwidth]{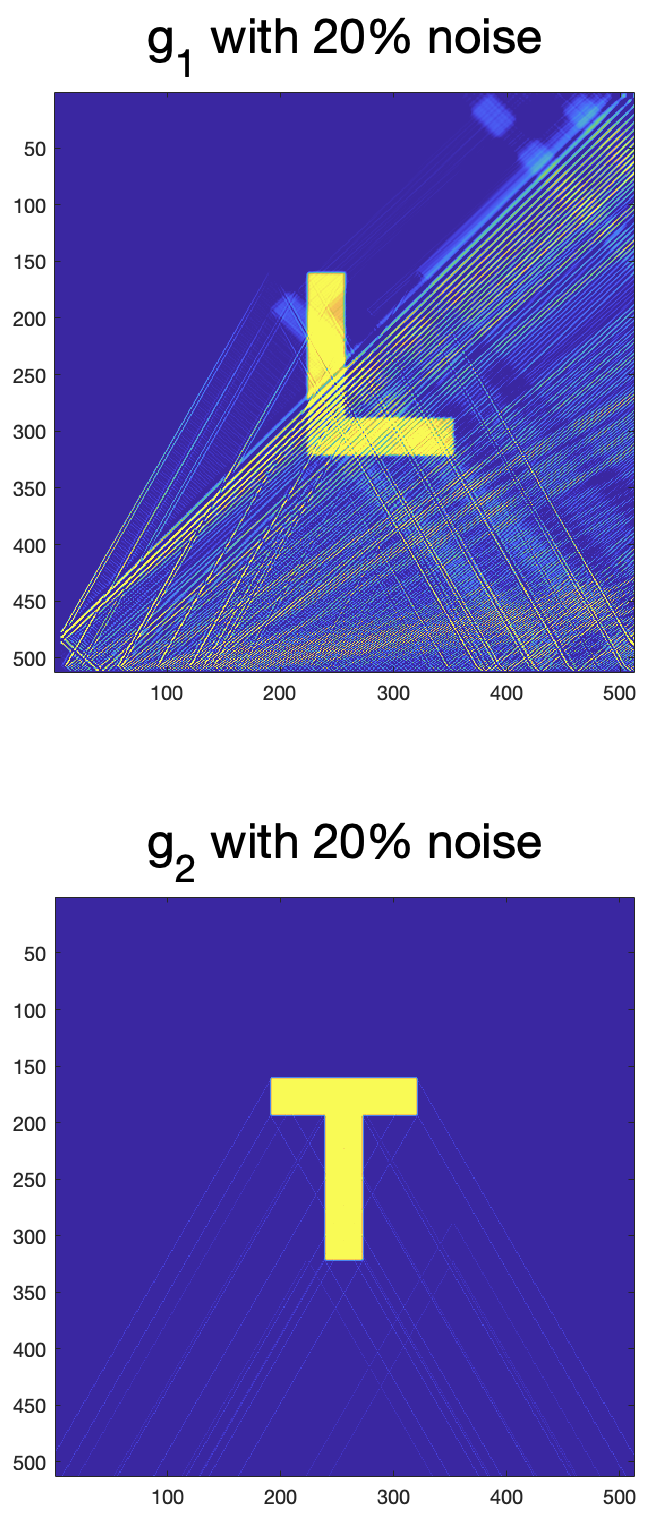}
     \end{subfigure}
     \caption{Column 1: original $g_1, g_2$; Column 2: reconstr. $g_1, g_2$, $u_1>u_2$ (elliptic), $20\%$ noise; Column 3: reconstr. $g_1, g_2$,  $u_1 = u_2$ (parabolic), $20\%$ noise; Column 4: reconstr. $g_1, g_2$, $u_1 < u_2$ (hyperbolic), $20\%$ noise.}\label{fig:ph2 dg with noise}
\end{figure}

In Figures \ref{fig:ph1 dg (elliptic)}, \ref{fig:ph1 dg (parabolic)}, \ref{fig:ph1 dg (hyperbolic)}, $g_{2}$ is reconstructed directly from $\Lc\vf$ (see \eqref{eq:g2}), while the recovery of $g_1$ requires either solving an elliptic ($\displaystyle\vu=(\cos{\pi/6},\sin{\pi/6})$)  or hyperbolic ($\displaystyle\vu=(\cos{\pi/3},\sin{\pi/3})$) PDE \eqref{eq:g1_ellip}, or double integration solving the parabolic PDE \eqref{eq:g1_para} (for $\displaystyle\vu=(\cos{\pi/4},\sin{\pi/4})$). Figure \ref{fig:ph1 dg with noise} presents the reconstructions from data with $20\%$ added noise. 
The same algorithms are used in Figures \ref{fig:ph2 dg (elliptic)}, \ref{fig:ph2 dg (parabolic)}, \ref{fig:ph2 dg (hyperbolic)}, \ref{fig:ph2 dg with noise} for the non-smooth phantom. 

\begin{rem}
    To recover $g_1$ when $\vf=\D\vg$, we are solving PDEs, in which the source term $D_\vu D_\vv h$ is known to us. Following the ideas discussed in Remark \ref{rem:5-pixels}, we assign the value of 0 to all entries of $D_\vu D_\vv h$ within 5 pixels from the boundary of the data domain.
\end{rem}

\begin{table}[h!]
\begin{center}
\begin{tabular}{ |c|c|c|c| } 
 \hline
  $\vg$ & 
 $20\%$ Noise - Elliptic & $20\%$ Noise - Parabolic & $20\%$ Noise - Hyperbolic\\ 
 \hline
 $g_1$  &33.75\%  & 107.46\%  &37.59\% \\ 
 \hline
 $g_2$  &13.80\%  & 8.72\% & 8.72\% \\  
 \hline
\end{tabular}
\caption{Relative errors of reconstruction of $\vg$ from $\Lc\vf$ and $\Mc\vf$ when $\vf= \D\vg$.}
\end{center}
\end{table}

The horizontal artifacts in the reconstruction of $g_1$ in Figure \ref{fig:ph2 dg (parabolic)} are due to the errors of numerical differentiation of $h(\vx)$ in the right-hand side of equation \eqref{eq:g1_para}, propagated by the double integration with respect to $x_1$. Although visually the horizontal artifacts do not seem to cause a severe distortion of the original image, their strength (i.e. the numerical value of the reconstructed function at the artifact pixels) is quite large, adversely affecting the relative error of the reconstruction (see Table \ref{tab:g_from_LM}).

\begin{table}[h!]
 \begin{center}
\begin{tabular}{ |c|c|c|c| } 
 \hline
 $\vg$  & $20\%$ Noise - Elliptic & $20\%$ Noise - Parabolic& $20\%$ Noise -Hyperbolic\\ 
 \hline
 $g_1$  & 30.54\%  & 820.16\% & 50.93\%\\  
 \hline
  $g_2$  &13.55\%  & 7.80\% & 15.73\%\\  
 \hline
\end{tabular}
\caption{Relative errors of reconstruction of $\vg$ from $\Lc\vf$ and $\Mc\vf$ when $\vf= \D\vg$.}\label{tab:g_from_LM}
\end{center}
\end{table}
 \subsection{Full Recovery of \texorpdfstring{$\vf$}{f} from \texorpdfstring{$\Lc\vf$}{Lf}, \texorpdfstring{$\Tc\vf$}{Tf} and \texorpdfstring{$\Mc\vf$}{Mf}}\label{subsec: full recovery Lf, Mf, and Tf}
 In this subsection, we use $\Lc\vf, \Tc\vf$, and $\Mc\vf$ to recover $\vf.$ As in the previous sections, we start by discussing the theoretical results for this setup in Theorem \ref{eq:full-rec from L,T and M}.

\begin{thr}\label{eq:full-rec from L,T and M}
 Consider a symmetric $2$-tensor field $\vf\in C_c^2\left(S^2;D_1\right)$, $\vev_1=(1,0)$ and $\vev_2=(0,1)$. 
\begin{itemize}
   \item  For $u_1=u_2,$ $\vf$ can be recovered from $\Lc\vf,\Tc\vf$, and $\Mc\vf$ by
    \begin{align}
       f_{11}(\vx) &= \frac{1}{4 u_2} D_{\vu} D_{\vv}\Xc_{\vev_2} \left(\Lc \vf + \Tc \vf \right)(\vx) - \frac{1}{2u_1} D_{\vu} D_{\vv} \Xc_{-\vev_1}\Mc \vf\, (\vx),\label{eq:f11(u1=u2)}\\
      f_{12}(\vx) &= \frac{1}{4 u_1} D_\vu D_\vv \Xc_{-\vev_1} \left(\Lc \vf - \Tc \vf \right)(\vx),\label{eq:f12(u1=u2)}\\
       f_{22}(\vx) &= \frac{1}{4 u_2} D_{\vu} D_{\vv}\Xc_{\vev_2} \left(\Lc \vf + \Tc \vf \right)(\vx) + \frac{1}{2u_1} D_{\vu} D_{\vv} \Xc_{-\vev_1}\Mc \vf\, (\vx)\label{eq:f22(u1=u2)}.
    \end{align}
\item 
    For $u_1\ne u_2,$  $\vf$ can be reconstructed from $\Lc\vf,\Tc\vf$, and $\Mc\vf$  as follows.\\
  $f_{12}$ can be found by solving the elliptic boundary value problem:
\begin{align}\label{eq:elliptic equation for f12}
    \left\{\begin{array}{rll}
   a \partial_{x_1}^2 f_{12}+ b\partial_{x_2}^2 f_{12}&= -g      & \mbox{ in } D_1,  \\
     f_{12} &= 0    & \mbox{ on } \partial D_1,
    \end{array}\right.
\end{align}
where  $a=2 u_1^2\left[1+(u_1^2 -u_2^2)\right ]  > 0 $, ~ $b=\left( u_1^2-u_2^2\right)^{2} >  0$, 
and 
$$g = \frac{1}{2u_2}\left[u_1^2 \partial_{x_1} D_\vu D_\vv \left( \Tc \vf-\Lc \vf\right)+(u_1^2 -u_2^2)\partial_{x_2} D_\vu D_\vv  \Mc \vf\,\right]. $$
$f_{11}$ can be recovered from the knowledge of $\Lc\vf,\Tc\vf$, and the reconstructed $f_{12}$ by
\begin{align}\label{eq:f_11}
f_{11}=-\Xc_{\vev_2}\left[\left(u_2^2 D_\vu D_\vv  \Tc \vf - u_1^2 D_\vu D_\vv  \Lc \vf + 4 u_1^2 u_2 \partial_{x_1}  f_{12}\right)/\left( 2 u_2(u_1^2 -u_2^2)\right)\right].
\end{align}
   $f_{22}$ can be recovered from $\Lc\vf,\Tc\vf$, and the reconstructed $f_{11}$ by
   \begin{align}\label{eq:f_22}
f_{22}(\vx) = \frac{1}{2 u_2} D_\vu D_\vv\Xc_{\vev_2} \left(\Lc \vf + \Tc \vf \right)(\vx) -f_{11}(\vx).
\end{align}
\end{itemize}
\end{thr}

Figure \ref{fig:ph2 (u1 = u2) Lf, Tf, and Mf} depicts the reconstructions of Phantom 2 using formulas \eqref{eq:f11(u1=u2)}, \eqref{eq:f12(u1=u2)}, \eqref{eq:f22(u1=u2)}. Figure \ref{fig:ph2 (u1 = u2) Lf, Tf, and Mf with noise} shows the corresponding reconstructions in the presence of different levels of noise.

\begin{table}[h!]
\begin{center}
\begin{tabular}{ |c|c|c|c|c| } 
 \hline
  $\vf$  & No noise & $5\%$ Noise & $10\%$ Noise& $20\%$ Noise \\ 
 \hline
  $f_{11}$  & 5.22\% & 85.56\% & 194.44\% & 452.72\%\\ 
 \hline
 $f_{12}$ & 9.42\% & 30.25\% & 68.14\% & 158.89\%\\  
 \hline
 $f_{22}$ & 8.15\% & 65.07\% & 146.95\% & 344.88\% \\  
 \hline
\end{tabular}
\caption{Relative errors of reconstructing $f_{11}, f_{12}, f_{22}$ from $\Lc\vf, \Mc\vf$, $\Tc\vf$, $\vu=(\cos{\pi/4,\sin{\pi/4}})$.}\label{tab:f_from_LMT-1}
\end{center}
\end{table}

\begin{figure}[H]
     \centering
     \begin{subfigure}[b]{0.56\textwidth}
         \centering
         \includegraphics[width=\textwidth]{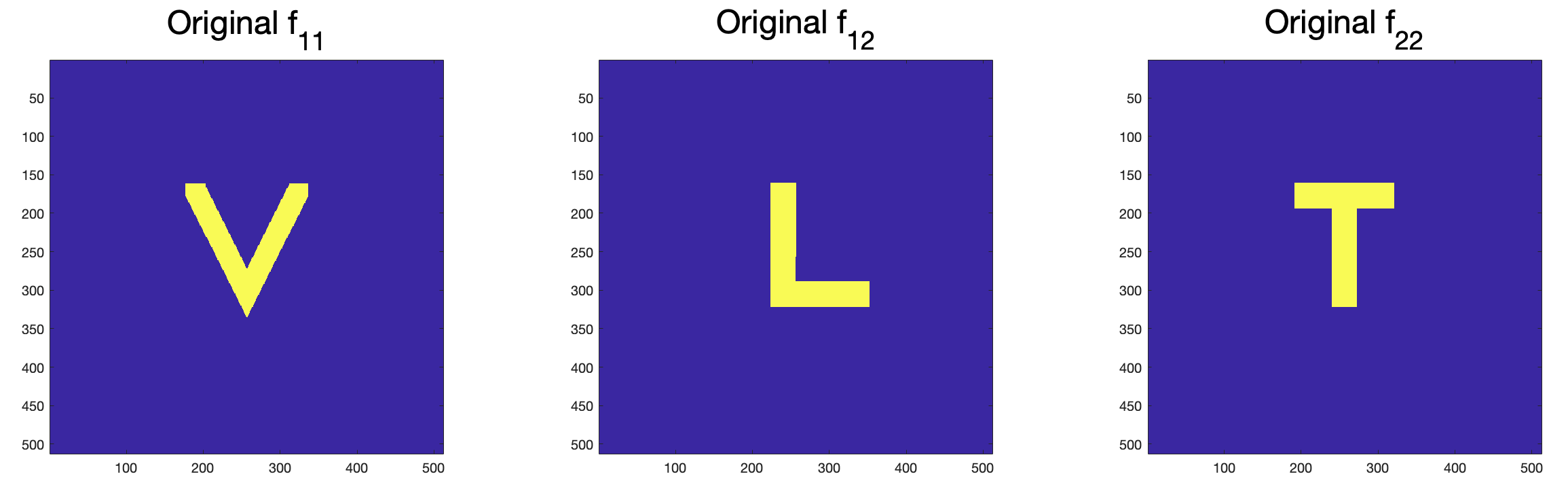}
     \end{subfigure}
     \vfill
     \begin{subfigure}[b]{0.56\textwidth}
         \centering
         \includegraphics[width=\textwidth]{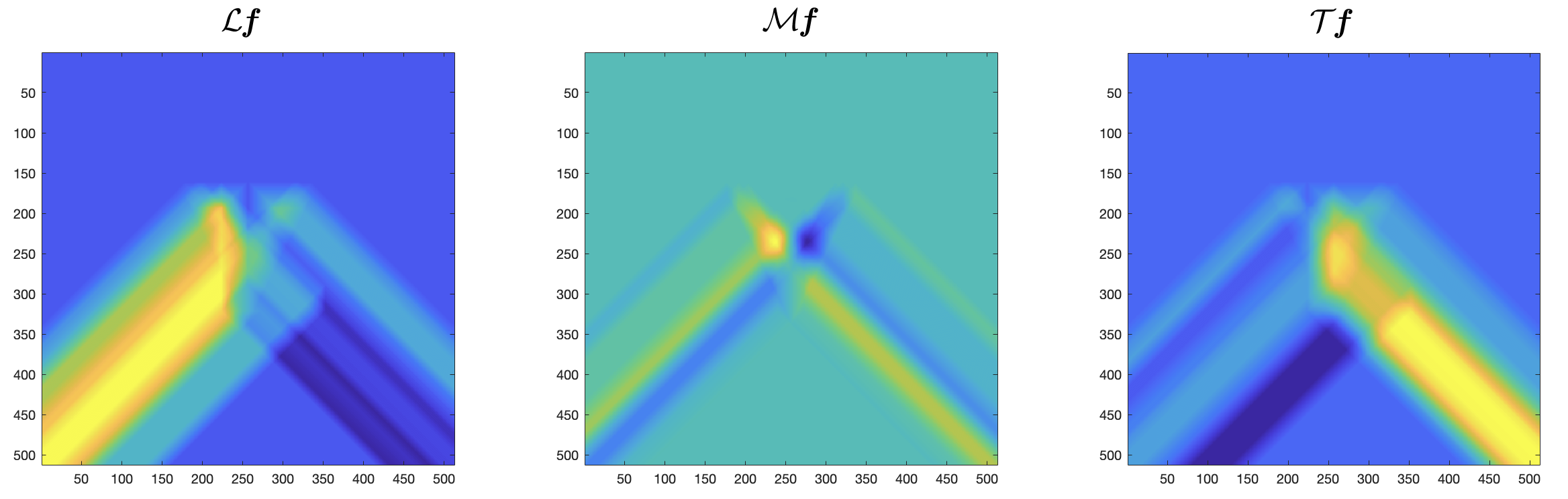}
     \end{subfigure}
     \vfill
     \begin{subfigure}[b]{0.56\textwidth}
         \centering
         \includegraphics[width=\textwidth]{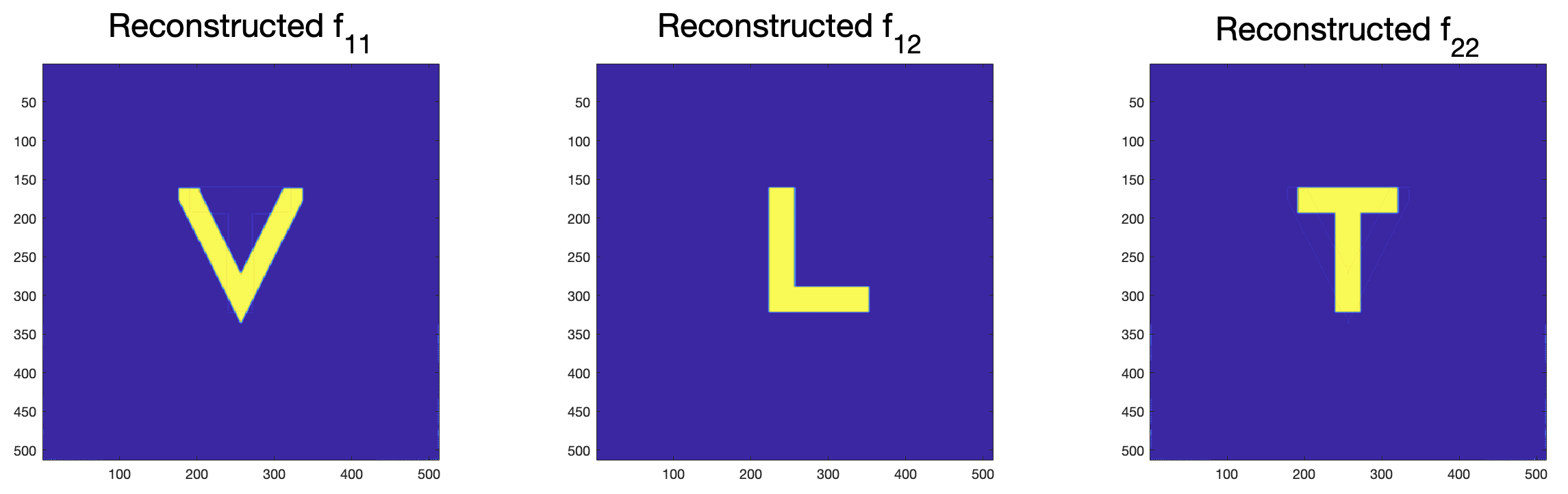}
     \end{subfigure}
\caption{Reconstruction of $\vf$ from $\Lc\vf,\Mc\vf,\Tc\vf$, $\vu=(\cos{\pi/4},\sin{\pi/4})$,
using formulas \eqref{eq:f11(u1=u2)}, \eqref{eq:f12(u1=u2)}, \eqref{eq:f22(u1=u2)}.}\label{fig:ph2 (u1 = u2) Lf, Tf, and Mf}
\end{figure}
 \begin{figure}[H]
     \centering
   \includegraphics[width= 0.9\textwidth]{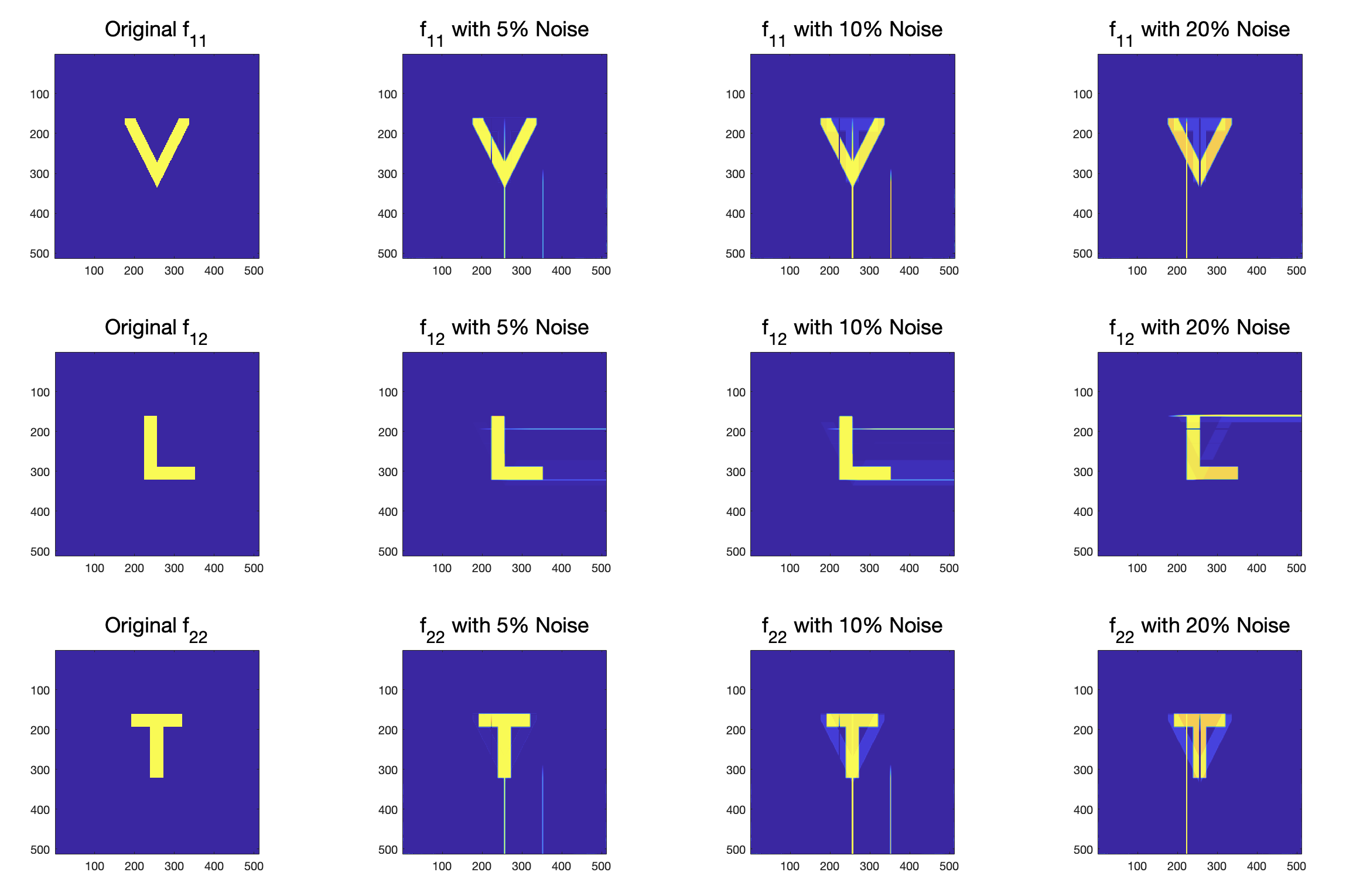}
    \caption{ Reconstructions from $\Lc\vf,\Mc\vf,\Tc\vf$, $\vu=(\cos{\pi/4},\sin{\pi/4})$ with $5\%, 10\%$, $20\%$ noise.}\label{fig:ph2 (u1 = u2) Lf, Tf, and Mf with noise}
 \end{figure}

We took $\vu=(\cos{\pi/3,\sin{\pi/3}})$ for implementations of \eqref{eq:elliptic equation for f12}, \eqref{eq:f_11}, and \eqref{eq:f_22} shown in Figures \ref{fig:ph1 (u1 neq u2) Lf, Tf, and Mf} and \ref{fig:ph2 (u1 neq u2) Lf, Tf, and Mf} for the smooth and non-smooth phantoms, respectively. Corresponding reconstructions in the presence of various levels of noise are presented in Figures \ref{fig:ph1 (u1 neq u2) Lf, Tf, and Mf with noise} and \ref{fig:ph2 (u1 neq u2) Lf, Tf, and Mf with noise}.
\begin{figure}[H]
     \centering
     \begin{subfigure}[b]{0.5\textwidth}
         \centering
         \includegraphics[width=\textwidth]{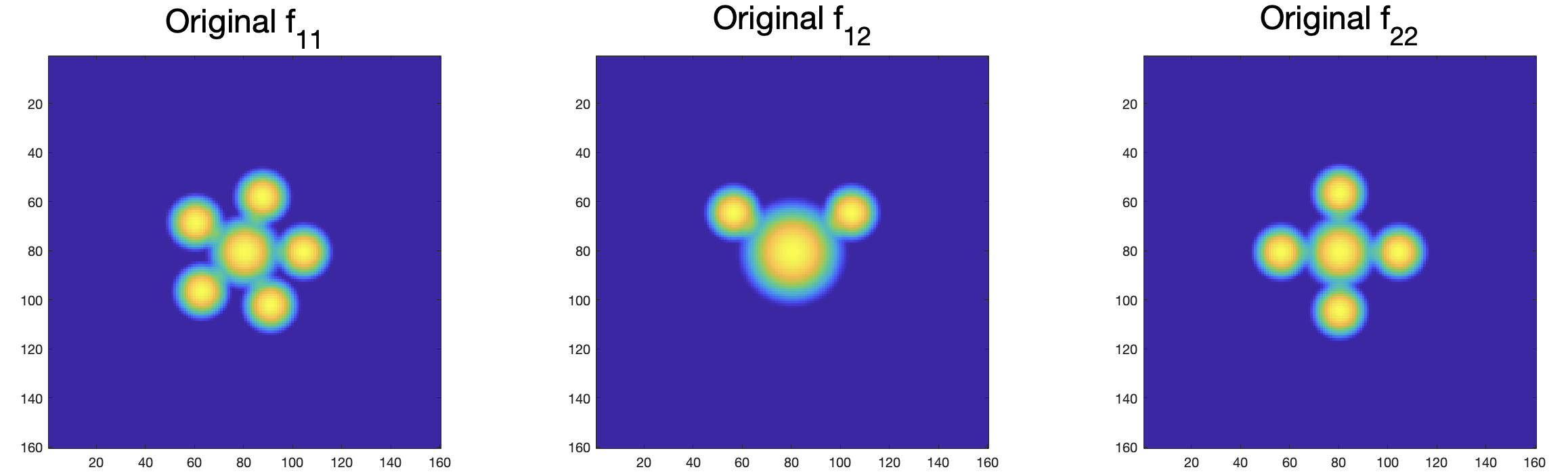}
     \end{subfigure}
     \vfill
     \begin{subfigure}[b]{0.5\textwidth}
         \centering
         \includegraphics[width=\textwidth]{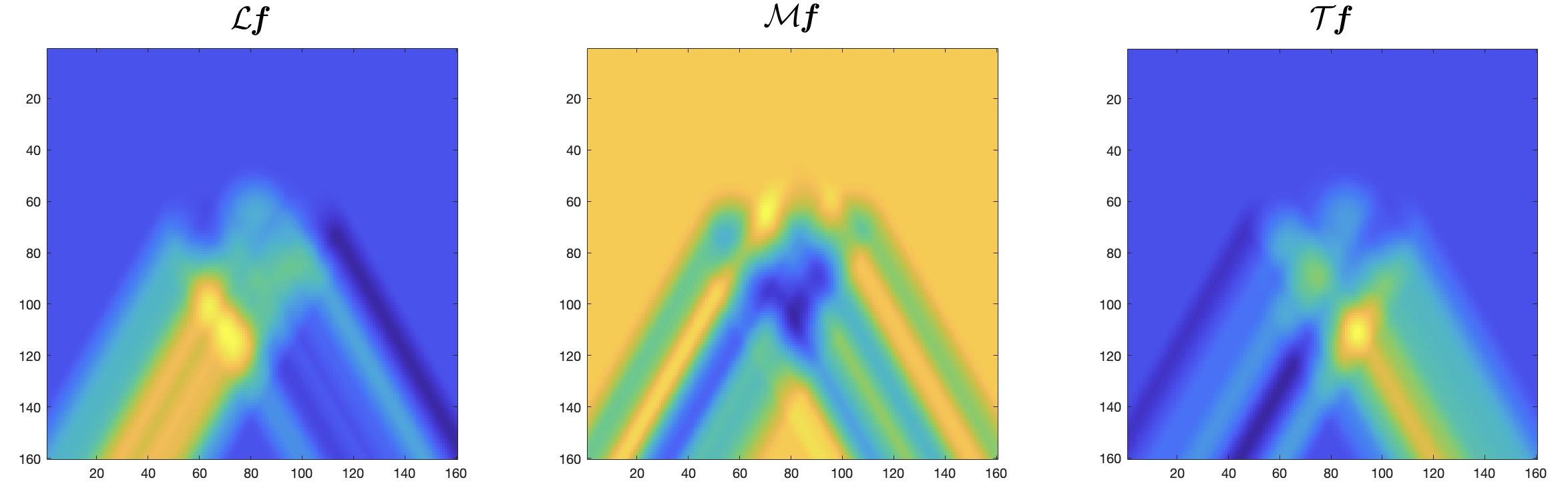}
     \end{subfigure}
     \vfill
     \begin{subfigure}[b]{0.5\textwidth}
         \centering
         \includegraphics[width=\textwidth]{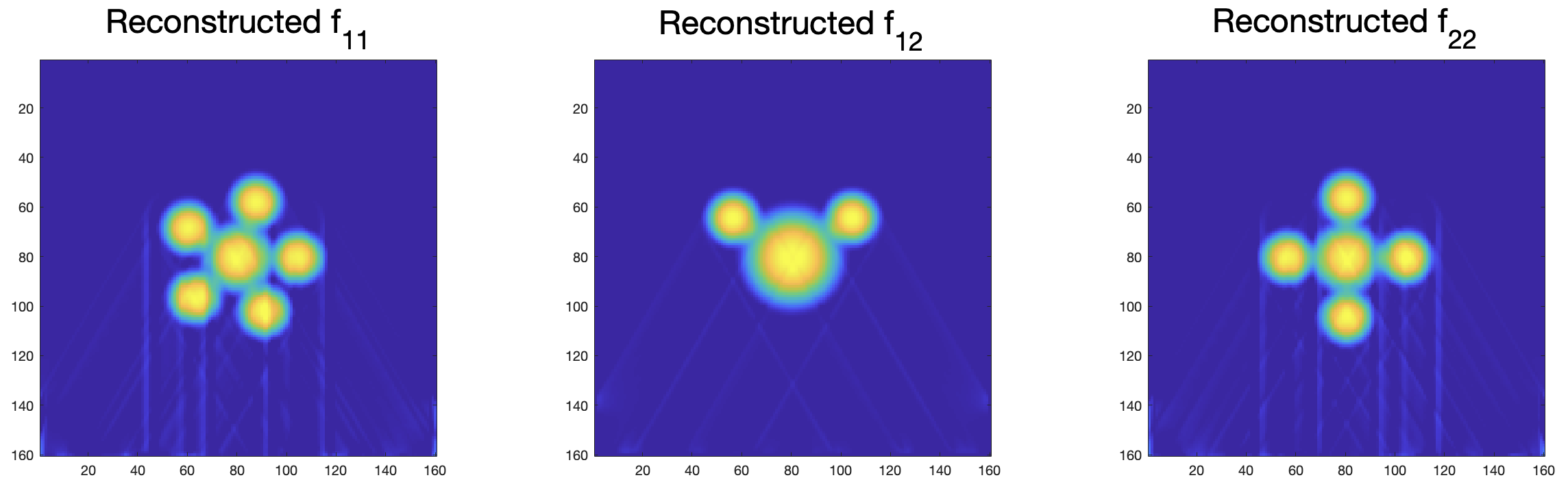}
     \end{subfigure}
\caption{Reconstructions of $\vf$ from $\Lc\vf,\Mc\vf,\Tc\vf$, $\vu=(\cos{\pi/3},\sin{\pi/3})$, using formulas \eqref{eq:elliptic equation for f12}, \eqref{eq:f_11}, \eqref{eq:f_22}.}\label{fig:ph1 (u1 neq u2) Lf, Tf, and Mf}
\end{figure}
 \begin{figure}[H]
     \centering
   \includegraphics[width= 0.83\textwidth]{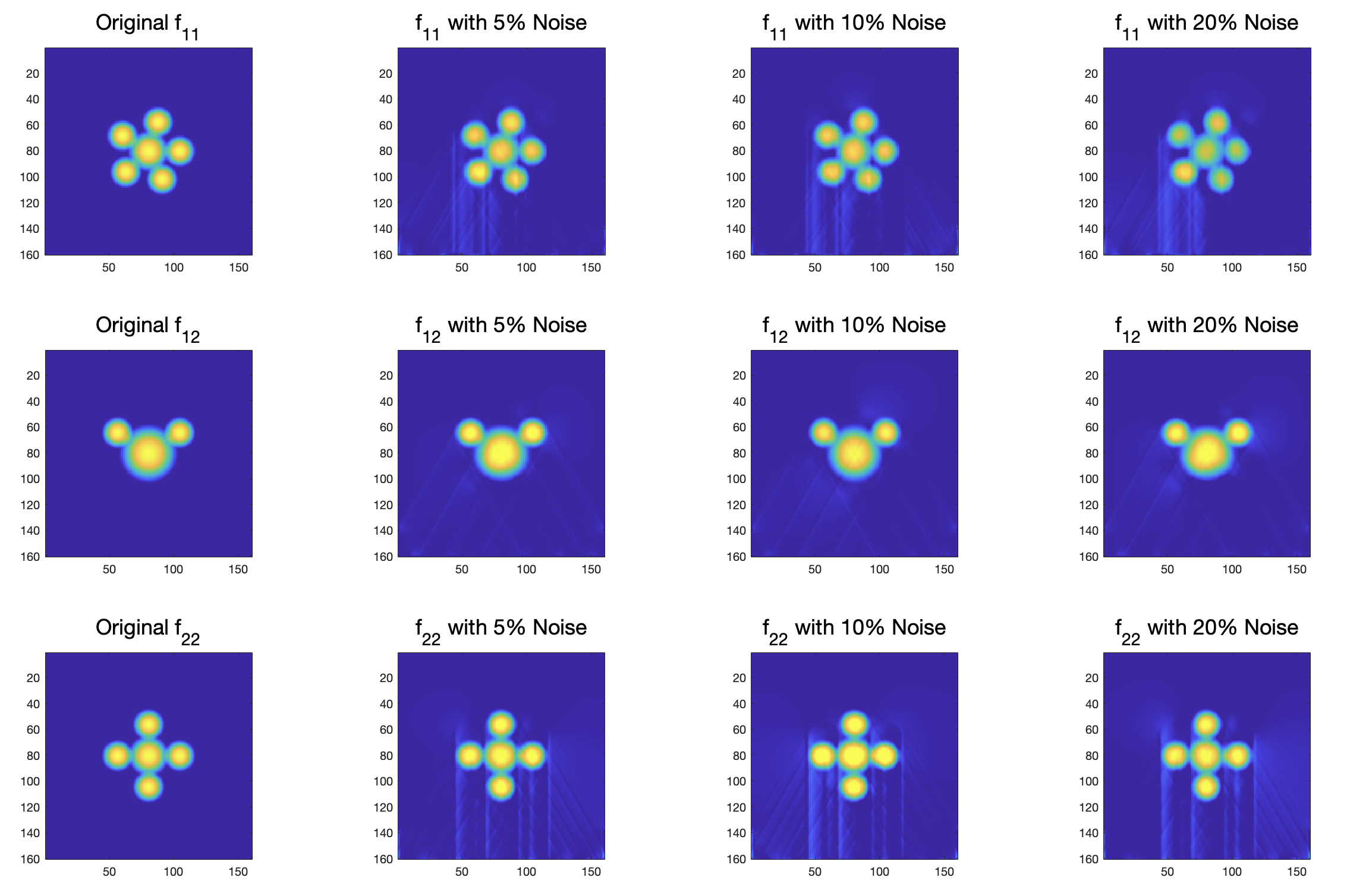}
    \caption{Reconstructions from $\Lc\vf,\Mc\vf,\Tc\vf$, $\vu=(\cos{\pi/3},\sin{\pi/3})$,  $5\%, 10\%$ and $20\%$ noise.}\label{fig:ph1 (u1 neq u2) Lf, Tf, and Mf with noise}
 \end{figure}
\begin{table}[h!]
 \begin{center}
\begin{tabular}{ |c|c|c|c|c| } 
 \hline
  $\vf$  & No noise & $5\%$ Noise & $10\%$ Noise& $20\%$ Noise \\ 
 \hline
  $f_{11}$  & 8.49\% & 16.27\% &18.60\% & 30.64\%\\ 
 \hline
 $f_{12}$ & 1.84\% & 8.62\% & 6.45\% &  9.54\%\\  
 \hline
 $f_{22}$ & 8.77\% & 15.20\% & 25.15\% & 21.34\%\\  
 \hline
\end{tabular}
\caption{Relative errors of reconstructing the components of the smooth Phantom 1 using $\Lc\vf, \Mc\vf, \Tc\vf$, and $\vu=(\cos{\pi/3,\sin{\pi/3}})$.}
\end{center}
\end{table}

 
\begin{table}[h!]
\begin{center}
\begin{tabular}{ |c|c|c|c|c| } 
 \hline
  $\vf$  & No noise & $5\%$ Noise & $10\%$ Noise& $20\%$ Noise \\ 
 \hline
 $f_{11}$  & 321.77\% &293.14\%  & 338.85\% & 450.29\% \\ 
 \hline
 $f_{12}$ & 16.26\% & 16.29\% & 18.79\% & 24.92\%\\  
 \hline
 $f_{22}$ & 234.58\% & 262.40\%  & 273.35\% & 257.46\%\\  
 \hline
\end{tabular}
\caption{Relative errors of reconstructing the components of the non-smooth Phantom 2, using  $\Lc\vf, \Mc\vf, \Tc\vf$, and $\vu=(\cos{\pi/3, \sin{\pi/3}})$.} \label{tab:f_from_LMT-2}
\end{center}
\end{table}

In our inversion formulas stated in Theorem \ref{eq:full-rec from L,T and M}, many terms include successive applications of differentiation and integration operators. It is easy to check that these operators commute, and one can choose the order in which they are applied to the data. In all results presented in this paper, the computations of terms like that are performed by applying the differentiation(s) first, followed by the integration. We have also numerically implemented some of these inversion formulas using the opposite order of the operators. However, the reconstructions using the latter approach were not much different, so we have not included them in the article. 

The alternatives discussed above are common in numerical inversions of Radon type transforms. Similar situations appear, for example, in the selection between filtered and $\rho$-filtered backprojection formulas for inversion of the (classical) Radon transform (e.g. see \cite{herman2009fundamentals}) and the spherical Radon transform (e.g. see \cite{ambartsoumian2007thermoacoustic}). The systematic analysis and comparison of the two approaches are beyond the scope of this work and are subject of future research. \\

The artifacts appearing in the reconstructions in Figures \ref{fig:ph2 (u1 = u2) Lf, Tf, and Mf with noise}, \ref{fig:ph2 (u1 neq u2) Lf, Tf, and Mf}, \ref{fig:ph2 (u1 neq u2) Lf, Tf, and Mf with noise} can be explained just like those in the reconstruction of $g_1$ in Figure \ref{fig:ph2 dg (parabolic)}. Namely, they are caused by the errors of numerical differentiation, propagated by integration with respect to $x_1$ (horizontal artifacts), or with respect to $x_2$ (vertical artifacts). Although visually they do not seem to alter the original image too much, their strengths are quite large, substantially increasing the relative errors of the reconstructions (see Tables \ref{tab:f_from_LMT-1}, \ref{tab:f_from_LMT-2}). It is also worth pointing out, that the reconstructions of the smooth phantoms using the same algorithms have much milder artifacts and much smaller relative errors.
 \vspace{10cm}
\begin{figure}[H]
     \centering
     \begin{subfigure}[b]{0.58\textwidth}
         \centering
  \includegraphics[width=\textwidth]{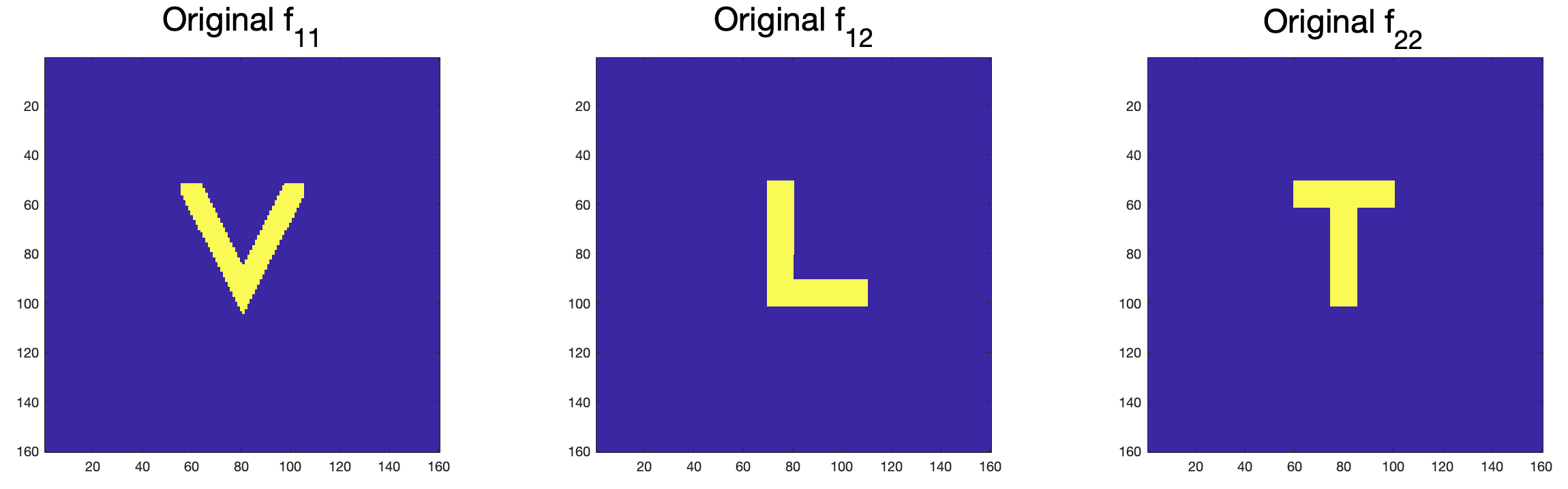}
     \end{subfigure}
     \vfill
     \begin{subfigure}[b]{0.58\textwidth}
         \centering
\includegraphics[width=\textwidth]{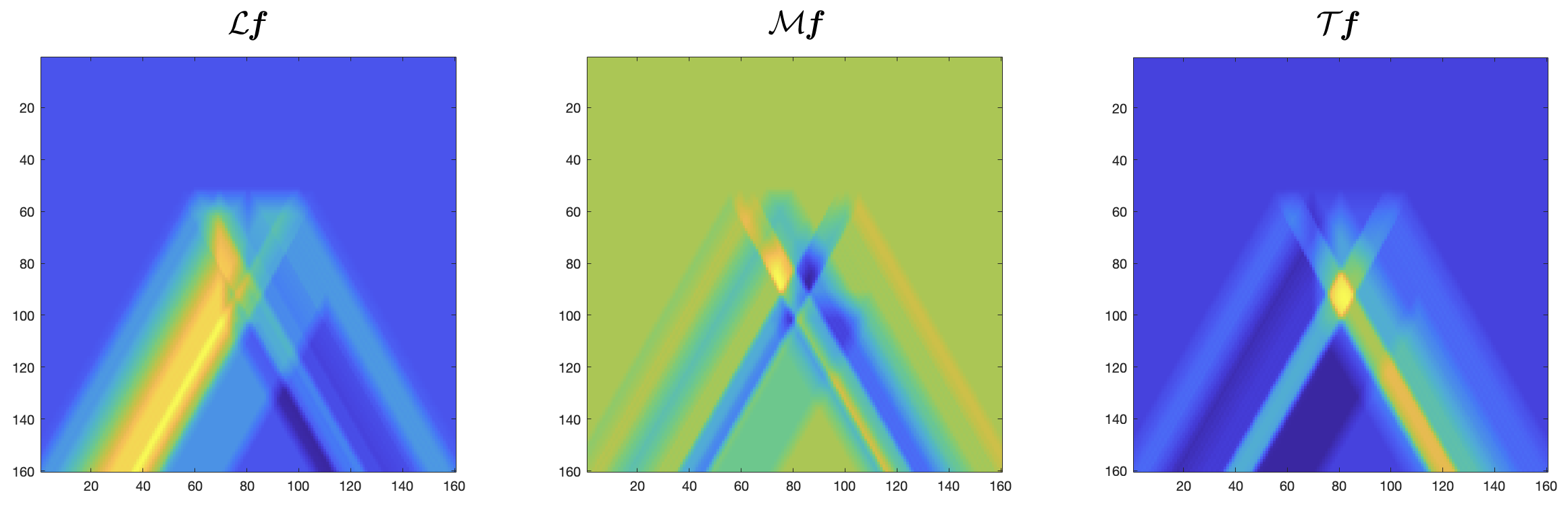}
     \end{subfigure}
     \vfill
     \begin{subfigure}[b]{0.60\textwidth}
         \centering
         \includegraphics[width=\textwidth]{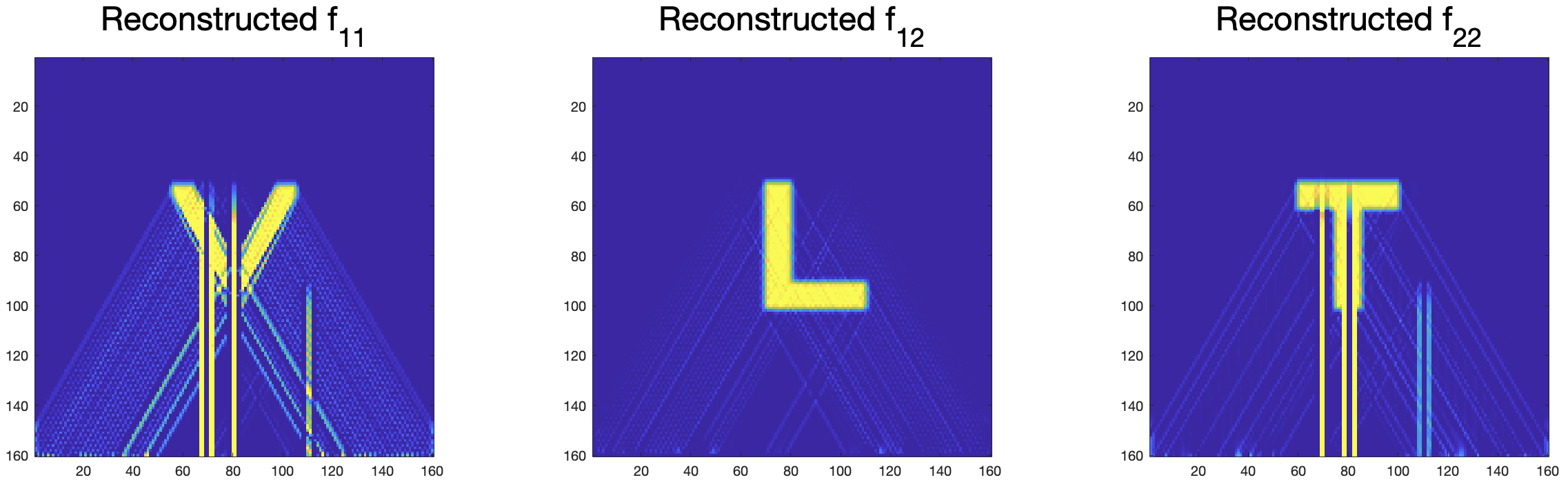}
\end{subfigure}
\caption{Reconstructions of $\vf$ from $\Lc\vf,\Mc\vf,\Tc\vf$, $\vu=(\cos{\pi/3,\sin{\pi/3}})$, using formulas \eqref{eq:elliptic equation for f12}, \eqref{eq:f_11}, \eqref{eq:f_22}.}\label{fig:ph2 (u1 neq u2) Lf, Tf, and Mf}
\end{figure}

 \begin{figure}[H]
     \centering
   \includegraphics[width= 0.89\textwidth]{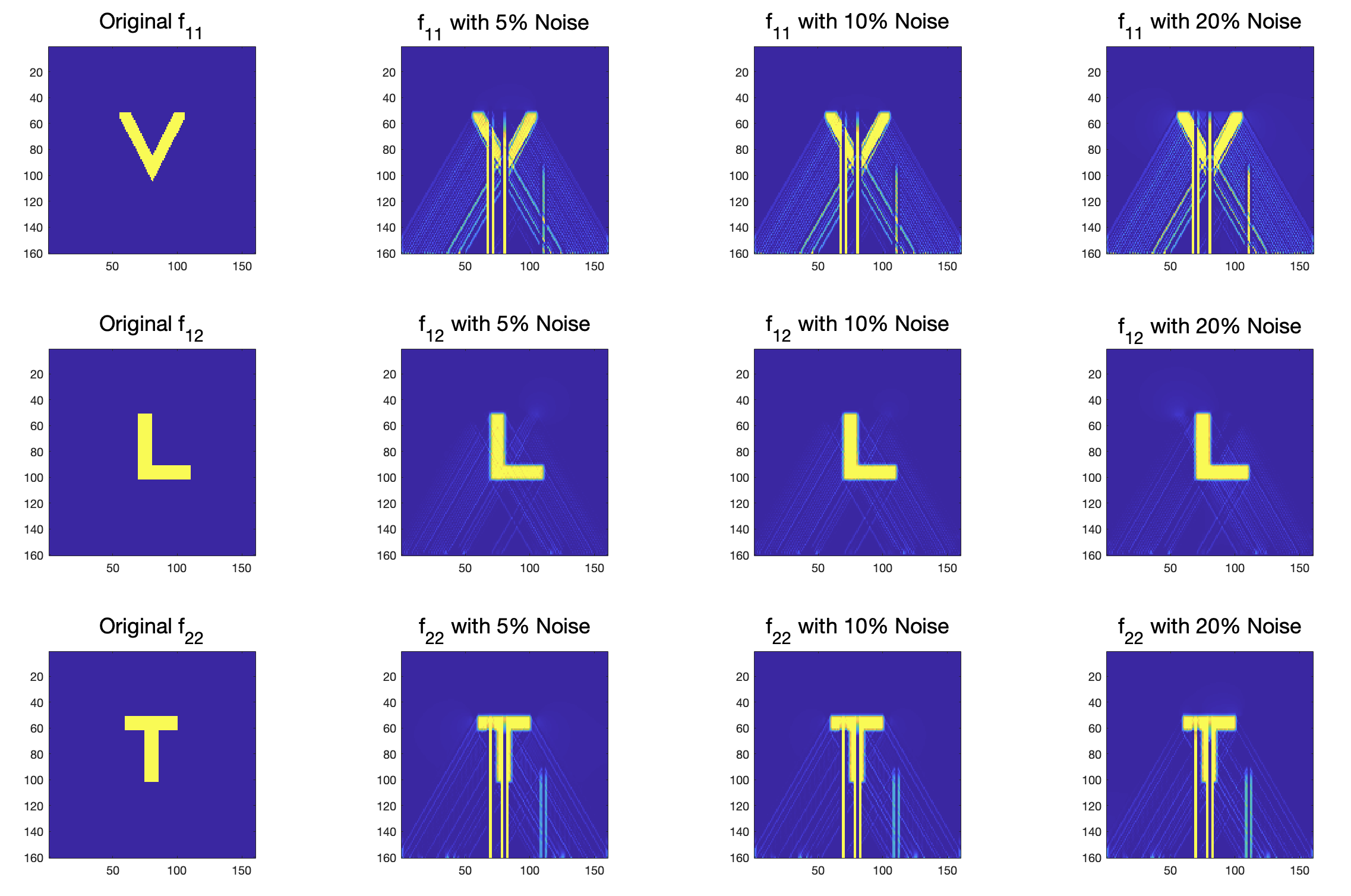}
    \caption{Reconstructions from $\Lc\vf,\Mc\vf,\Tc\vf$, $\vu=(\cos{\pi/3,\sin{\pi/3}})$,  $5\%, 10\%$ and $20\%$ noise.}\label{fig:ph2 (u1 neq u2) Lf, Tf, and Mf with noise}
 \end{figure}

\subsection{Recovery of tensor fields from first moment VLTs} \label{subsec: full recovery with moments}

This subsection is devoted to the full recovery of $\vf$ from various combinations of its V-line transforms and their first moments. The reconstruction of the unknown tensor field can be obtained from different subsets of transforms involving first moments, such as 
$\{\Lc\vf,\Lc^1\vf,\Tc\vf\}$, $\{\Tc\vf,\Tc^1\vf,\Lc\vf\}$, $\{\Mc\vf,\Mc^1\vf,\Lc\vf\}$, $\{\Mc\vf,\Mc^1\vf,\Tc\vf\}$, $\{\Lc\vf,\Lc^1\vf,\Mc\vf\}$, or $\{\Tc\vf,\Tc^1\vf,\Mc\vf\}$. Here, we discuss only two subsets $\{\Lc\vf,\Lc^1\vf,\Tc\vf\}$ and $\{\Lc\vf,\Lc^1\vf,\Mc\vf\}$, as the other combinations yield results similar to one of these cases. More specifically, we discuss the recovery of $\vf$ from $\{\Lc\vf,\Lc^1\vf,\Tc\vf\}$ when $\vu=(\cos{\pi/3},\sin{\pi/3})$, and from $\{\Lc\vf,\Lc^1\vf,\Mc\vf\}$ when $\vu=(\cos{\pi/4},\sin{\pi/4}).$

\begin{thr}\label{eq:full-rec from L,L^1,T}
If $u_1\ne u_2$, then $\vf\in C_c^2\left(S^2;D_1\right)$  
can be recovered from $\Lc \vf$, $\Lc^1 \vf$, and $\Tc \vf$ by
\begin{align}
f_{12}(\vx) &= \frac{1}{4u_1^2} \partial_{x_2}\Xc_{\vev_1}\left[D_\vu D_\vv \Lc^1\vf+\left( D_\vu + D_\vv + \frac{1}{u_2}D_\vu D_\vv \Xc_{\vev_2} \right)\Lc\vf \, \right](\vx),\label{eq: f12 from moments}\\
f_{11}(\vx)&=\frac{1}{2(u_2^2-u_1^2)}\left[ D_\vu D_\vv\Lc^1\vf  +  (D_\vu + D_\vv)\Lc\vf  +u_2 D_\vu D_\vv  \Xc_{\vev_2}\left(\Lc \vf + \Tc \vf \right)\,\right](\vx),\label{eq: f11 from moments}\\
f_{22}(\vx)&=\frac{1}{2(u_1^2-u_2^2)}\left[D_\vu D_\vv\Lc^1\vf  + (D_\vu + D_\vv)\Lc\vf  +\frac{u_1^2}{u_2} D_\vu D_\vv  \Xc_{\vev_2}\left(\Lc \vf + \Tc \vf \right)\,\right](\vx).\label{eq: f22 from moments}
 \end{align}
\end{thr}
\begin{thr}\label{eq:full rec from L, L^1, M}
The tensor field $\vf\in C_c^2\left(S^2;D_1\right)$  can be recovered from $\Lc \vf$, $\Lc^1 \vf$, and $\Mc\vf$ by
\begin{align}
f_{12}(\vx) &= \frac{1}{4u_1^2} \partial_{x_2}\Xc_{\vev_1}\left[D_\vu D_\vv \Lc^1\vf+\left( D_\vu + D_\vv + \frac{1}{u_2}D_\vu D_\vv \Xc_{\vev_2} \right)\Lc\vf \, \right](\vx),\label{eq: f12 from moments L, L^1}\\
f_{11}(\vx) &=  - \frac{1}{2} \left\{D_\vu D_\vv\Lc^1\vf  +  (D_\vu + D_\vv)\Lc\vf - \frac{u_2}{u_1^2}D_\vu D_\vv \Xc_{-\vev_1} \left[(u_1^2 -u_2^2)\Vc(f_{12}) - \Mc \vf\, \right]\right\}(\vx),\label{eq: f11 from moments L, L^1}\\
f_{22}(\vx) &=  - \frac{1}{2} \left\{D_\vu D_\vv\Lc^1\vf  +  (D_\vu + D_\vv)\Lc\vf + \frac{1}{u_2}D_\vu D_\vv \Xc_{-\vev_1} \left[(u_1^2 -u_2^2)\Vc(f_{12}) - \Mc \vf\,\right]\right\}(\vx).\label{eq: f22 from moments L, L^1}
 \end{align}
\end{thr}

\begin{rem}
As we have already mentioned in Remark \ref{rem:5-pixels}, the compactness of the support of $\vf$ has certain implications about the support of the VLTs of $\vf$ and their derivatives. In particular, it is easy to check that outside of a disc contained in the square domain of reconstruction, the terms $D_\vu D_\vv\Lc^1\vf  +  (D_\vu + D_\vv)\Lc\vf$ and $D_\vu D_\vv \Xc_{\vev_2}\Lc\vf$ are zero. The same is true for the expressions containing the transverse and mixed transforms in place of the longitudinal VLTs. In our numerical implementations, we assign the value of 0 to all entries of those quantities within 5 pixels from the boundary of the data domain. This reduces the artifacts in the reconstruction due to the errors of numerical differentiation at the boundary of the domain.
\end{rem}

 In the case of inversion algorithms applicable to VLTs with different opening angles, the quality of reconstruction often varies based on the opening angle (e.g. see \cite{Gaik_Mohammad_Rohit_2024_numerics}). In particular, the artifacts decrease when the opening angle gets closer to $\pi/2$, which corresponds to the case $\vu= (\cos{\pi/4},\sin{\pi/4})$. This fact is one of the primary reasons of the discrepancy in the quality of reconstructions between Figures \ref{fig:ph1 from L, L1, T}-\ref{fig:ph2 from L, L1, T with noise} and Figures
\ref{fig:ph1 from L, L1, M}-\ref{fig:ph2 from L, L1, M with noise}.
 
Figures \ref{fig:ph1 from L, L1, T} and \ref{fig:ph2 from L, L1, T} are generated by implementing formulas \eqref{eq: f11 from moments}, \eqref{eq: f12 from moments}, \eqref{eq: f22 from moments}, respectively, for Phantom 1 and 2. Figures \ref{fig:ph1 from L, L1, T with noise} and \ref{fig:ph2 from L, L1, T with noise} show the effect of noise on the reconstructions. 
\begin{figure}[H]
\centering
\includegraphics[width= 0.82\textwidth]{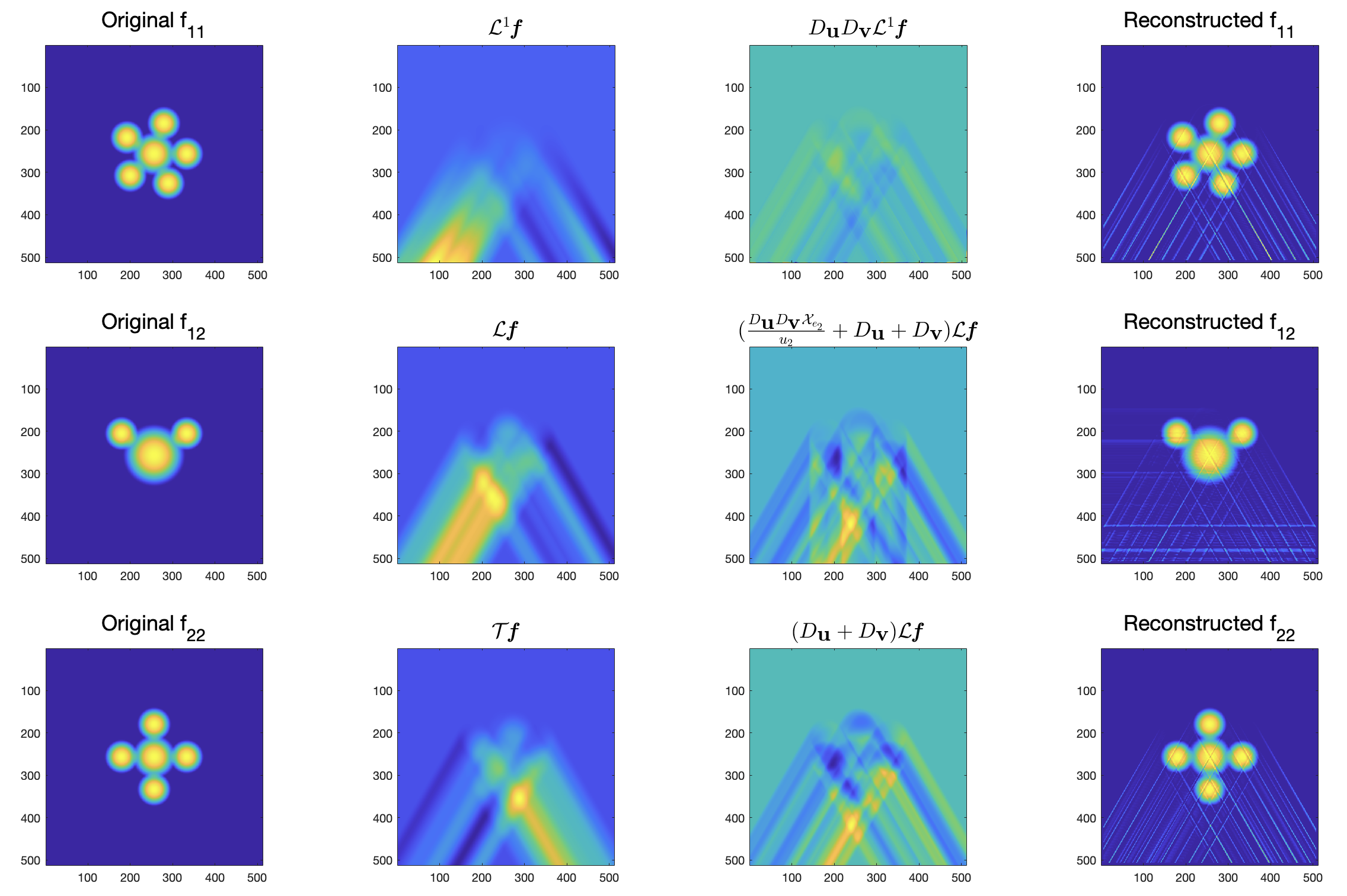}
\caption{Reconstruction of $\vf$ from $\Lc\vf,\Lc^1\vf,\Tc\vf$, $\vu= (\cos{\pi/3},\sin{\pi/3})$ using formulas \eqref{eq: f11 from moments}, \eqref{eq: f12 from moments}, \eqref{eq: f22 from moments}.}\label{fig:ph1 from L, L1, T}
 \end{figure}

\begin{figure}[H]
     \centering
   \includegraphics[width= 0.82\textwidth]{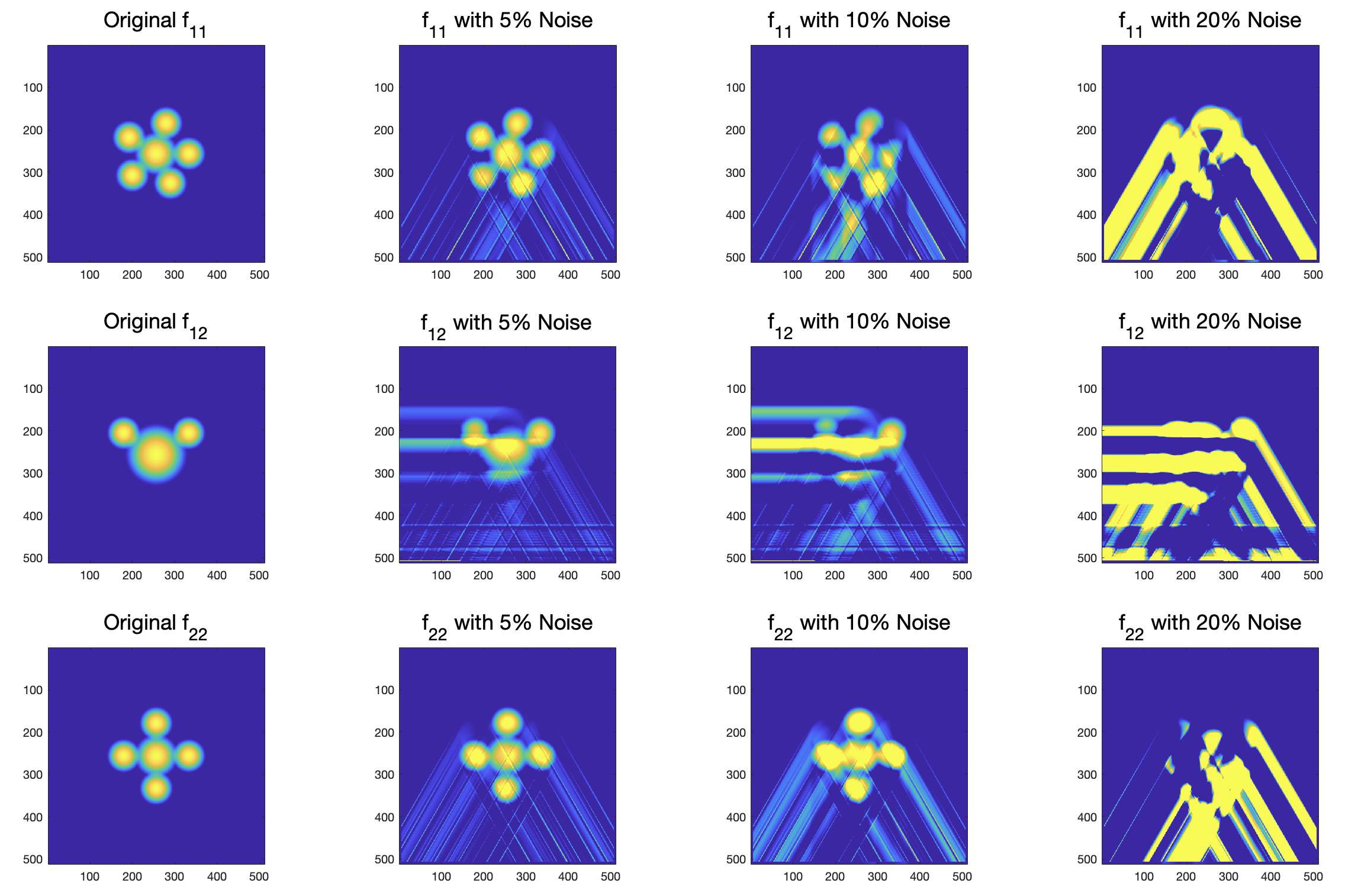}
    \caption{ Reconstruction of $\vf$ from $\Lc\vf,\Lc^1\vf,\Tc\vf$ when $\vu=(\cos{\pi/3,\sin{\pi/3}})$ with $5\%, 10\%$ and $20\%$ noise.} \label{fig:ph1 from L, L1, T with noise}
 \end{figure}
\begin{figure}[H]
\centering
\includegraphics[width= 0.82\textwidth]{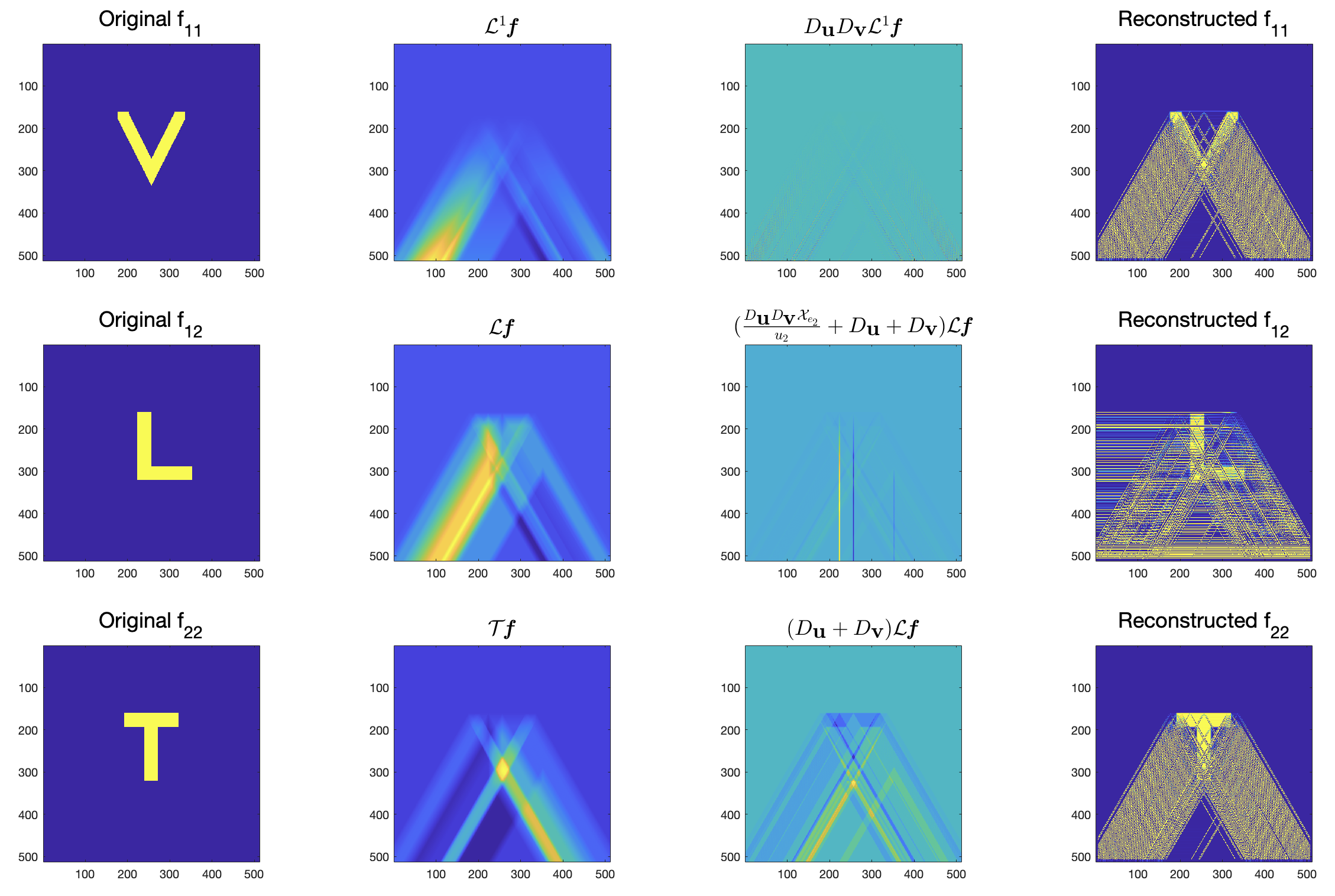}
\caption{
Reconstruction of $\vf$ from $\Lc\vf,\Lc^1\vf,\Tc\vf$, $\vu=(\cos{\pi/3,\sin{\pi/3}})$ 
using formulas \eqref{eq: f11 from moments},\eqref{eq: f12 from moments},\eqref{eq: f22 from moments}.}\label{fig:ph2 from L, L1, T}
\end{figure}

\begin{figure}[H]
     \centering
   \includegraphics[width= 0.82\textwidth]{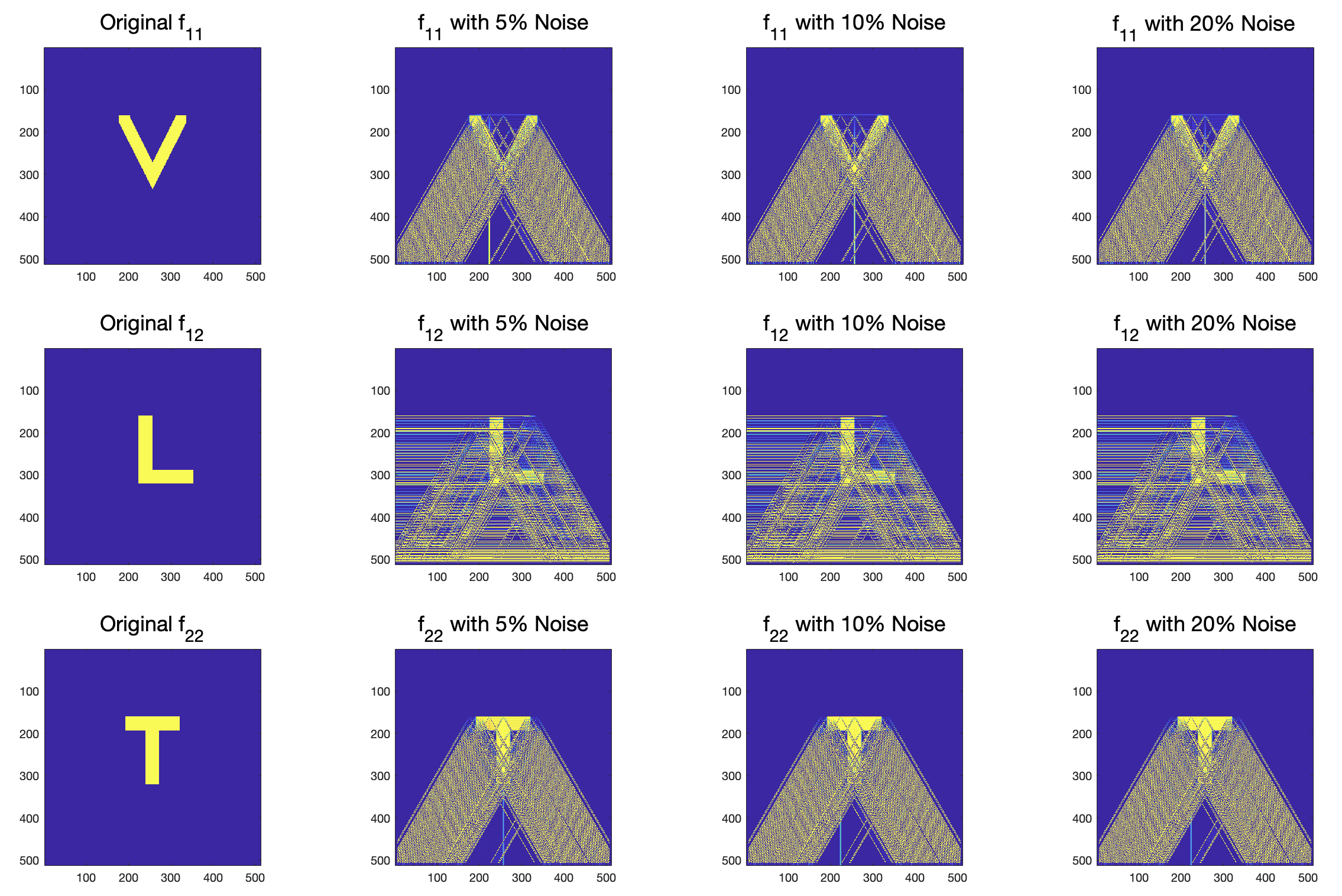}
    \caption{Reconstruction of $\vf$ from $\Lc\vf,\Lc^1\vf,\Tc\vf$ when $\vu=(\cos{\pi/3,\sin{\pi/3}})$ with $5\%, 10\%$ and $20\%$ noise.}\label{fig:ph2 from L, L1, T with noise}
 \end{figure}

\begin{table}[h!]
 \begin{center}
\begin{tabular}{ |c|c|c|c|c| } 
 \hline
  $\vf$  & No noise & $5\%$ Noise & $10\%$ Noise& $20\%$ Noise \\ 
 \hline
  $f_{11}$  & 13.77\% & 17.29\% & 61.49\% & 428.98\%\\ 
 \hline
 $f_{12}$ &806.32\% & 833.58\% & 990.84\% &  2404.38\%\\  
 \hline
 $f_{22}$ & 14.10\% & 18.43\% & 63.49\% & 445.50\%\\  
 \hline
\end{tabular}
\caption{Relative errors of reconstructing Phantom 1 from $\Lc\vf, \Lc^1\vf$, $\Tc\vf$, $\vu=(\cos{\pi/3,\sin{\pi/3}})$.}\label{tab:L_L1_T-1}
\end{center}
\end{table}

\begin{table}[h!]
\begin{center}
\begin{tabular}{ |c|c|c|c|c| } 
 \hline
  $\vf$  & No noise & $5\%$ Noise & $10\%$ Noise& $20\%$ Noise \\ 
 \hline
  $f_{11}$  & 4009.19\% & 4036.22\% & 4215.82\% & 4218.84\%\\ 
 \hline
 $f_{12}$ & 2408.93\% & 2422.88\% & 2534.95\% &  2532.71\%\\  
 \hline
 $f_{22}$ & 3041.98 \% & 3062.53\% & 3198.72\% & 3201.04\%\\  
 \hline
\end{tabular}
\caption{Relative errors of reconstructing Phantom 2 from $\Lc\vf, \Lc^1\vf$, $\Tc\vf$, $\vu=(\cos{\pi/3,\sin{\pi/3}})$.} \label{tab:L_L1_T-2}
\end{center}
\end{table}
The inversion algorithms using the first moment transforms are highly ill-conditioned. The primary source of instability of the associated inversion formulas lies in the presence of an additional differentiation in the terms involving the first moment (e.g. see formulas \eqref{eq: f12 from moments}-\eqref{eq: f22 from moments L, L^1}). As a result, the reconstructions based on these approaches have strong artifacts, which become worse when applied to non-smooth phantoms (compare Figures \ref{fig:ph1 from L, L1, T},  \ref{fig:ph1 from L, L1, T with noise} with Figures \ref{fig:ph2 from L, L1, T}, \ref{fig:ph2 from L, L1, T with noise}). The situation is more severe in the setups where the V-lines have an opening angle far from $\pi/2$ (compare Figures \ref{fig:ph1 from L, L1, T}-\ref{fig:ph2 from L, L1, T with noise} with Figures \ref{fig:ph1 from L, L1, M}-\ref{fig:ph2 from L, L1, M with noise}). 

Although visually the patterns of the artifacts in Figure \ref{fig:ph2 from L, L1, T with noise} do not seem to change much with the addition of noise to the data, the strengths of the artifacts change, which can be noticed from Table \ref{tab:L_L1_T-2} of relative errors. 

A similar phenomenon, related to the ill-conditioned inversions involving first moment VLTs due to the extra derivatives, was observed in the study of VLTs defined on vector fields in $\mathbb{R}^2$ (see \cite{Gaik_Mohammad_Rohit_2024_numerics}).
\vspace{1cm}

Next, we consider the recovery of $\vf$ from the combination $\Lc\vf,\Lc^1\vf$ and $\Mc\vf.$ In Figures \ref{fig:ph1 from L, L1, M} and \ref{fig:ph2 from L, L1, M}, we use the inversion formulas \eqref{eq: f12 from moments L, L^1}, \eqref{eq: f11 from moments L, L^1}, and \eqref{eq: f22 from moments L, L^1} to generate $f_{12}$, $f_{11}$, and $f_{22}$ respectively for the smooth and non-smooth phantoms. Figures \ref{fig:ph1 from L, L1, M with noise} and \ref{fig:ph2 from L, L1, M with noise} show the reconstructions after adding different levels of noise for smooth and non-smooth phantoms, respectively. 

\vspace{10cm}
\begin{figure}[H]
\centering
\includegraphics[width= 0.83\textwidth]{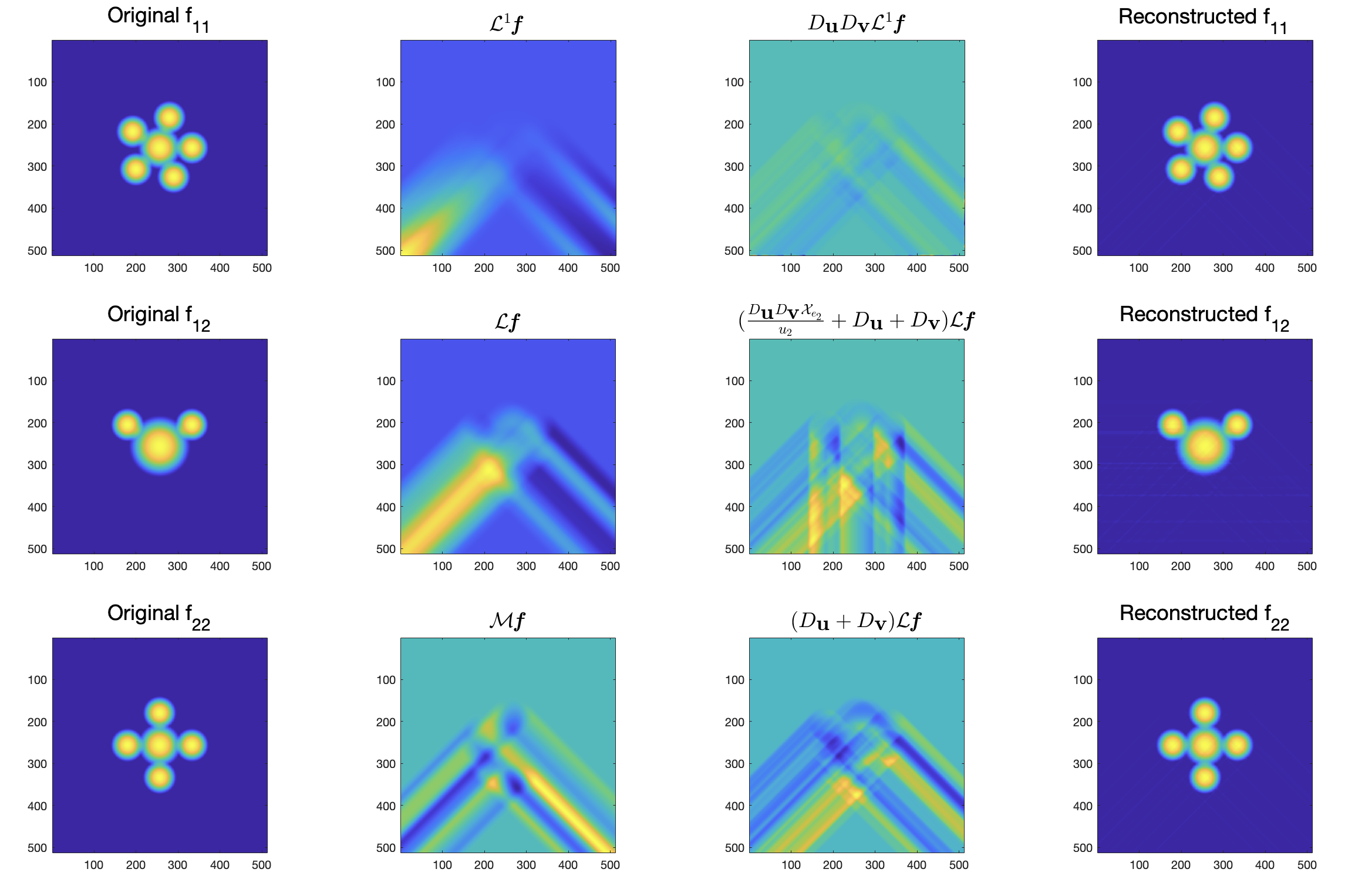}
\caption{Reconstruction of $\vf$ from $\Lc\vf,\Lc^1\vf,\Mc\vf$, $\vu= (\cos{\pi/4},\sin{\pi/4})$ using formulas \eqref{eq: f12 from moments L, L^1}, \eqref{eq: f11 from moments L, L^1}, \eqref{eq: f22 from moments L, L^1}.}\label{fig:ph1 from L, L1, M}
\end{figure}

\begin{figure}[H]
     \centering
   \includegraphics[width= 0.83\textwidth]{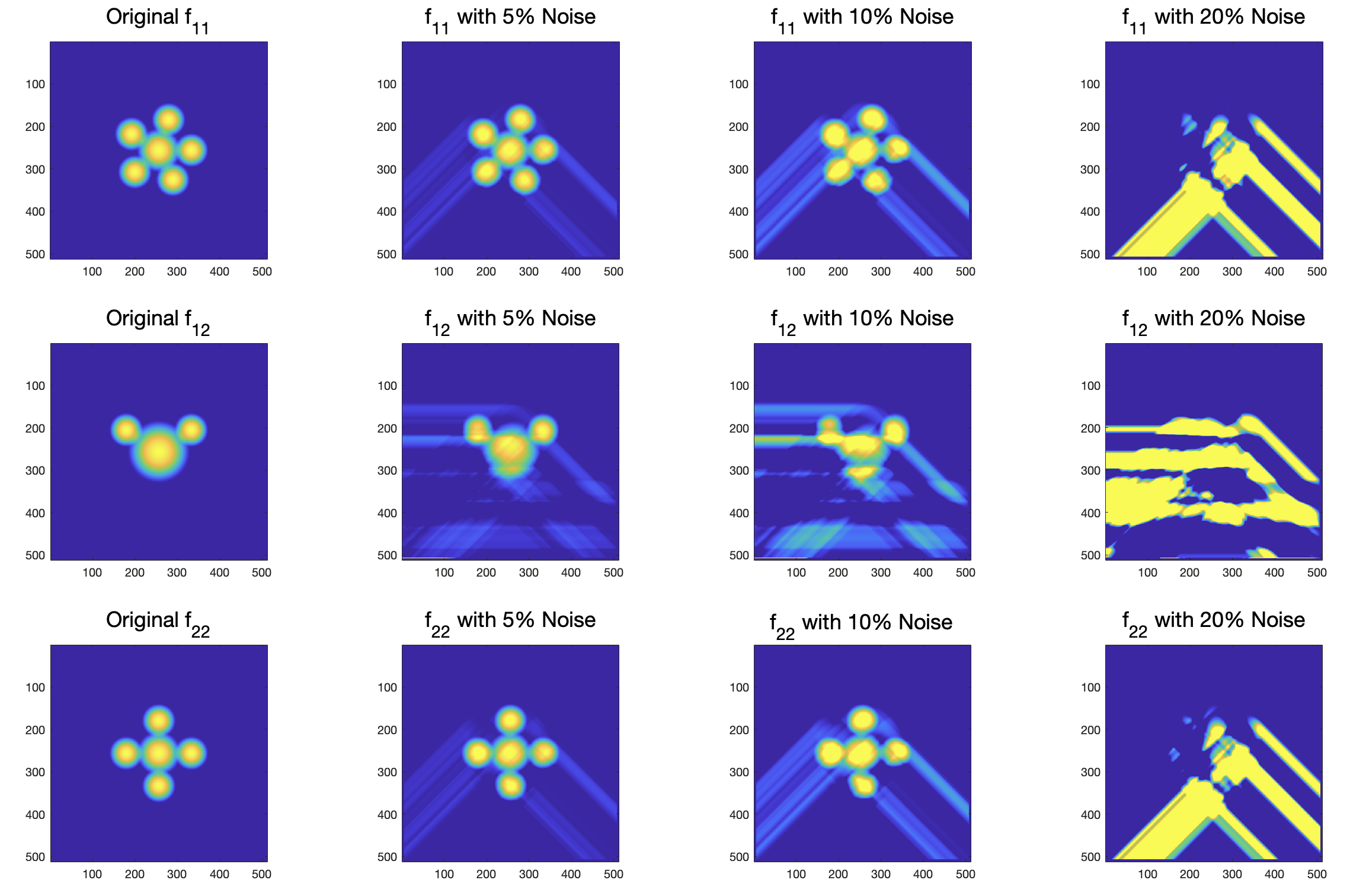}
    \caption{Reconstruction of $\vf$ from $\Lc\vf,\Lc^1\vf,\Mc\vf$, $\vu=(\cos{\pi/4,\sin{\pi/4}})$ with $5\%, 10\%$ and $20\%$ noise.}\label{fig:ph1 from L, L1, M with noise}
 \end{figure}

\begin{figure}[H]
\centering
\includegraphics[width= 0.87\textwidth]{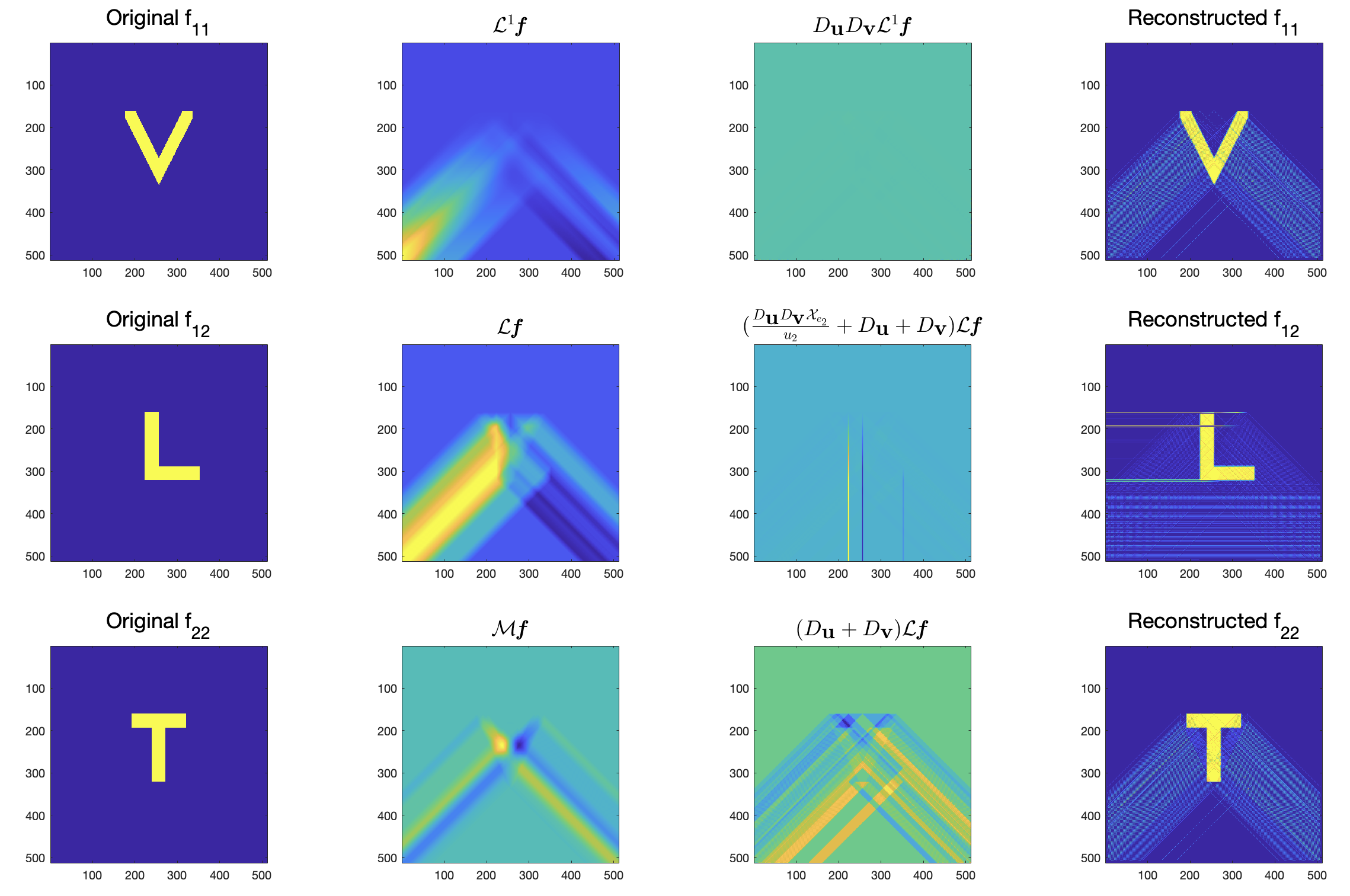}
\caption{ 
Reconstructions of $\vf$ from $\Lc\vf, \Lc^1\vf, \Mc\vf$, $\vu=(\cos{\pi/4,\sin{\pi/4}})$,
using formulas \eqref{eq: f12 from moments L, L^1}, \eqref{eq: f11 from moments L, L^1}, \eqref{eq: f22 from moments L, L^1}.}\label{fig:ph2 from L, L1, M}
\end{figure}

\begin{figure}[H]
     \centering
   \includegraphics[width= 0.87\textwidth]{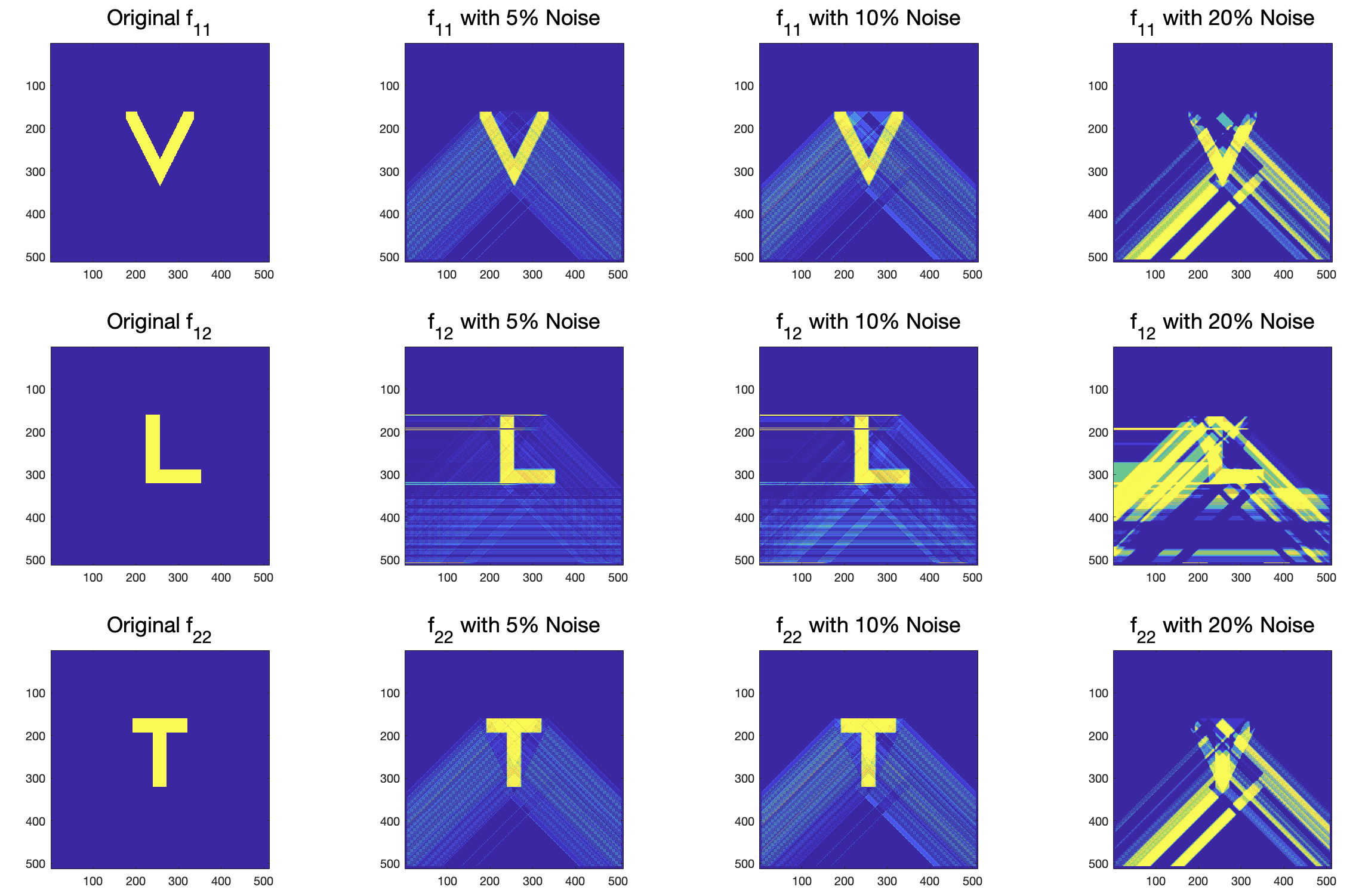}
    \caption{Reconstructions from $\Lc\vf,\Lc^1\vf,\Mc\vf$, $\vu=(\cos{\pi/4,\sin{\pi/4}})$,  $5\%, 10\%$ and $20\%$ noise.}\label{fig:ph2 from L, L1, M with noise}
 \end{figure}

 \begin{table}[h!]
\begin{center}
\begin{tabular}{ |c|c|c|c|c| } 
 \hline
 $\vf$  & No noise & $5\%$ Noise & $10\%$ Noise& $20\%$ Noise \\ 
 \hline
  $f_{11}$  & 0.68\% & 11.72\% & 37.44\% & 273.99\%\\ 
 \hline
 $f_{12}$ & 803.03\% & 840.92\% & 1008.68\% &  1325.82\%\\  
 \hline
 $f_{22}$ & 0.70\% & 11.70\% & 38.17\%  & 278.72\%   \\  
 \hline
\end{tabular}
\caption{Relative errors of reconstructing Phantom 1 from $\Lc\vf, \Lc^1\vf$, $\Mc\vf$, $\vu=(\cos{\pi/4,\sin{\pi/4}})$.}
\end{center}
\end{table}

 \begin{table}[h!]
\begin{center}
\begin{tabular}{ |c|c|c|c|c| } 
 \hline
 $\vf$  & No noise & $5\%$ Noise & $10\%$ Noise& $20\%$ Noise \\ 
 \hline
  $f_{11}$  & 81.83\% & 82.69\% & 92.73\% & 204.12\%\\ 
 \hline
 $f_{12}$ & 857.81\% & 883.65\% & 1026.22\% &  1378.59\%\\  
 \hline
 $f_{22}$ & 62.14\% & 62.80\% & 70.42\% & 156.51\%\\  
 \hline
\end{tabular}
\caption{Relative errors of reconstructing Phantom 2 from $\Lc\vf, \Lc^1\vf$, $\Mc\vf$, $\vu=(\cos{\pi/4,\sin{\pi/4}})$.}
\end{center}
\end{table}

\vspace{-1cm}

\subsection{Numerical implementation of inverting the tensor star transform}\label{subsec: star transform}
Before starting our discussion about the star transform of symmetric 2-tensor fields, let us quickly recall the definition of the star transform. Assume $\vf\in C_c^2\left(S^2;\Rb^2\right)$,
and let $\vgamma_{1},\vgamma_{2},...,\vgamma_{m}$ be distinct unit vectors in $\mathbb{R}^{2}$. The star transform of $\vf$ is defined as
\begin{align*}
\Sc\vf &=  \sum_{i=1}^{m} c_{i} {\Xc}_{{\vgamma}_{i}} 
\begin{bmatrix}
 \vf\cdot {{\vgamma}_{i}}^{2}\\
 \vf\cdot {{\vgamma}_{i} \odot {\vgamma}_{i}^{\perp} }\\
 \vf\cdot ({\vgamma}_{i}^{\perp})^{2}
\end{bmatrix},
\end{align*} 
where $c_{1}, c_{2},\dots, c_{m}$ are non-zero constants in $\mathbb{R}$.
\begin{defn}
Consider the  star transform $\Sc \vf$  of a symmetric 2-tensor field $\vf$ with branches along directions $\vgamma_{1},\vgamma_{2},...,\vgamma_{m}$. We call
$$ \Zc_{1}= \cup_{i=1}^{m}\{ \vxi : \vxi \cdot {\vgamma}_{i} =0\} $$ the set of singular directions of type 1 for $\Sc$.
\end{defn}
Now, let us define three vectors in $\Rb^3$, which will be important for further calculations. For $\vxi \in \mathbb{S}^{1}\setminus \Zc_{1}$, we define 
\begin{align}
\vgamma(\vxi) = - \sum_{i=1}^{m} \frac{c_{i}{\vgamma}_{i}^{2}}{\vxi \cdot {\vgamma}_{i}},\qquad
\vgamma^\dag(\vxi) = - \sum_{i=1}^{m} \frac{c_{i}{\vgamma}_{i}\odot \vgamma_i^\perp}{\vxi \cdot {\vgamma}_{i}},\qquad
\vgamma^\perp(\vxi) = - \sum_{i=1}^{m} \frac{c_{i}({\vgamma_i}^\perp)^2 }{\vxi \cdot {\vgamma}_{i}}\label{eq: gamma_xi perp}.
\end{align}
\begin{defn}
We call $$\Zc_{2} = \left\{ \vxi :  \vgamma^\dag(\vxi) = 0\right\}\ \bigcup \ \left\{ \vxi :\sum_{i=1}^m\frac{c_i}{\vxi \cdot \vgamma_i} = 0 \right\}$$ 
the set of singular directions of type 2 for $\Sc$.
\end{defn}
\begin{thr}\label{th:inversion of star transform S}
 Let $\vf\in C_c^2\left(S^2;D_1\right)$, and $\vgamma_{1},\vgamma_{2}, \dots ,\vgamma_{m}$ be the branch directions of the star transform. Let  \begin{align*}
       \Qc(\vxi) =  \begin{bmatrix}
             {\vgamma}(\vxi)\\
 {\vgamma}^{\dag}(\vxi) \\
 {\vgamma}^{\perp}(\vxi)
        \end{bmatrix}.
    \end{align*}
    Then for any $\vxi \in \Sb^{1}\setminus(\Zc_{1}\cup\Zc_{2}) $ and any $ s\in \Rb$ we have 
\begin{equation} \label{eq: star inversion}
\left[\Qc(\vxi)\right]^{-1}\frac{d}{ds} \Rc ({\Sc}\vf) (\vxi,s)=  \Rc \vf (\vxi,s)
    \end{equation}
where $\Rc \vf$ is the component-wise Radon transform of $\vf$.
\end{thr}

To reconstruct numerically the tensor from its star transform, we consider a star with the branches along $\vgamma_i=(\cos\phi_i,\sin\phi_i)$, where $\phi_1=0, \phi_2=2\pi/3,\phi_3=4\pi/3$, and weights $c_i=1$ for $i=1,2,3$. It can be easily verified that in this case $\Zc_{2} =\varnothing$ and $\Zc_{1}$ consists of six directions. Therefore, formula \eqref{eq: star inversion} will recover $\Rc \vf (\vxi,s)$ for all but six $\xi\in\mathbb{S}^1$, allowing the recovery of $\vf$ using any standard inversion technique for the (classical) Radon transform.\\

\begin{figure}[H]
     \centering
     \begin{subfigure}[b]{0.55\textwidth}
         \centering
         \includegraphics[width=\textwidth]{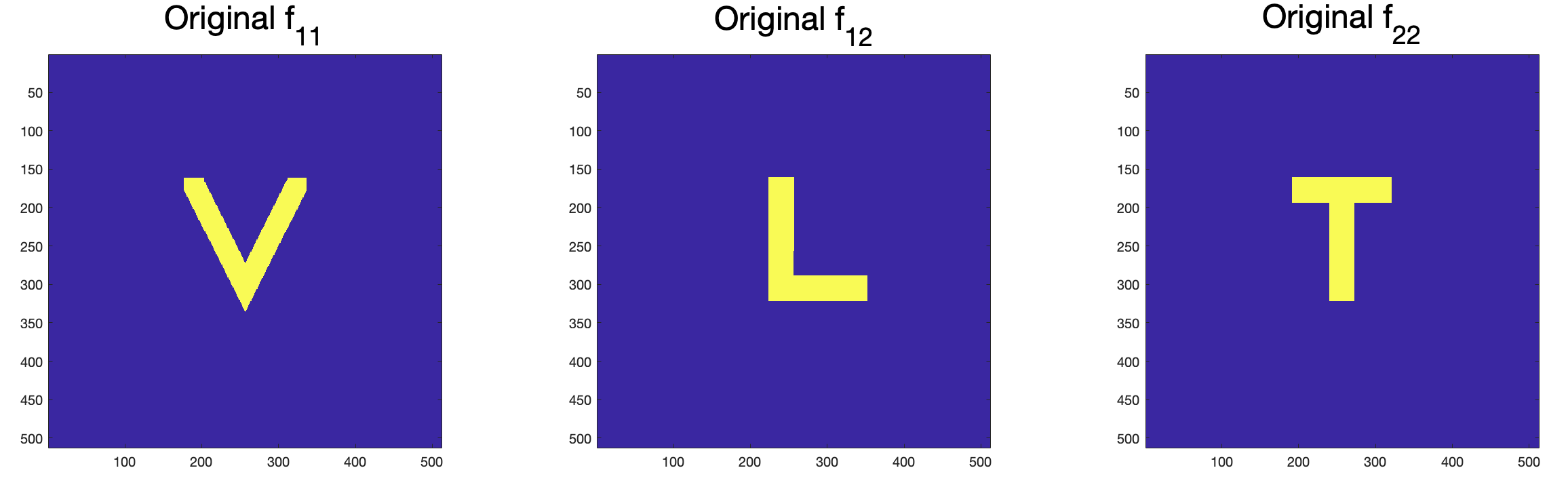}
     \end{subfigure}
     \vfill
     \begin{subfigure}[b]{0.55\textwidth}
         \centering
         \includegraphics[width=\textwidth]{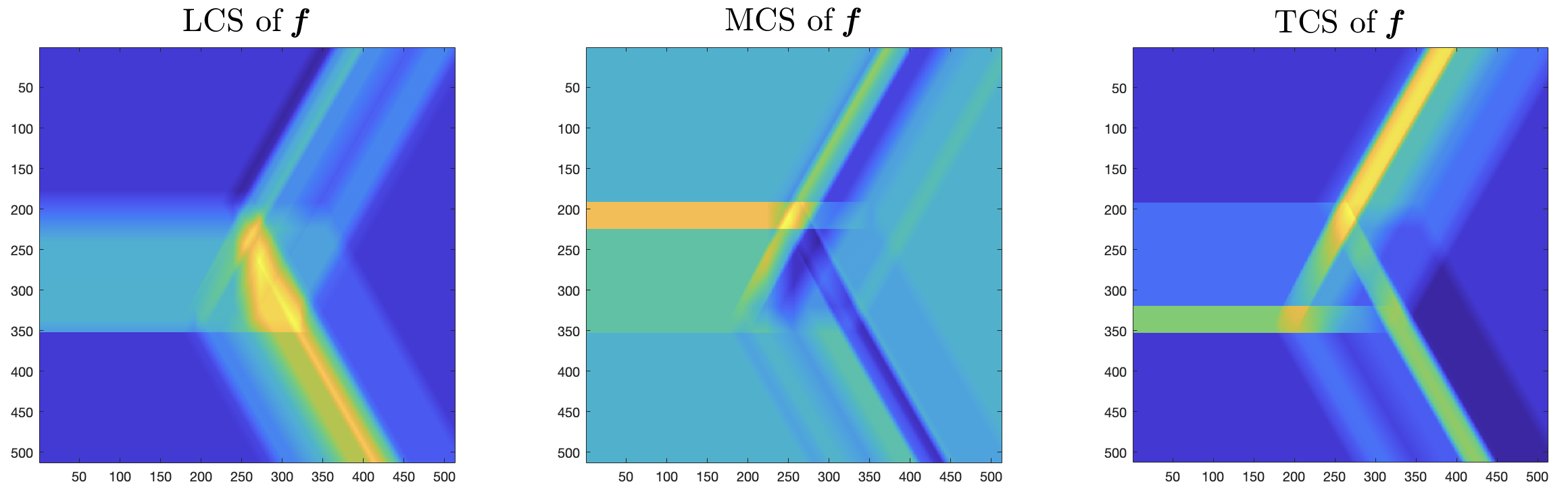}
     \end{subfigure}
     \vfill
     \begin{subfigure}[b]{0.55\textwidth}
         \centering
         \includegraphics[width=\textwidth]{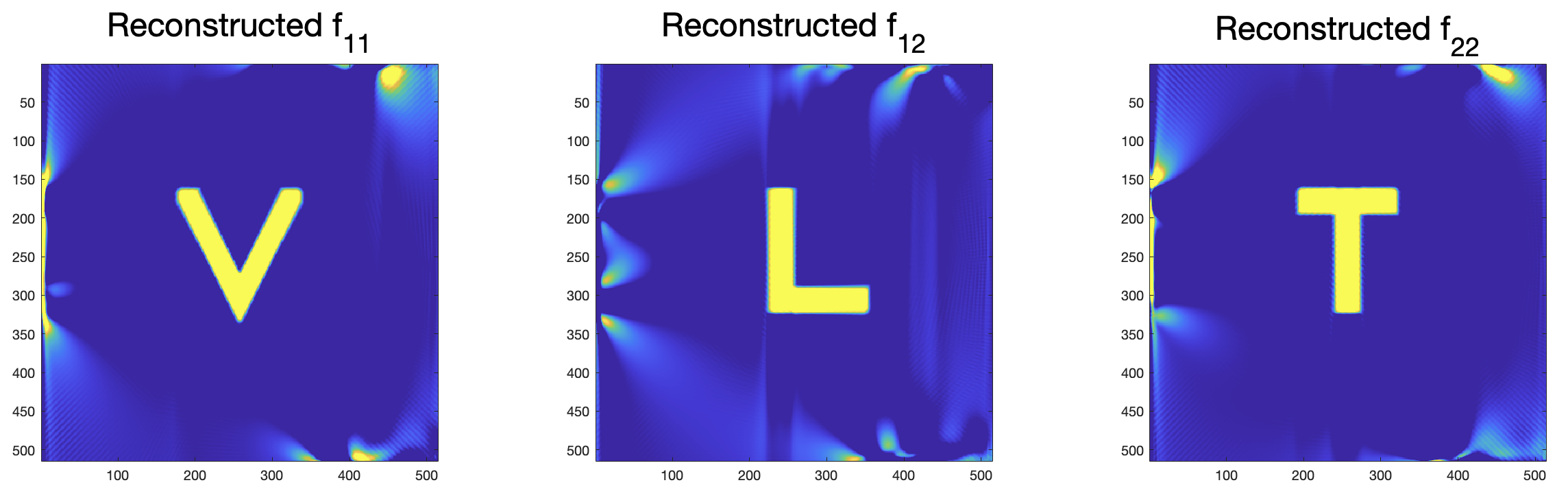}
     \end{subfigure}
        \caption{ Row 1: components of $\vf$; Row 2:  longitudinal (LCS), mixed(MCS), and transverse (TCS) components of $\Sc\vf$; Row 3: reconstructed components of $\vf$ using formula \eqref{eq: star inversion}.}\label{fig:ph2 from star transform}
\end{figure}

\begin{figure}[H]
     \centering
   \includegraphics[width= 0.92\textwidth]{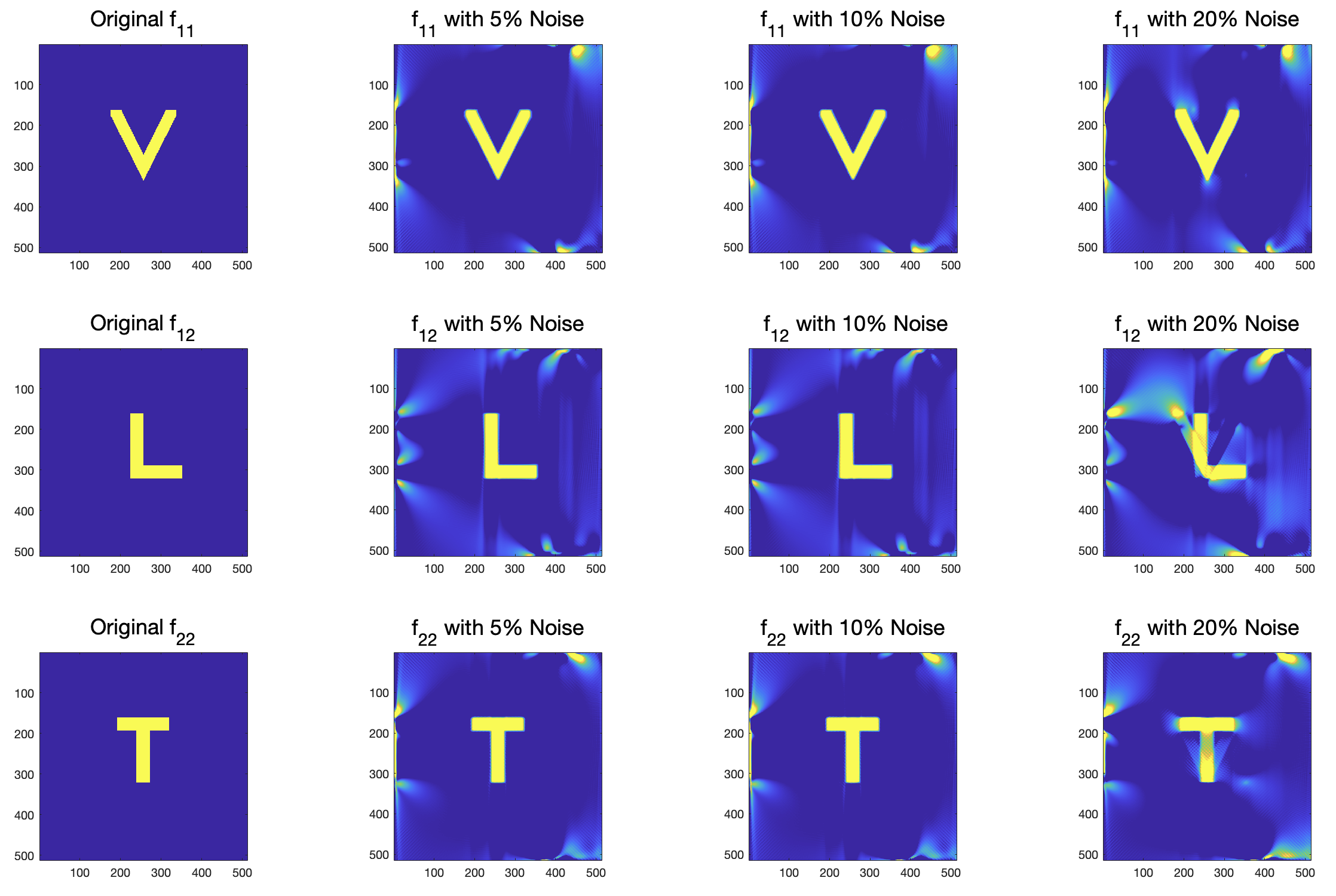}
    \caption{Reconstruction of $\vf$ from $\Sc\vf$ with $5\%, 10\%$ and $20\%$ noise }\label{fig:ph2 from star transform with noise}
 \end{figure}

The numerical reconstruction is performed using the following steps:
\begin{itemize}
    \item   We calculate  $ \Xc_{\vgamma_i}(\vf\cdot {{\vgamma}_{i}}^{2}),
 \Xc_{\vgamma_i}(\vf\cdot {{\vgamma}_{i} \odot {\vgamma}_{i}^{\perp} }),
 \Xc_{\vgamma_i}(\vf\cdot ({\vgamma}_{i}^{\perp})^{2})$. Then, combining these according to the definition, we compute the (forward) data corresponding to $\Sc\vf.$
 \item Using the Matlab function $\textbf{radon}$, we generate the data for $\Rc(\Sc\vf)(\vxi,s)$ and then calculate $\frac{d}{ds} \Rc ({\Sc}\vf) (\vxi,s).$
 \item Next, for each value of $\vxi$ we compute $[\Qc(\vxi)]^{-1}$ and multiply it with $\frac{d}{ds} \Rc ({\Sc}\vf) (\vxi,s).$
 \item Finally, using \textbf{iradon} to  $\displaystyle\left[\Qc(\vxi)\right]^{-1}\frac{d}{ds} \Rc ({\Sc}\vf) (\vxi,s)$, we recover $\vf.$
\end{itemize}
In Figure \ref{fig:ph2 from star transform}, we show the reconstruction using the algorithm discussed above for the non-smooth phantom. Figure \ref{fig:ph2 from star transform with noise} demonstrates the results of the reconstruction in the presence of noise. The relative errors of the reconstructions are summarized in Table \ref{tab:star}.\\

The artifacts present in the reconstructions using the star transform are due to the cutting of the data on the edges of the image domain. Notice, that almost all of those artifacts are located outside of the field of view, i.e. the a priori known disc of support of the phantoms. Therefore, one can eliminate those artifacts through image post-processing, by assigning the value of zero to all pixels outside of that disc. In that case, the relative errors presented in Table \ref{tab:star} will also substantially decrease. 

\begin{table}[h!]
\begin{center}
\begin{tabular}{ |c|c|c|c|c| } 
 \hline
  $\vf$  & No noise & $5\%$ Noise & $10\%$ Noise& $20\%$ Noise \\ 
 \hline
  $f_{11}$  & 114.01\% & 114.72\% & 128.88\% & 175.17\%\\ 
 \hline
 $f_{12}$ & 92.41\% & 96.36\% & 96.81\% &  115.22\%\\  
 \hline
 $f_{22}$ & 95.56\% & 96.02\% & 109.25\% & 139.61\%\\  
 \hline
\end{tabular}
\caption{Relative errors of the reconstruction of $\vf$ from $\Sc\vf$.} \label{tab:star}
\end{center}
\end{table}

\section{Acknowledgements}\label{sec:acknowledge}
GA was partially supported by the NIH grant U01-EB029826. RM was partially supported by SERB SRG grant No. SRG/2022/000947. IZ was supported by the Prime Minister's Research Fellowship from the Government of India.

\bibliography{references}

\begin{thebibliography}{10}

\bibitem{amb_2012}
Gaik Ambartsoumian.
\newblock Inversion of the {V}-line {R}adon transform in a disc and its
  applications in imaging.
\newblock {\em Computers \& Mathematics with Applications}, 64(3):260--265,
  2012.

\bibitem{amb-chapter}
Gaik Ambartsoumian.
\newblock {V}-line and conical {R}adon transforms with applications in imaging.
\newblock In Ronny Ramlau and Otmar Scherzer, editors, {\em The {R}adon
  Transform: The First 100 Years and Beyond}, pages 143--168. De Gruyter,
  Berlin, Boston, 2019.

\bibitem{amb_2023_book}
Gaik Ambartsoumian.
\newblock {\em Generalized Radon Transforms and Imaging by Scattered Particles:
  Broken Rays, Cones, and Stars in Tomography}.
\newblock World Scientific, 2023.

\bibitem{amb-lat_2019}
Gaik Ambartsoumian and Mohammad~J. Latifi.
\newblock The {V}-line transform with some generalizations and cone
  differentiation.
\newblock {\em Inverse Problems}, 35(3):034003, 2019.

\bibitem{Amb_Lat_star}
Gaik Ambartsoumian and Mohammad~J. Latifi.
\newblock Inversion and symmetries of the star transform.
\newblock {\em The Journal of Geometric Analysis}, 31(11):11270--11291, 2021.

\bibitem{Gaik_Mohammad_Rohit_2020}
Gaik Ambartsoumian, Mohammad~J. Latifi, and Rohit~K. Mishra.
\newblock Generalized {V}-line transforms in {2D} vector tomography.
\newblock {\em Inverse Problems}, 36(10):104002, 2020.

\bibitem{Gaik_Mohammad_Rohit_2024_numerics}
Gaik Ambartsoumian, Mohammad~J. Latifi, and Rohit~K. Mishra.
\newblock Numerical implementation of generalized {V}-line transforms on {2D}
  vector fields and their inversions.
\newblock {\em SIAM Journal on Imaging Sciences}, 17(1):595--631, 2024.

\bibitem{Gaik_Indrani_Rohit_24}
Gaik Ambartsoumian, Rohit~K. Mishra, and Indrani Zamindar.
\newblock V-line 2-tensor tomography in the plane.
\newblock {\em Inverse Problems}, 40(3):035003, 2024.

\bibitem{Gaik_Moon}
Gaik Ambartsoumian and Sunghwan Moon.
\newblock {A series formula for inversion of the V-line Radon transform in a
  disc}.
\newblock {\em Computers \& Mathematics with Applications}, 66(9):1567--1572,
  2013.

\bibitem{ambartsoumian2007thermoacoustic}
Gaik Ambartsoumian and Sarah~K. Patch.
\newblock Thermoacoustic tomography: numerical results.
\newblock In {\em Photons Plus Ultrasound: Imaging and Sensing 2007: The Eighth
  Conference on Biomedical Thermoacoustics, Optoacoustics, and Acousto-optics},
  volume 6437, pages 346--355. SPIE, 2007.

\bibitem{Amb_Roy}
Gaik Ambartsoumian and Souvik Roy.
\newblock Numerical inversion of a broken ray transform arising in single
  scattering optical tomography.
\newblock {\em IEEE Transactions on Computational Imaging}, 2(2):166--173,
  2016.

\bibitem{Basko_et_al-97}
Roman Basko, Gengsheng~L. Zeng, and Grant~T. Gullberg.
\newblock Analytical reconstruction formula for one-dimensional {C}ompton
  camera.
\newblock {\em IEEE Transactions on Nuclear Science}, 44(3):1342--1346, Jun
  1997.

\bibitem{Basko_V-proj}
Roman Basko, Gengsheng~L. Zeng, and Grant~T. Gullberg.
\newblock Fully three dimensional image reconstruction from ``{V}''-projections
  acquired by {C}ompton camera with three vertex electronic collimation.
\newblock In {\em 1997 IEEE Nuclear Science Symposium Conference Record},
  volume~2, pages 1077--1081, Nov 1997.

\bibitem{Rohit_Rahul_Manmohan}
Rahul Bhardwaj, Rohit~K. Mishra, and Manmohan Vashisth.
\newblock Inversion of generalized {V}-line transforms of vector fields in
  $\mathbb{R}^2$.
\newblock {\em Preprint, arXiv:2404.12479}, 2024.

\bibitem{derevtsov3}
Evgeny~Yu. Derevtsov and Ivan~E. Svetov.
\newblock Tomography of tensor fields in the plane.
\newblock {\em Eurasian Journal of Mathematical and Computer Applications},
  3(2):24--68, 2015.

\bibitem{durran2013numerical}
Dale~R. Durran.
\newblock {\em Numerical Methods For Wave Equations In Geophysical Fluid
  Dynamics}, volume~32.
\newblock Springer Science \& Business Media, 2013.

\bibitem{FMS-PhysRev-10}
Lucia Florescu, Vadim~A. Markel, and John~C. Schotland.
\newblock Single-scattering optical tomography: Simultaneous reconstruction of
  scattering and absorption.
\newblock {\em Physical Review E}, 81:016602, Jan 2010.

\bibitem{Florescu-Markel-Schotland_2011}
Lucia Florescu, Vadim~A. Markel, and John~C. Schotland.
\newblock Inversion formulas for the broken-ray {R}adon transform.
\newblock {\em Inverse Problems}, 27(2):025002, 2011.

\bibitem{florescu2018}
Lucia Florescu, Vadim~A Markel, and John~C Schotland.
\newblock Nonreciprocal broken ray transforms with applications to fluorescence
  imaging.
\newblock {\em Inverse Problems}, 34(9):094002, 2018.

\bibitem{FMS-09}
Lucia Florescu, John~C. Schotland, and Vadim~A. Markel.
\newblock Single-scattering optical tomography.
\newblock {\em Physical Review E}, 79:036607, Mar 2009.

\bibitem{gouia2014analytical}
Rim Gouia-Zarrad.
\newblock Analytical reconstruction formula for $n$-dimensional conical {R}adon
  transform.
\newblock {\em Computers \& Mathematics with Applications}, 68(9):1016--1023,
  2014.

\bibitem{gouia2014exact}
Rim Gouia-Zarrad and Gaik Ambartsoumian.
\newblock Exact inversion of the conical {R}adon transform with a fixed opening
  angle.
\newblock {\em Inverse Problems}, 30(4):045007, 2014.

\bibitem{herman2009fundamentals}
Gabor~T Herman.
\newblock {\em Fundamentals of Computerized Tomography: Image Reconstruction
  From Projections}.
\newblock Springer Science \& Business Media, 2009.

\bibitem{ilmavirta2022broken}
Joonas Ilmavirta and Gabriel~P. Paternain.
\newblock Broken ray tensor tomography with one reflecting obstacle.
\newblock {\em Communications in Analysis and Geometry}, 30(6):1269--1300,
  2022.

\bibitem{Ilmavirta_Salo_2016}
Joonas Ilmavirta and Mikko Salo.
\newblock Broken ray transform on a {R}iemann surface with a convex obstacle.
\newblock {\em Communications in Analysis and Geometry}, 24(2):379--408, 2016.

\bibitem{Shubham_Manas}
Shubham~R. Jathar, Manas Kar, and Jesse Railo.
\newblock Broken ray transform for twisted geodesics on surfaces with a
  reflecting obstacle.
\newblock {\em The Journal of Geometric Analysis}, 34(7):212, 2024.

\bibitem{Kats_Krylov-13}
Alexander Katsevich and Roman Krylov.
\newblock Broken ray transform: inversion and a range condition.
\newblock {\em Inverse Problems}, 29(7):075008, 2013.

\bibitem{krylov2015inversion}
Roman Krylov and Alexander Katsevich.
\newblock Inversion of the broken ray transform in the case of energy-dependent
  attenuation.
\newblock {\em Physics in Medicine \& Biology}, 60(11):4313, 2015.

\bibitem{morvidone2010}
Marcela Morvidone, Mai~K. Nguyen, Tuong~T. Truong, and Habib Zaidi.
\newblock On the {V}-line {R}adon transform and its imaging applications.
\newblock {\em International Journal of Biomedical Imaging}, 2010:208179, 2010.

\bibitem{Palamodov2017}
Victor Palamodov.
\newblock Reconstruction from cone integral transforms.
\newblock {\em Inverse Problems}, 33(10):104001, 2017.

\bibitem{Rigaud2013}
Ga{\"{e}}l Rigaud, R{\'{e}}mi R{\'{e}}gnier, Mai~K. Nguyen, and Habib Zaidi.
\newblock Combined modalities of {C}ompton scattering tomography.
\newblock {\em IEEE Transactions on Nuclear Science}, 60(3):1570--1577, 2013.

\bibitem{Sharafutdinov_Book}
Vladimir~A. Sharafutdinov.
\newblock {\em Integral Geometry of Tensor Fields}.
\newblock De Gruyter, Berlin, New York, 1994.

\bibitem{Sherson}
Brian Sherson.
\newblock {\em Some Results in Single-Scattering Tomography}.
\newblock PhD thesis, Oregon State University, 2015.
\newblock {PhD} Advisor: D. Finch.

\bibitem{Truong_2011-Vline}
Tuong~T. Truong and Mai~K. Nguyen.
\newblock On new {V}-line {R}adon transforms in $\mathbb{R}^2$ and their
  inversion.
\newblock {\em Journal of Physics A: Mathematical and Theoretical},
  44(7):075206, jan 2011.

\bibitem{walker2021iterative}
Michael~R Walker and Joseph~A. O’Sullivan.
\newblock Iterative algorithms for joint scatter and attenuation estimation
  from broken ray transform data.
\newblock {\em IEEE Transactions on Computational Imaging}, 7:361--374, 2021.

\bibitem{ZSM-star-14}
Fan Zhao, John~C. Schotland, and Vadim~A. Markel.
\newblock Inversion of the star transform.
\newblock {\em Inverse Problems}, 30(10):105001, 2014.

\end{thebibliography}
\bibliographystyle{plain}

\end{document}